\documentclass[letterpaper, 11pt]{amsart}

\usepackage{amsmath,amsthm,amsfonts,amssymb,amscd}
\usepackage{bbm}
\usepackage{bm}
\usepackage{tikz}
\usepackage{tikz-cd}
\usepackage{appendix}
\usepackage{BOONDOX-calo}
\usepackage{standalone}
\usepackage[hidelinks]{hyperref}
\usetikzlibrary{arrows,chains,matrix,positioning,scopes}
\usepackage[letterpaper, left=2.8cm,right=2.8cm, top=2.8cm, bottom=2.8cm]{geometry}
\usepackage{adjustbox}
\usepackage{color}

\newcommand{\CC}{{\mathbb C}}

\newcommand{\ZZ}{{\mathbb Z}}
\newcommand{\QQ}{{\mathbb Q}}

\newcommand{\inEnd}{{\underline{End}}}
\newcommand{\fu}{{\mathfrak{u}}}

\newcommand{\inHom}{{\underline{Hom}}}

\newcommand\dual{\raise0.9ex\hbox{$\scriptscriptstyle\vee$}}
\newcommand{\sN}{{\mathcal{N}}}
\newcommand{\sL}{{\mathcal{L}}}

\newcommand{\sE}{{\mathcal{E}}}
\newcommand{\sX}{{\mathcal{X}}}
\newcommand{\sY}{{\mathcal{Y}}}
\newcommand{\sZ}{{\mathcal{Z}}}
\newcommand{\bT}{{\mathbf{T}}}

\newcommand{\ffu}{\underline{\fu}}
\newcommand{\dbullet}{{\bullet,\bullet}}
\newcommand{\db}{{\bullet,\bullet}}

\theoremstyle{plain}
\newtheorem{thm}{Theorem} 
\newtheorem{prop}[thm]{Proposition}
\newtheorem{lemma}[thm]{Lemma}  
\newtheorem{cor}[thm]{Corollary}
\numberwithin{thm}{subsection}

\newtheorem{manualtheoreminner}{Theorem}
\newenvironment{thm'}[1]{%
  \renewcommand\themanualtheoreminner{#1}%
  \manualtheoreminner
}{\endmanualtheoreminner}

\theoremstyle{definition}
\newtheorem{defn}[thm]{Definition}
\newtheorem{construction}[thm]{Construction}

\newtheorem{notation}[thm]{Notation}

\theoremstyle{remark}
\newtheorem{rem}[thm]{Remark}
\newtheorem*{example}{Example}

\tikzset{>=stealth}

\makeatletter
\def\@seccntformat#1{%
  \protect\textup{\protect\@secnumfont
    \ifnum\pdfstrcmp{subsection}{#1}=0 \bfseries\fi
    \csname the#1\endcsname
    \protect\@secnumpunct
  }%
}  
\makeatother

\makeatletter
\@namedef{subjclassname@2020}{%
  $2020$ Mathematics Subject Classification}
\makeatother

\begin{document}

\title[Blended extensions and motives with maximal unipotent radicals]{On blended extensions in filtered abelian categories and motives with maximal unipotent radicals}
\author{Payman Eskandari}
\address{Department of Mathematics and Statistics, University of Winnipeg, Winnipeg MB, Canada }
\email{p.eskandari@uwinnipeg.ca}
\subjclass[2020]{14F42 , 18M25 (primary); 11F67, 11M32, 14G10  (secondary).}
\begin{abstract}
Grothendieck's theory of blended extensions ({\it extensions panach\'{e}es}) gives a natural framework to study 3-step filtrations in abelian categories. We give a generalization of this theory that is suitable for filtrations with an arbitrary finite number of steps. We use this generalization to study two natural classification problems for objects with a fixed associated graded in an abelian category equipped with a filtration similar to the weight filtration on rational mixed Hodge structures. We then give an application to the study of mixed motives with a given associated graded and maximal unipotent radicals of motivic Galois groups. We prove a homological classification result for the isomorphism classes of such motives when the given associated graded is ``graded-independent", a condition defined in the paper. The special case of this result for motives with 3 weights was proved with K. Murty in \cite{EM2} under some extra hypotheses.

\vspace*{-.2in}
\end{abstract}
\maketitle
\setcounter{tocdepth}{2}
\hypersetup{bookmarksdepth = 2}
\tableofcontents

\section{Introduction}\label{sec: intro}

Let $\bT$ be an abelian category equipped with a filtration $W_\bullet$ similar to the weight filtration on the category of rational mixed Hodge structures or any reasonable category of mixed motives. That is, $W_\bullet$ is indexed by $\ZZ$, functorial, increasing, exact, and finite on every object. Consider a graded object
\[
A = \bigoplus\limits_{r=1}^k A_r,
\]
where $k\geq 2$ and the $A_r$ are nonzero pure objects in an increasing order of weights. One may consider the following two classification problems:
\medskip\par 
\begin{itemize}
\item[(1)] Classify the equivalence classes of all pairs $(X,\phi)$ of an object $X$ of $\bT$ whose associated graded
\[
Gr^WX = \bigoplus\limits_{n} W_nX/W_{n-1}X
\]
is isomorphic to $A$, and an isomorphism $\phi: Gr^WX \rightarrow A$. Here, the equivalence relation for such pairs is defined as follows: two pairs $(X,\phi)$ and $(X',\phi')$ are considered equivalent if there exists an isomorphism $f: X\rightarrow X'$ such that $\phi'\circ Gr^Wf = \phi$. Denote the set of equivalence classes of such pairs by $S'(A)$.

\item[(2)] Classify, up to isomorphism in $\bT$, all objects $X$ whose associated graded $Gr^WX$ is isomorphic to $A$. Note that here the data of a choice of isomorphism $Gr^W X \rightarrow A$ is not recorded at all. Denote the set of isomorphism classes of such $X$ by $S(A)$.
\end{itemize}
\medskip\par 
When $k=2$ and $A=A_1\oplus A_2$, the homological concept that classifies the pairs $(X,\phi)$ up to equivalence of such pairs is the Ext group $Ext^1(A_2,A_1)$. As for the isomorphism classes of $X$, the answer is given by the quotient of $Ext^1(A_2,A_1)$ by the group 
\[Aut(A_1) \times Aut(A_2),\] 
where the actions of automorphisms of $A_1$ and $A_2$ on extension classes is by pushforward and pullback.
\medskip\par 
When $k=3$, the homological concept related to problems (1) and (2) is the concept of a blended extension\footnote{The original French term for the concept is {\it extension panach\'{e}e}. The English translation to the term {\it blended extension}, which we first found in \cite{Ber98}, is attributed by Bertrand to L. Breen.}, introduced by Grothendieck in SGA 7.I \cite[\S 9.3]{Gr68} to study 3-step filtrations. A nice detailed discussion of the relation between problem (1) in the case $k=3$ and blended extensions can be found in the appendix of \cite{RSZ}.
\medskip\par 
The set $S'(A)$ of problem (1) for an arbitrary $k$ has been studied previously in the setting of difference modules\footnote{There is no interesting categorical weight filtration in this setting. The set $S'(A)$ ($\mathcal{F}(A_1,\ldots,A_k)$ is the notation of \cite{RSZ}) is the set of equivalence classes of pairs $(X,\phi)$ of a finitely filtered difference module $X$ and a graded isomorphism $\phi: Gr(X)\rightarrow A$, where the equivalence relation is as discussed above.} over difference rings by Ramis, Salouy and Zhang \cite{RSZ} and in the setting of real mixed Hodge structures in Ferrario's PhD thesis\footnote{I thank Peter Jossen for bringing this to my attention. Note that Ferrario calls a pair $(X,Gr^WX \xrightarrow{\phi, \simeq} A)$ an amalgam of $A$, see Definition 3.2.4 of \cite{Fe20}. Our $S'(A)$ is the same as $\text{Am}(A)$ in his notation.} \cite{Fe20}. In the former setting, Ramis et al  show that $S'(A)$ is a scheme. In the latter setting, Ferrario studies $S'(A)$ to obtain some interesting results on a complex analogue of Grothendieck's section conjecture. In both works, the authors study $S'(A)$ by an inductive approach on the number of graded components $k$.
\medskip\par 
The first and main goal of the present paper is to give a new approach towards problems (1) and (2) in the general setting of an arbitrary $k$ and $\bT$. We introduce a new concept, which we call a {\it generalized extension of a given level $\ell<k$} (Definition \ref{def: gen ext level l}), that is a natural generalization of the notion of a blended extension and provides a natural homological framework for the study of problems (1) and (2). This concept naturally leads to an inductive approach (namely, induction on the level) towards the problems that has better naturalness properties than the more obvious approach of induction on the number of weights. All of this is done in detail in \S \ref{sec: objects in filtered tan cats with given gr}. The main result is Theorem \ref{thm: main thm 1}. A comparison of some features of our approach with the approach of induction on the number of weights $k$ can be found in Remark \ref{rem: comparison with induction on k}.
\medskip\par 
We should point out that we find the classification problem (2) of $S(A)$ more interesting, even though one sometimes prefers to work with $S'(A)$ because of its better moduli properties, as the work \cite{RSZ} in the setting of difference modules illustrates. Note that in passing from $S'(A)$ to $S(A)$ we do not content ourselves with a cursory description of $S(A)$ as the quotient of $S'(A)$ by the obvious action of
\[Aut(A)= \prod_r Aut(A_r).\] 
In fact, one of the advantages of the approach proposed in this paper is that it allows us to trace this action more explicitly, resulting in a particularly simplified outcome in the important  ``totally nonsplit case", discussed in \S \ref{sec: tot nonsplit case}. This case is relevant for the application to motives given in the last part of the paper, which we now proceed to discuss.
\medskip\par 
The second goal of the paper is to consider the problem of classifying motives $X$ with associated graded isomorphic to a given semisimple motive $A$ and with maximal unipotent radicals of motivic Galois groups. This is the subject of \S \ref{sec: mixed motives with maximal unipotent radicals}. By maximality of the unipotent radical here we mean that the Lie algebra of the unipotent radical of the motivic Galois group of $X$ is equal to $W_{-1}\inHom(X,X)$, where $\inHom$ is the internal Hom (that is, the unipotent radical is as large as it can be under the constraints imposed by the weight filtration). The interest in motives with maximal unipotent radicals (of motivic Galois groups) is inspired by Grothendieck's period conjecture, which predicts that the transcendence degree of the field generated by the periods of any motive over $\overline{\QQ}$ should be equal to the dimension of the motivic Galois group of the motive (see for instance, Andr\'{e}'s letter to Bertolin published at the end of \cite{Be20}). In view of this conjecture, among motives over $\overline{\QQ}$ with the same associated graded, the field generated by the periods of a motive with a maximal unipotent radical (if such a motive exists) should have the largest transcendence degree.
\medskip\par 
When $A$ satisfies a certain property (what we call {\it graded-independence}, see Definition \ref{defn: graded independence} and the comments after), we give a particularly simple homological answer to the classification problem of isomorphism classes of motives with associated graded isomorphic to $A$ and maximal unipotent radicals (Theorem \ref{thm: classification of motives with give gr in graded independent case}). The use of the term motive here is merely because of its suggestiveness: in fact, all of \S \ref{sec: mixed motives with maximal unipotent radicals} with the sole exception of \S \ref{sec: examples} takes place in the generality of a filtered tannakian category over a field of characteristic 0 where all pure objects are semisimple. The special case of Theorem \ref{thm: classification of motives with give gr in graded independent case} when $A=A_1\oplus A_2\oplus A_3$ with the $A_r$ pure in an increasing order of weights, $A_3=\mathbbm{1}$, and $Ext^1(\mathbbm{1},A_1)=0$ was proved with K. Murty in \cite{EM2}.

\subsection{Conventions}\label{sec: conventions}
Throughout, by an abelian category we mean one where for every objects $X,Y$ and every integer $i$, Yoneda's $Ext^i(X,Y)$ is a set (rather than just a proper class). By a tannakian category we mean a neutral tannakian category over a field in the sense of \cite{DM82}. By a weight filtration on an abelian category $\bT$ we mean a $\ZZ$-indexed, increasing, functorial, exact filtration $W_\bullet$ that is finite on every object. More explicitly, this means that for each integer $n$, we have an exact (in particular, additive) functor $W_n: \bT\rightarrow \bT$, such that for each object $X$ of $\bT$, $(W_\bullet X)$ is an increasing finite filtration of $X$, i.e.
\begin{align*}
W_{n-1}X \ &\subset \ W_nX \hspace{.3in}(\forall n)\\
W_nX \ &= \ 0 \hspace{.3in}(\forall n\ll 0)\\
W_nX \ &= \ X \hspace{.3in}(\forall n\gg 0),
\end{align*}
and the inclusions $W_{n-1}X\subset W_n X$ for various $X$ give a morphism of functors from $W_{n-1}$ to $W_{n}$ for each $n$. In the context of tannakian categories, a weight filtration is also assumed to be linear and respectful of the tensor structure, i.e. for every $n\in\ZZ$ and objects $X$ and $Y$ of the category,
\[
W_n(X\otimes Y) = \sum\limits_{p+q=n}W_p(X)\otimes W_q(Y).
\]
By a filtered abelian (resp. tannakian) category we mean an abelian (resp. tannakian) category equipped with a weight filtration. An object $X$ of such a category will be called pure if there is an integer $n$ such that $W_{n-1}X=0$ and $W_nX=X$. By a weight of an object $X$ we mean an integer $n$ such that $W_{n-1}X\neq W_nX$. An object $X$ is nonzero and pure if and only if it has exactly one weight.
\medskip\par 
For the purposes of this paper it is important to carefully distinguish between some terms that are sometimes abused in the literature. By an extension we mean a 1-extension, i.e. a short exact sequence. We will have to distinguish between an extension and the object that sits in its middle. We will often use upper case English letters in script font (e.g. $\sN$, $\sX$) for an extension or its class in the $Ext^1$ group (which of the two is the intended interpretation will be clear from the context or explicitly mentioned), and use upper case English letters in ordinary font (e.g. $N$, $X$) for objects of $\bT$. The notation $EXT^i(X,Y)$ (or simply $EXT(X,Y)$ when $i=1$) will be used for the collection of $i$-extensions (or simply, extensions when $i=1$) of $X$ by $Y$. As usual, $Ext^i(X,Y)$ denotes the group (or vector space, when our category is tannakian) of the equivalence classes of $i$-extensions with respect to the standard equivalence given by isomorphism of $i$-extensions, i.e. commuting morphisms between the middle objects that induce identity on $X$ and $Y$. Note that throughout, our notations for Ext and Hom groups as well as for the proper classes EXT do not include a mention of the category under consideration (simply denoting these by $Ext^i$, $Hom$, etc). In a few occasions where the category is not clear from the context, it will be included in the notation as a subscript.
\medskip\par 
Finally, internal Homs are denoted by $\inHom$ and all group actions are left actions.
\medskip\par 
\noindent {\bf Acknowledgements.} I would like to thank Kumar Murty for many helpful conversations. This work is a natural continuation of the joint work \cite{EM2} with him. I would also like to thank Peter Jossen and Daniel Bertrand for some helpful correspondences, and in particular for bringing possible connections to \cite{Fe20} and \cite{BVK16} to my attention. Finally, I would like to thank the anonymous referee for a careful reading of the paper and for providing many valuable suggestions and comments, which helped improve the organization and exposition of the paper.

\section{Recollections on blended extensions}\label{sec: background on blended extensions}
\numberwithin{thm}{section}
In this section we recall some of the basics of the theory of blended extensions. The original reference is \S 9.3 of Grothendieck's \cite{Gr68}. Another excellent reference is Bertrand's \cite{Ber13}. 

Let $\bT$ be an abelian category. Let $A_1$, $A_2$, and $A_3$ be objects of $\bT$. Fix two extensions
\[
\begin{tikzcd}[row sep = small]
\sL: \hspace{0.3in}0 \arrow[r] & A_1 \arrow[r] & L \arrow[r] & A_2 \arrow[r] & 0\\
\sN: \hspace{0.3in}0 \arrow[r] & A_2 \arrow[r] & N \arrow[r] & A_3 \arrow[r] & 0
\end{tikzcd}
\]
in $\bT$.  A blended extension of $\sN$ by $\sL$ is by definition a diagram of the form 
\begin{equation}\label{eq53}
\begin{tikzcd}
   & & 0 \arrow{d} & 0 \arrow{d} &\\
   0 \arrow[r] & A_1 \ar[equal]{d} \arrow[r, ] & L  \arrow[d] \arrow[r, ] &  A_2 \arrow{d} \arrow[r] & 0 \\
   0 \arrow[r] & A_1 \arrow[r] & X \arrow[d] \arrow[r] &  N \arrow{d}  \arrow[r] & 0 \\
   & & A_3 \arrow{d} \ar[equal]{r} & A_3 \arrow{d} & \\
   & & 0 & 0 &   
\end{tikzcd}
\end{equation}
with exact rows and columns. The collection of all blended extensions of $\sN$ by $\sL$ is denoted by  $EXTPAN(\sN,\sL)$ (for {\it extension panach\'{e}es}). We will refer to the object $X$ in the diagram as the middle object.
\medskip\par 
The standard notion of a morphism of blended extensions of $\sN$ by $\sL$ is a morphism in $\bT$ between the middle objects which induces identity on $L$ and $N$ (and hence on $A_1$, $A_2$, and $A_3$). Via this notion of morphisms, $EXTPAN(\sN,\sL)$ is a category in which every morphism is an isomorphism (i.e. is a groupoid category). The collection of isomorphism classes of blended extensions of $\sN$ by $\sL$ is denoted by $Extpan(\sN,\sL)$.
\medskip\par 
We recall three basic results about blended extensions here, which together form the contents of Proposition 9.3.8 of \cite{Gr68}. The first is that when $Extpan(\sN,\sL)$ is nonempty, it has a natural structure of a torsor over $Ext^1(A_3,A_1)$.\footnote{In fact, one can do this at the level of the categories and make $EXTPAN(\sN,\sL)$ a torsor over $EXT(A_3,A_1)$. See page 105 of \cite{Gr68}.} The action of $Ext^1(A_3,A_1)$ on $Extpan(\sN,\sL)$ can be described as follows. Denote the map $N\twoheadrightarrow A_3$ in $\sN$ by $\omega$. Let $\sX\in EXTPAN(\sN,\sL)$ be the blended extension \eqref{eq53}. Let $\sX^h \in EXT(N,A_1)$ be its second row (here $h$ is for horizontal). Given an element $\sE\in EXT(A_3,A_1)$, consider the Baer sum
\[
\sX^h+\omega^\ast\sE \in EXT(N,A_1).
\]
There is a canonical map from $L$ to the middle object of $\sX^h+\omega^\ast\sE$ and a canonical map from this middle object to $A_3$, and these make $\sX^h+\omega^\ast\sE$ the second row of an element of $EXTPAN(\sN,\sL)$. Denote this element by $\sE\ast\sX$ (we may call it the translation of $\sX$ by $\sE$). One can check that the map
\[
EXT(A_3,A_1) \times EXTPAN(\sN,\sL) \rightarrow EXTPAN(\sN,\sL) \hspace{.5in} (\sE,\sX)\mapsto \sE\ast\sX
\]
descends to a map
\begin{equation}\label{eq64}
Ext^1(A_3,A_1) \times Extpan(\sN,\sL) \rightarrow Extpan(\sN,\sL).
\end{equation}
When $Extpan(\sN,\sL)$ is nonempty, the map above makes it a torsor over $Ext^1(A_3,A_1)$. We also use the symbol $\ast$ for the descended action. 
\medskip\par 
Of course, there is an alternative way to try to define the translation of $\sX$ by $\sE$, namely by pushing an extension of $A_3$ by $A_1$ forward along the injection $A_1\hookrightarrow L$ in $\sL$, and then adding it in $EXT(A_3,L)$ to the first vertical extension in $\sX$ (in fact, this is the original construction given in \cite{Gr68}). Bertrand \cite[Appendix]{Ber13} has checked that the two constructions coincide after passing to the level of equivalence classes \eqref{eq64}.
\medskip\par 
For referencing purposes we record the other two basic results as a lemma below. Before stating the lemma, recall the Yoneda product
\[
EXT(A_2,A_1) \times EXT(A_3,A_2) \rightarrow EXT^2(A_3,A_1)
\]
given by splicing. With $\sL$ and $\sN$ as above, the Yoneda product of $\sL$ and $\sN$ is given by
\[
\begin{tikzcd}
0 \arrow[r] & A_1 \arrow[r] & L \arrow[r] & N \arrow[r] & A_3 \arrow[r] & 0
\end{tikzcd}
\]
where the maps $A_1\rightarrow L$ and $N\rightarrow A_3$ come from $\sL$ and $\sN$, and $L\rightarrow N$ is the composition of the map $L\rightarrow A_2$ of $\sL$ with the map $A_2\rightarrow N$ of $\sN$. (See \cite[\S 3]{Yo60}.)

\begin{lemma}\label{lem: criteria for compatibility of extension pairs}
(a) Given $\sN\in EXT(A_3,A_2)$ and $\sL\in EXT(A_2,A_1)$, there exists a blended extension of $\sN$ by $\sL$ if and only if the Yoneda product of $\sL$ and $\sN$ vanishes in $Ext^2(A_3,A_1)$. 
\medskip\par 
\noindent (b) The automorphism group of a blended extension of $\sN$ by $\sL$ is in a canonical bijection with $Hom(A_3,A_1)$. 
\end{lemma}

For proofs, see Proposition 9.3.8(a,c) of \cite{Gr68} (or Lemma 6.4.1 of \cite{EM2} for part (a)). Note that throughout the paper, we shall only deal with blended extensions for which 
\[ Hom(A_3,A_1) \cong 0, \]
so that they always have a trivial automorphism group.
\medskip\par 
We end this subsection with a remark about various equivalence relations for blended extensions and an observation. Throughout the paper, by the standard equivalence relation on blended extensions of $\sN$ by $\sL$ we mean the one considered above (where two blended extensions are considered equivalent if there exists a morphism between their middle objects that induces identity on the first rows and second columns), and the notation $Extpan(\sN,\sL)$ is always used in reference to this relation. There are however two coarser equivalence relations that one may alternatively consider. In the first alternative relation, two blended extensions of $\sN$ by $\sL$ are considered equivalent if there is a morphism between the middle objects that induces morphisms on the analogous objects of the two diagrams such that the induced morphisms on $A_1$, $A_2$ and $A_3$ (but not necessarily on $L$ and $N$) are identity. In the second alternative relation, which is the coarsest of all three equivalence relations, one declares two blended extensions to be equivalent if there exists a morphism between the middle objects that induces isomorphisms between the analogous objects of the two diagrams.
\medskip\par 
We should note that in the paper, our blended extensions will also always satisfy 
\begin{equation}\label{eq59}
Hom(A_2,A_1)\cong Hom(A_3,A_2) \cong 0.
\end{equation}
In this case, the only automorphism of $L$ (resp. $N$) that induces identity and $A_2$ and $A_1$ (resp. $A_3$) is the identity, so that the standard notion of equivalence on $EXTPAN(\sN,\sL)$ coincides with the one only requiring morphisms to be identity on the $A_i$. Later in the paper, when we introduce the notion of generalized extensions, the equivalence that requires inducing identity on the $A_i$ is denoted by $\sim'$, whereas the one that allows for arbitrary automorphisms of the $A_i$ is denoted by $\sim$.\footnote{To avoid confusion, the reader is warned that in our definition of morphisms of generalized extensions (to be given in \S \ref{sec: gen exts, defn}) we will only require commutativity of diagrams. The morphisms on the $A_i$ will be allowed to be arbitrary (even zero).}
\medskip\par 
Finally, we make a further observation about blended extensions in the case where \eqref{eq59} holds. In general, without assuming \eqref{eq59}, thanks to the fact that isomorphisms in $EXTPAN(\sN,\sL)$ are identity on $\sN$ and $\sL$, by sending the isomorphism class of a blended extension to the class of the second row of a representative we have a well-defined map
\[
-^h: Extpan(\sN,\sL) \rightarrow Ext^1(N,A_1).
\]
As before, let $\omega$ be the map $N\twoheadrightarrow A_3$ in $\sN$. Let $\iota$ be the map $A_2\hookrightarrow N$. Given $\sX\in Extpan(\sN,\sL)$ and $\sE\in Ext^1(A_3,A_1)$, by definition of the torsor structure on $Extpan(\sN,\sL)$ we have
\[
(\sE\ast \sX)^h = \omega^\ast\sE+\sX^h.
\]
Thus $(\sE\ast \sX)^h=\sX^h$ in $Ext^1(N,A_1)$ if and only if 
\[\sE \in \ker\bigm(Ext^1(A_3,A_1)\xrightarrow{ \ \omega^\ast \ } Ext^1(N,A_1)\bigm) \ = \ \text{Image}\bigm(Hom(A_2,A_1)\rightarrow Ext^1(A_3,A_1)\bigm),\]
where the latter map is the connecting homomorphism in the long exact sequence obtained by applying $Hom(-,A_1)$ to $\sN$. Since the action of $Ext^1(A_3,A_1)$ on $Extpan(\sN,\sL)$ is transitive, we obtain the following statement:

\begin{lemma}\label{lem: X mapsto X^h is injective when Hom(A_2,A_1)=0}
Suppose that $Hom(A_2,A_1)\cong 0$. Then the map
\[
Extpan(\sN,\sL) \rightarrow Ext^1(N,A_1)
\]
which sends the class of a blended extension to the class of its second row is injective.
\end{lemma}

\section{Objects with a prescribed associated graded}\label{sec: objects in filtered tan cats with given gr}
\numberwithin{thm}{subsection}
\subsection{Two classification problems}\label{sec: statement of problems for S(A) and S'(A)}
From here until the end of \S \ref{sec: objects in filtered tan cats with given gr} we shall fix an abelian category $\bT$ equipped with a weight filtration $W_\bullet$ (see \S \ref{sec: conventions}). We shall also fix a positive integer $k$, and nonzero pure objects $A_1,\ldots, A_k$ of $\bT$, respectively of weights $p_1,\ldots, p_k$ with $p_1<\cdots<p_k$. Set
\[
A:= \bigoplus\limits_{1\leq r\leq k} A_r.
\]
For now, we assume $k\geq 2$, but the real case of interest is when $k\geq 3$ (in fact, the new content is for $k>3$, which goes beyond the case of classical blended extensions).
\medskip\par 
Our goal is to study the sets $S(A)$ and $S'(A)$ introduced in \S \ref{sec: intro}. Let us recall the definitions.

\begin{defn}\label{def: S(A) and S'(A)}
\noindent (a) We denote by $S(A)$ the set of isomorphism classes of objects of $\bT$ whose associated graded (with respect to the weight filtration) is isomorphic to $A$:
\[
S(A):=~\{X\in \text{Obj}(\bT): \text{$Gr^WX$ is isomorphic to $A$}\}\bigm/ \text{isomorphism in $\bT$}.
\]
Note that here, we do not keep track of the data of the isomorphisms between the associated gradeds and $A$.

\noindent (b) We denote by $S'(A)$ the set of equivalence classes of pairs 
\[
(X, \, Gr^WX\xrightarrow{\phi, \, \simeq }A ) 
\]
of an object $X$ of $\bT$ whose associated graded is isomorphic to $A$ together with a choice of an isomorphism $\phi: Gr^WX\rightarrow A$. Here, two pairs $(X, \phi )$ and $(X', \phi')$ are declared to be equivalent if there exists an isomorphism $f: X\rightarrow X'$ for which we have $\phi' \, Gr^Wf=\phi$, where $Gr^Wf:  Gr^WX\rightarrow Gr^WX'$ is the isomorphism induced by $f$.
\end{defn}

The group
\[ Aut(A)= \prod\limits_r Aut(A_r)\]
acts on $S'(A)$ by twisting the isomorphism between the associated graded and $A$. More precisely, given a pair $(X, \phi )$ as in (b) and $\sigma\in Aut(A)$, we set 
\[
\sigma\cdot (X, \phi ) = (X, \sigma\phi ).
\]
This defines an action of $Aut(A)$ on the collection of pairs $(X,\phi)$ as in (b) which is easily seen to descend to an action on $S'(A)$.

There is an obvious surjection 
\[
S'(A)\twoheadrightarrow S(A)
\]
induced by forgetting the data of $\phi$. The reader can easily see that two elements of $S'(A)$ are mapped to the same element of $S(A)$ if and only if they belong to the same orbit of $Aut(A)$. For future referencing, we record the conclusion:

\begin{lemma}\label{lem: relation between S'(A) and S(A)}
The natural surjection $S'(A)\twoheadrightarrow S(A)$ induced by $(X,\phi)\mapsto X$ descends to a bijection
\[
S'(A)/Aut(A) \cong S(A).
\]
\end{lemma}

To motivate our approach towards the study of $S'(A)$ and $S(A)$ and and put the results in a better context, in \S \ref{sec: k=3 case} we will briefly discuss the picture for the case $k=3$ without proofs. We will then state our main result for an arbitrary $k$ in \S \ref{sec: statement of the main theorem} (Theorem \ref{thm: main thm 1}). Sections \S \ref{sec: gen exts, defn} - \S \ref{sec: fibers 3} are then dedicated to the proof of this result.

\subsection{The case $k=3$}\label{sec: k=3 case}
The claims that will be made in this subsection will be proved later as we prove the general result for an arbitrary $k$. Given an object $X$ of $\bT$ with
\[Gr^WX\simeq A= A_1\oplus A_2\oplus A_3\]
where $A_r$ is nonzero pure of weight $p_r$ and $p_1<p_2<p_3$, set $X_r:=W_{p_r}X$ so that the $X_r$ form a 3-step filtration. Choosing an isomorphism $\phi: Gr^W X \rightarrow A$ we obtain a blended extension 
\begin{equation}\label{eq58}
\begin{tikzcd}
   & & 0 \arrow{d} & 0 \arrow{d} &\\
   0 \arrow[r] & A_1 \ar[equal]{d} \arrow[r] & X_2  \arrow[d] \arrow[r] &  A_2 \arrow{d} \arrow[r] & 0 \\
   0 \arrow[r] & A_1 \arrow[r] & X \arrow[d] \arrow[r] &  X/X_1 \arrow{d}  \arrow[r] & 0 \\
   & & A_3 \arrow{d} \ar[equal]{r} & A_3 \arrow{d} & \\
   & & 0 & 0 &   
\end{tikzcd}
\end{equation}
with obvious arrows. It is easy to see that by sending the equivalence class of $(X,\phi)$ to the classes of the extensions of the top row and the right column we obtain a (well-defined) map
\begin{equation}\label{eq52}
S'(A) \rightarrow Ext^1(A_2,A_1)\times Ext^1(A_3,A_2).
\end{equation}
Let $\sN\in Ext^1(A_3,A_2)$ and $\sL\in Ext^1(A_2,A_1)$. Fix representative extensions for $\sN$ and $\sL$. In view of the fact that
\[Hom(A_2,A_1)\cong Hom(A_3,A_2) \cong 0,\]
one easily sees that the fiber of the map \eqref{eq52} above $(\sL, \sN)$ is in a canonical bijection with $Extpan(\sN,\sL)$. This bijection sends the equivalence class of a pair $(X,\phi)$ above $(\sL,\sN)$ to the class of the blended extension \eqref{eq58}, with the top row and right column replaced by the chosen representatives of $\sL$ and $\sN$ via the canonical isomorphisms. (Note that since $Hom(A_2,A_1)$ and $Hom(A_3,A_2)$ vanish, there is a canonical isomorphism between any two extensions representing $\sL$ or $\sN$.) Thanks to the general theory of blended extensions (see \S \ref{sec: background on blended extensions}), each fiber of \eqref{eq52} is thus either empty or a torsor over $Ext^1(A_3,A_1)$. Moreover, one can see that the torsor structure on the fiber of \eqref{eq52} above $(\sL, \sN)$ (when nonempty) is canonical, in the sense that it does not depend on the choice of representative extensions for $\sN$ and $\sL$. Finally, it follows also from the theory of blended extensions that the image of \eqref{eq52} is the kernel of the Yoneda composition
\[
Ext^1(A_2,A_1)\times Ext^1(A_3,A_2) \rightarrow Ext^2(A_3,A_1).
\]

Let us turn our attention to $S(A)$. The map \eqref{eq52} descends to a map
\begin{equation}\label{eq51}
S(A) \rightarrow \bigm(Ext^1(A_2,A_1)\times Ext^1(A_3,A_2)\bigm)/Aut(A),
\end{equation}
where the action of $Aut(A)$ is by pushforwards and pullbacks of extensions: an element
\[(\sigma_1,\sigma_2, \sigma_3)\in Aut(A_1)\times Aut(A_2)\times Aut(A_3) = Aut(A)\]
sends $(\sL, \sN)$ to $((\sigma_1)_\ast (\sigma_2^{-1})^\ast\sL, (\sigma_2)_\ast (\sigma_3^{-1})^\ast\sN)$. The fiber of \eqref{eq51} above the $Aut(A)$-orbit of $(\sL,\sN)$ is the image in $S(A)$ of the fiber of \eqref{eq52} above $(\sL,\sN)$. One can show that there is a group $\Gamma(\sL,\sN)$ with a natural action on the set $Extpan(\sN,\sL)$ such that the fiber of \eqref{eq51} above the orbit of $(\sL,\sN)$ can be identified with
\[
Extpan(\sN,\sL)/\Gamma(\sL,\sN).
\]

We briefly include the description of $\Gamma(\sL,\sN)$ and its action as it will give some intuition for the analogous group action in the case of an arbitrary $k$. Fixing representatives 
\[
\begin{tikzcd}
  0 \arrow[r] &  A_1 \arrow[r, ] & L \arrow[r, ] &  A_2  \arrow[r] & 0
\end{tikzcd}
\]
and 
\[
\begin{tikzcd}
  0 \arrow[r] &  A_2 \arrow[r, ] & N \arrow[r, ] &  A_3  \arrow[r] & 0
\end{tikzcd}
\]
for $\sL$ and $\sN$, the group $\Gamma(\sL,\sN)$ is the subgroup of $Aut(L)\times Aut(N)$ consisting of pairs $(\sigma_L,\sigma_N)$ such that the automorphisms of $A_2$ induced by $\sigma_L$ and $\sigma_N$ coincide. Then the action of $\Gamma(\sL,\sN)$ on $Extpan(\sN,\sL)$ is by twisting the arrows: the class of the blended extension of the left below is sent to the class of the one on the right (here, $\sigma_{A_r}$ is the automorphism of $A_r$ induced by $\sigma_L$ or $\sigma_N$):
\begin{equation}\label{eq1}
\begin{tikzcd}
   & & 0 \arrow{d} & 0 \arrow{d} &\\
   0 \arrow[r] & A_1 \ar[equal]{d} \arrow[r, ] & L  \arrow[d, "\iota"] \arrow[r, ] &  A_2 \arrow{d} \arrow[r] & 0 \\
   0 \arrow[r] & A_1 \arrow[r, "j"] & X \arrow[d, "\omega"] \arrow[r, "\pi"] &  N \arrow{d}  \arrow[r] & 0 \\
   & & A_3 \arrow{d} \ar[equal]{r} & A_3 \arrow{d} & \\
   & & 0 & 0 &   
\end{tikzcd} 
\hspace{.2in}
\begin{tikzcd}
   & & 0 \arrow{d} & 0 \arrow{d} &\\
   0 \arrow[r] & A_1 \ar[equal]{d} \arrow[r, ] & L  \arrow[d, "\iota\sigma_L^{-1}"] \arrow[r, ] &  A_2 \arrow{d} \arrow[r] & 0 \\
   0 \arrow[r] & A_1 \arrow[r, "j\sigma_{A_1}^{-1}"] & X \arrow[d, "\sigma_{A_3}\omega"] \arrow[r, "\sigma_N\pi"] &  N \arrow{d}  \arrow[r] & 0. \\
   & & A_3 \arrow{d} \ar[equal]{r} & A_3 \arrow{d} & \\
   & & 0 & 0 &   
\end{tikzcd}
\end{equation}
Note that the arrows on the top row and the right column remain unchanged. The stabilizer of the class of the blended extension on the left under this action is the image of the natural injection
\[
Aut(X) \hookrightarrow Aut(L)\times Aut(N).
\]

\subsection{Statement of the main result and the idea of the proof}\label{sec: statement of the main theorem}
We shall generalize the above picture for $k=3$ to the case of an arbitrary number of steps in the weight filtration. The main result is the following:

\begin{thm}\label{thm: main thm 1}
Let $k\geq 2$. There exist sets $S'_\ell(A)$ for $1\leq \ell\leq k-1$ and maps
\begin{equation}\label{eq maps S'_l main thm}
S'_{k-1}(A)\xrightarrow{\Theta_{k-1}} S'_{k-2}(A) \xrightarrow{\Theta_{k-2}} S'_{k-3}(A) \rightarrow \cdots \rightarrow S'_2(A) \xrightarrow{\Theta_{2}} S'_{1}(A)
\end{equation}
with the following properties:
\begin{itemize}
\item[(a)] There are canonical bijections $S'_{k-1}(A)\cong S'(A)$ \ and \ $S'_1(A)\cong \prod\limits_{r} Ext^1(A_{r+1},A_{r})$.

\item[(b)] Let $2\leq \ell\leq k-1$. Every nonempty fiber of the map $\Theta_\ell: S'_{\ell}(A)\rightarrow S'_{\ell-1}(A)$ is canonically a torsor for
\[
\prod\limits_{r}Ext^1(A_{r+\ell},A_r).
\]
\item[(c)] Let $2\leq \ell\leq k-1$. If the $Ext^2$ groups 
\[Ext^2(A_{r+\ell},A_r) \hspace{.2in} (1\leq r\leq k-\ell)\]
vanish, then the map $\Theta_\ell$ is surjective.

\item[(d)] There is a natural action of $Aut(A)$ on each $S'_\ell(A)$ ($1\leq \ell\leq k-1$) such that setting $S_\ell(A)=S'_\ell(A)/Aut(A)$, the maps in \eqref{eq maps S'_l main thm} and part (a) descend to maps
\[S_{k-1}(A)\rightarrow S_{k-2}(A) \rightarrow S_{k-3}(A) \rightarrow \cdots \rightarrow S_2(A)\rightarrow S_{1}(A)\]
and
\[
S_{k-1}(A)\cong S(A) \hspace{.2in}\text{and}\hspace{.2in} S_1(A)\cong \bigm( \prod\limits_{r} Ext^1(A_{r+1},A_{r})\bigm)\bigm/Aut(A)
\]
(the action of $Aut(A)=\prod\limits_{1\leq j\leq k}Aut(A_j)$ on $\prod\limits_{r} Ext^1(A_{r+1},A_{r})$ being by pushforwards and pullbacks; see \S \ref{sec: levels 1 and k-1} for more details).

\item[(e)] Let $2\leq \ell\leq k-1$. Denote the induced map $S_{\ell}(A)\rightarrow S_{\ell-1}(A)$ also by $\Theta_\ell$. For every $\epsilon \in S_{\ell-1}(A)$ and every $\epsilon'\in S'_{\ell-1}(A)$ above $\epsilon$, the fiber $\Theta_\ell^{-1}(\epsilon)$ is the image of the fiber $\Theta_\ell^{-1}(\epsilon')$ under the quotient map $S'_\ell(A)\twoheadrightarrow S_\ell(A)$. Moreover, there exists a group $\Gamma(\epsilon')$ acting on the fiber $\Theta_\ell^{-1}(\epsilon')$ such that the map $\Theta_\ell^{-1}(\epsilon') \twoheadrightarrow \Theta_\ell^{-1}(\epsilon)$ induces a bijection
\[
\Theta_\ell^{-1}(\epsilon')/\Gamma(\epsilon') \cong  \Theta_\ell^{-1}(\epsilon).
\]
\end{itemize}
\end{thm}

The proof of the theorem will be given in \S \ref{sec: gen exts, defn} - \S \ref{sec: fibers 3}. We will also describe the stabilizers of the actions of part (e) (Lemma \ref{lem: stablizer of the action by automorphisms}). In the rest of this subsection we will first give a rough idea of the proofs and some of the constructions involved. We then compare some features of the approach of Theorem \ref{thm: main thm 1} to the inductive approach on the number of weights (Remark \ref{rem: comparison with induction on k}). The subsection ends with a map of the proof of different parts of Theorem \ref{thm: main thm 1}.

Let $X$ be an object of $\bT$ whose associated graded is isomorphic to $A$. Fix an isomorphism $Gr^WX\rightarrow A$ to identify the two. For any integers $m,n$ with $0\leq m<n\leq k$, set
\[
X_{m,n}:= W_{p_n}X/W_{p_m}X,
\] 
where we have set $p_0=p_1-1$ (so that $W_{p_0}X=0$). It is convenient to introduce the following notation: for $m,n$ with $0\leq m< n\leq k$ and $n-m\geq 2$, let $\sX^h_{m,n}$ and $\sX^v_{m,n}$ be the following two extensions with $X_{m,n}$ in the middle:
\[
\sX^h_{m,n}: \hspace{.3in}
\begin{tikzcd}
   0 \arrow[r] &  A_{m+1}\arrow[r, ] & X_{m,n} \arrow[r, ] &  X_{m+1,n}  \arrow[r] & 0
\end{tikzcd}
\]
\[
\sX^v_{m,n}: \hspace{.3in}
\begin{tikzcd}
  0 \arrow[r] &  X_{m, n-1}\arrow[r, ] & X_{m,n} \arrow[r, ] &  A_n  \arrow[r] & 0
\end{tikzcd}
\]
where we have used our fixed isomorphism $Gr^WX\rightarrow A$ to identify each $X_{r-1,r}$ with $A_r$. Here, the superscripts $h$ and $v$ stand for horizontal and vertical, respectively; the reason for the choice of notation is that these will be considered respectively as horizontal and vertical extensions in diagrams of blended extensions.  

The $X_{m,n}$ fit into the commutative diagram
\begin{equation}\label{eq gen ext for M}
\begin{tikzcd}[column sep=small, row sep=small]
A_1         \arrow[d, hookrightarrow] & & & &&& \\
X_{0,2} \arrow[r, twoheadrightarrow] \arrow[d, hookrightarrow] & A_2 \arrow[d, hookrightarrow]     & & & &&\\
X_{0,3} \arrow[r, twoheadrightarrow]  \arrow[d, hookrightarrow] & X_{1,3} \arrow[r, twoheadrightarrow] \arrow[d, hookrightarrow]   & A_3 \arrow[d, hookrightarrow] & & &&\\
X_{0,4} \arrow[r, twoheadrightarrow] \arrow[d, hookrightarrow]  & X_{1,4} \arrow[r, twoheadrightarrow]  \arrow[d, hookrightarrow]    & X_{2,4} \arrow[r, twoheadrightarrow]  \arrow[d, hookrightarrow] & A_4 \arrow[d, hookrightarrow] & &&\\
\vdots& \vdots & \vdots & \vdots &\ddots &&\\
&&&&&&\\
X_{0,k-1} \arrow[d, hookrightarrow] \arrow[r, twoheadrightarrow] & X_{1,k-1} \arrow[d, hookrightarrow] \arrow[r, twoheadrightarrow] &\cdots & ~ & X_{k-3,k-1} \arrow[l, twoheadleftarrow] \arrow[d, hookrightarrow] \arrow[r, twoheadrightarrow] &A_{k-1} \arrow[d,hookrightarrow] &\\
X_{0,k}  \arrow[r, twoheadrightarrow] & X_{1,k} \arrow[r, twoheadrightarrow] &\cdots & ~ & X_{k-3,k} \arrow[l, twoheadleftarrow] \arrow[r, twoheadrightarrow] & X_{k-2,k}  \arrow[r, twoheadrightarrow] &A_k. 
\end{tikzcd}
\end{equation}
Every $X_{m,n}$ appears in the diagram exactly once. Each horizontal arrow is surjective and is given by modding out by the first step in the weight filtration on the domain. The vertical arrows are all injective and are the inclusions $X_{m,n-1}\hookrightarrow X_{m,n}$ given by the weight filtration. 

Roughly speaking, our goal is to obtain all $X$ or $(X, Gr^WX\xrightarrow{\simeq} A)$ up to the appropriate equivalence relation. Our approach is to do this step by step as follows: First consider possibilities for $(X_{r-1,r+1})_r$, i.e. the first diagonal below the $A_r$. The object $X_{r-1,r+1}$ is an extension of $A_{r+1}$ by $A_{r}$. So we must look at 
\[ \prod\limits_r EXT(A_{r+1},A_{r})\]
up to some equivalence. Fixing $(X_{r-1,r+1})_r$, we now consider the possibilities for $(X_{r-1,r+2})_r$, i.e. the second diagonal below the $A_r$. The object $X_{r-1,r+2}$ will be the middle object of the blended extension
\[
\begin{tikzcd}
   & & 0 \arrow{d} & 0 \arrow{d} &\\
   0 \arrow[r] & A_{r} \ar[equal]{d} \arrow[r, ] & X_{r-1,r+1}  \arrow[d] \arrow[r, ] &  X_{r,r+1} \arrow{d} \arrow[r] & 0 \\
   0 \arrow[r] & A_{r} \arrow[r] & X_{r-1,r+2} \arrow[d] \arrow[r] & X_{r,r+2}\arrow{d}  \arrow[r] & 0 \\
   & & A_{r+2} \arrow{d} \ar[equal]{r} & A_{r+2} \arrow{d} & \\
   & & 0 & 0 &   
\end{tikzcd}
\]
of $\sX^v_{r,r+2}$ by $\sX^h_{r-1,r+1}$. We must thus look at
\[
\prod_r EXTPAN(\sX^v_{r,r+2}, \sX^h_{r-1,r+1})
\]
up to some equivalence. We continue in the same fashion until we get to the possibilities for $X_{0,k}=X$. Of course, one also has to keep track of the appropriate equivalence relations in each step.

To make this approach precise, we introduce the notion of a {\it generalized extension of $A$ of a given level $\ell$} ($1\leq \ell\leq k-1$) (Definitions \ref{defn: generalized extensions} and \ref{def: gen ext level l}). A generalized extension of level $\ell$ of $A$ is the abstract data of a diagram as in \eqref{eq gen ext for M}, but only with the first $\ell$ diagonals below the $A_r$ included. Thus a generalized extension of level $k-1$ is an abstract version of the full diagram \eqref{eq gen ext for M}. When $k=2$, a generalized extension of level 1 of $A$ is simply an extension of $A_2$ by $A_1$. When $k=3$, a generalized extension of level 2 of $A$ will simply be the data of a blended extension as in the left diagram of \eqref{eq1}, with varying $L$ and $N$, but $A_1, A_2, A_3$ fixed. For any $k$, the data of a generalized extension of level 1 of $A$ consists of an extension of $A_{r+1}$ by $A_r$ for each $1\leq r\leq k-1$. We highlight that our notion of a generalized extension becomes interesting when the level is less than $k-1$ (as in level $k-1$, everything is determined by the bottom left object).

The sets $S'_\ell(A)$ and $S_\ell(A)$ in Theorem \ref{thm: main thm 1} are the quotients of the collection of all generalized extensions of level $\ell$ of $A$ by suitable equivalence relations.\footnote{This characterization of $S_\ell(A)$ will be equivalent to taking $S_\ell(A)$ to be the quotient of $S'_\ell(A)$ by $Aut(A)$.} A pair $(X,\phi)$ gives rise to a generalized extension of level $k-1$, inducing the identification $S'(A)\cong S_{k-1}(A)$. The maps $S'_\ell(A)\rightarrow S'_{\ell-1}(A)$ and $S_\ell(A)\rightarrow S_{\ell-1}(A)$ are simply induced by truncation.

\begin{rem}\label{rem: comparison with induction on k}
Of course, there is a somewhat more straightforward inductive approach towards the study of $S'(A)$ and $S(A)$, namely to induct on the number of weights of $A$. Setting
\[
A_{\leq \ell} := \bigoplus\limits_{r\leq \ell} A_r,
\]
the weight filtration gives rise to maps
\[
S'(A)=S'(A_{\leq k}) \rightarrow S'(A_{\leq k-1}) \rightarrow \cdots \rightarrow S'(A_{\leq 3}) \rightarrow S'(A_{\leq 2})
\]
and 
\[
S(A)=S(A_{\leq k}) \rightarrow S(A_{\leq k-1}) \rightarrow \cdots \rightarrow S(A_{\leq 3}) \rightarrow S(A_{\leq 2}).
\]
The fiber of $S'(A_{\leq \ell})\rightarrow S'(A_{\leq \ell-1})$ above the equivalence class of $(X,\phi)$ is in a canonical bijection with $Ext^1(A_\ell, X)$ (this is written in detail in \cite{Fe20}, see Proposition 3.2.9 and Remark 3.2.10 therein).

The inductive approach using the level proposed in this paper is not merely more elegant. The application to mixed motives with maximal unipotent radicals given in \S \ref{sec: mixed motives with maximal unipotent radicals} illustrates the usefulness of the approach proposed in this paper. We also draw the reader's attention to the better naturalness properties of the approach of induction on the level compared to induction on the number of weights: Every nonempty fiber of $S'_\ell(A)\rightarrow S'_{\ell-1}(A)$ is canonically a torsor for the same group $\prod_r Ext^1(A_{r+\ell}, A_r)$ (compare with the structure of the fibers of $S'(A_{\leq \ell})\rightarrow S'(A_{\leq \ell-1})$). Also, the fiber of $S(A_{\leq \ell})\rightarrow S(A_{\leq \ell-1})$ above the isomorphism class of $X$ is in bijection with
\[
Ext^1(A_\ell, X) \bigm/ Aut(X)\times Aut(A_\ell),
\]
where the actions of $Aut(X)$ and $Aut(A_\ell)$ are by pushforward and pullback of extensions. Unless $Ext^1(A_\ell, X)$ is trivial, the action of $Aut(X)\times Aut(A_\ell)$ on $Ext^1(A_\ell, X)$ is never trivial. Compare this situation with the one in \S \ref{sec: thm B, tot nonsplit case}.
\end{rem}

\noindent {\it Map of the proof of Theorem \ref{thm: main thm 1}:} The machinery of generalized extensions of $A$ of a given level is developed in \S \ref{sec: gen exts, defn} - \S \ref{sec: basic properties of gen exts}. This is used in \S \ref{sec: equiv rels on gen exts} to formally define the sets $S'_\ell(A)$ and $S_\ell(A)$ and the maps $\Theta_\ell$. The characterization of the sets $S'_\ell(A)$ and $S'_\ell(A)$ in levels $\ell=1$ and $\ell= k-1$ asserted in parts (a) and (d) of the theorem is established in \S \ref{sec: levels 1 and k-1}. The remaining assertions are addressed in the remainder of \S \ref{sec: objects in filtered tan cats with given gr}. Parts (b) and (c) are the subject of \S \ref{sec: fibers 1}-\ref{sec: fibers of truncations II}. Part (e) is the subject of \S \ref{sec: fibers 3}. In the same subsection we will also give a characterization of the stabilizers of the group actions of part (e).

\subsection{Generalized extensions - Definitions}\label{sec: gen exts, defn}
In this subsection we define our notion of generalized extensions of a given level. This notion will be the key to our approach to the classification problems of interest in the paper.

As mentioned in \S \ref{sec: statement of the main theorem}, given any pair $(X, \phi)$ of an object $X$ of $\bT$ and an isomorphism $\phi: Gr^WX\rightarrow A$, setting $X_{m,n}=W_{p_n}X/W_{p_m}X$ for $0\leq m<n\leq k$ with $p_0:=p_1-1$, the natural inclusions and quotient maps among the $X_{m,n}$ give rise to the diagram \eqref{eq gen ext for M}. The following definition, modelled based on this diagram, formalizes the situation.

\begin{defn}[Generalized extensions of full level of $A$] \label{defn: generalized extensions}~\\
\noindent (a) By a generalized extension of full level of $A$ we mean the data of a collection of objects
\[(X_{m,n})_{0\leq m<n\leq k}\]
of $\bT$ with $X_{r-1,r}=A_r$ for all $1\leq r\leq k$, together with a surjective morphism $X_{m,n}\rightarrow X_{m+1,n}$ and an injective morphism $X_{m,n-1}\rightarrow X_{m,n}$ for every $m,n$ in the eligible\footnote{Here and elsewhere throughout, by the adjective {\it eligible} in the context of indices we mean the range in which the indices in question make sense. So here, for instance, we have an injective map $X_{m,n-1}\rightarrow X_{m,n}$ for every pair of integers $(m,n)$ with $0\leq m<n\leq k$ and $m<n-1$.} ranges, such that the following axioms hold:
\begin{itemize}
\item[(i)] Every diagram of the form
\begin{equation}\label{eq9}
\begin{tikzcd}
X_{m,n-1} \arrow[d, hookrightarrow] \arrow[r, twoheadrightarrow] & X_{m+1,n-1} \arrow[d, hookrightarrow] \\
X_{m,n} \arrow[r, twoheadrightarrow] & X_{m+1,n} 
\end{tikzcd}
\end{equation}
(with the maps as in the given data) commutes.
\item[(ii)] The diagram
\begin{equation}\label{eq11}
\begin{tikzcd}
   0 \arrow[r] &  X_{m,n-1}\arrow[r, ] & X_{m,n} \arrow[r, ] &  A_n  \arrow[r] & 0
\end{tikzcd}
\end{equation}
is an exact sequence for every $m,n$ in the eligible range. Here, the morphism $X_{m,n}\rightarrow A_n$ is the composition
\[
X_{m,n} \twoheadrightarrow X_{m+1,n} \twoheadrightarrow X_{m+2,n} \twoheadrightarrow \cdots \twoheadrightarrow X_{n-1,n}=A_n.  
\]
\end{itemize}
\noindent (b) The collection of all generalized extensions of full level of $A$ is denoted by $D_{k-1}(A)$.
\end{defn}

The qualification ``of full level" is included to highlight the contrast with the more general version of this notion, to be introduced shortly below (Definition \ref{def: gen ext level l}). The reason for including axiom (ii) in the definition is to make sure that the $X_{m,n}$ cannot be larger than what we like. With (ii) included as a requirement, one is guaranteed to also get exact sequences
\[
\begin{tikzcd}
   0 \arrow[r] &  A_{m+1}\arrow[r, ] & X_{m,n} \arrow[r, ] &  X_{m+1,n}  \arrow[r] & 0
\end{tikzcd}
\]
(see Lemma \ref{lem: weight filtration of objects in generalized extensions}).

Given an object $X$ of $\bT$ whose associated graded is isomorphic to $A$, choosing an isomorphism $\phi:Gr^WX\rightarrow A$ to identify the two, the subquotients $(X_n/X_m)_{0\leq m<n\leq k}$ with $X_r=W_{p_r}X$ together with the natural successive inclusion and projection maps between them form a generalized extension of full level of $A$. We call this the generalized extension of $A$ associated with $(X, \phi)$ and denote it by $ext(X,\phi)$.

In general, a generalized extension of full level of $A$ can be visualized by a diagram as in \eqref{eq gen ext for M} of \S \ref{sec: statement of the main theorem}. We will simply speak of a generalized extension $(X_{m,n})_{0\leq m<n\leq k}$, or often merely $(X_{m,n})$ or $(X_{\bullet, \bullet})$ without including the arrows or range of indices in the notation. For simplicity and to save space we might sometimes drop the arrows even from our diagrams.

Note that while the definition of a generalized extension $(X_\db)$ only includes maps between objects in adjacent positions in the diagram, for every pairs $(m,n)$ and $(m',n')$ with $m'\geq m$ and $n'\geq n$
by composing the morphisms along any path from $X_{m,n}$ to $X_{m',n'}$ we get a map 
\[
X_{m,n} \rightarrow X_{m',n'}.
\] 
Commutativity of the diagram for $(X_\dbullet)$ guarantees that the outcome does not depend on the choice of the path.

\begin{example} Let $k=3$. Then the data of a generalized extension of full level is the same as the data of a blended extension whose top row is an extension of $A_2$ by $A_1$ and whose right column is an extension of $A_3$ by $A_2$. Indeed, given a generalized extension $(X_{m,n})$ (with $0\leq m<n\leq 3$), we have a blended extension
\[
\begin{tikzcd}
   & & 0 \arrow{d} & 0 \arrow{d} &\\
   0 \arrow[r] & A_1=X_{0,1} \ar[equal]{d} \arrow[r] & X_{0,2}  \arrow[d] \arrow[r] &  X_{1,2}=A_2 \arrow{d} \arrow[r] & 0 \\
   0 \arrow[r] & A_1 \arrow[r] & X_{0,3} \arrow[d] \arrow[r] &  X_{1,3} \arrow{d}  \arrow[r] & 0 \\
   & & A_3 \arrow{d} \ar[equal]{r} & A_3=X_{2,3}, \arrow{d} & \\
   & & 0 & 0 &   
\end{tikzcd}
\]
where the maps are all compositions of the structure maps. The passage from blended extensions to generalized extensions is also clear from this.
\end{example}

Back to working with an arbitrary $k$, we crucially also need truncated versions of the notion, which only include the data of a number of top left to bottom right diagonals of the diagram of equation \eqref{eq gen ext for M} of \S \ref{sec: statement of the main theorem}.

\begin{defn}[Generalized extensions of level $\ell\leq k-1$ of $A$]\label{def: gen ext level l} ~ \\
\noindent (a) Let $1\leq \ell\leq k-1$. By a generalized extension of level $\ell$ of $A$ we mean the data of an object $X_{m,n}$ of $\bT$ for each pair $(m,n)$ of integers with $0\leq m<n\leq k$ and $n-m\leq \ell+1$, with $X_{r-1,r}=A_r$ for all $1\leq r\leq k$, together with the data of a surjective morphism $X_{m,n}\rightarrow X_{m+1,n}$ and an injective morphism $X_{m,n-1}\rightarrow X_{m,n}$ for every $m$ and $n$ in the eligible range such that axioms (i) and (ii) of Definition \ref{defn: generalized extensions}(a) hold.

\noindent (b) The collection of all generalized extensions of level $\ell$ of $A$ is denoted by $D_\ell(A)$. For convenience, we set $D_{k'}(A)=D_{k-1}(A)$ for $k'\geq k$.
\end{defn}

A generalized extension of full level of $A$ is the same as a generalized extension of level $k-1$ of $A$. We visualize a generalized extension of level $\ell$ ($1\leq \ell\leq k-1$) by a (possibly) truncated version of \eqref{eq gen ext for M}, with only $\ell$ diagonals below and parallel to the diagonal consisting of the $A_r$.

\begin{example} A generalized extension of level 1 of $A$ can be visualized as a diagram of the form
\[
\begin{tikzcd}[column sep=0in, row sep=0in]
A_1         & & & &&& \\
X_{0,2} & A_2  & & & &&\\
  & X_{1,3}  & A_3 & & &&\\
  &  & \ddots &\ddots &&\\
&&&&&\\
   && ~ & \hspace{-.2in} X_{k-3,k-1}  &\hspace{-.2in}A_{k-1} &\\
    & & ~ &  & \hspace{-.2in} X_{k-2,k}  &\hspace{-.1in} A_k 
\end{tikzcd}
\]
where the arrows, dropped from the writing for convenience, satisfy axiom (ii) of the definition. This is simply the data of $k-1$ extensions 
\[
\begin{tikzcd}
0 \arrow[r] & A_r \arrow[r] & X_{r-1,r+1} \arrow[r] & A_{r+1} \arrow[r] &  0
\end{tikzcd}
\hspace{.3in}  (1\leq r\leq k-1).
\]
Note that these are short exact sequences, rather than elements of $Ext^1$, since we have not yet introduced any equivalence relations on $D_1(A)$.
\end{example}

We make each $D_\ell(A)$ ($1\leq \ell\leq k-1$) the collection of objects of a category by defining the notion of morphisms of generalized extensions as follows: Let $(X_\dbullet)$ and $(X'_\dbullet)$ be generalized extensions of $A$ of the same level. A morphism of generalized extensions $(X_\dbullet)\rightarrow (X'_\dbullet)$ is a collection of morphisms $f_{m,n}: X_{m,n}\rightarrow  X'_{m,n}$ (one for each pair $(m,n)$ in the eligible range) that
commute with the structure morphisms of $(X_\dbullet)$ and $(X'_\dbullet)$; that is, such that each diagram below commutes for all eligible $(m,n)$:
\begin{equation}\label{eq diagrams for morphisms of gen exts}
\begin{tikzcd}
X_{m,n} \arrow[d, "f_{m,n}"] \arrow[r, twoheadrightarrow] & X_{m+1,n} \arrow[d, "f_{m+1,n}"]\\
X'_{m,n} \arrow[r, twoheadrightarrow] & X'_{m+1,n}
\end{tikzcd}
\hspace{.3in}
\begin{tikzcd}
X_{m,n-1} \arrow[d, "f_{m,n-1}"] \arrow[r, hookrightarrow] & X_{m,n} \arrow[d, "f_{m,n}"]\\
X'_{m,n-1} \arrow[r, hookrightarrow] & X'_{m,n}
\end{tikzcd}
\end{equation}
Note that the morphisms 
\[
f_{r-1,r}: A_r\rightarrow A_r 
\]
here are {\it not} necessarily isomorphisms. With abuse of notation, the category of generalized extensions of level $\ell$ of $A$ is also denoted by $D_\ell(A)$. A morphism $(f_\db)$ in $D_\ell(A)$ is an isomorphism if and only if all the $f_{m,n}$ are isomorphisms. (It will follow from Lemma \ref{lem: weight filtration of objects in generalized extensions} below that $(f_\db)$ is an isomorphism if and only if every $f_{r-1,r}:A_r\rightarrow A_r$ is an isomorphism.)

\medskip\par
\noindent \underline{Truncation and cropping functors}: One can naturally define two types of forgetful functors between categories of generalized extensions. The first are the functors 
\[
\Theta_\ell: D_\ell(A) \rightarrow D_{\ell-1}(A) \hspace{.2in}(2\leq \ell\leq k-1)
\]
defined when $k\geq 3$ as follows: For each $\ell$, the functor $\Theta_\ell$ merely erases the lowest diagonal of each generalized extension of level $\ell$ (i.e. the objects $(X_{r-1,r+\ell})_r$ together with the arrows going into and coming out of these objects). Its action on morphisms of generalized extensions is by restricting the data to the part of the diagrams that survive. We refer to the functors $\Theta_\ell$ as {\it truncation} functors.

The second are the functors that crop diagrams horizontally and vertically to include only the part between two particular graded pieces of $A$. Given integers $i$ and $j$ with $0\leq i,j\leq k$ and $i+1<j$, we have a functor
\begin{equation}\label{eq8}
D_\ell(A) \rightarrow D_\ell(\bigoplus\limits_{i< r\leq j} A_r)
\end{equation}
which only keeps the part of diagrams that lie in the intersection of the columns between $A_{i+1}$ and $A_j$ inclusively, and the rows between $A_{i+1}$ and $A_j$ inclusively. The action on morphisms is again by restricting the data. We call the second type of functors {\it cropping} functors. 

\subsection{Basic properties of generalized extensions}\label{sec: basic properties of gen exts}
In this section we gather some basic results about generalized extensions which will be used in the remainder of the paper. Throughout, unless otherwise indicated, a generalized extension means a generalized extension of $A$ (with $A$ as introduced earlier in \S \ref{sec: statement of problems for S(A) and S'(A)}). We visualize a generalized extension by a diagram as in equation \eqref{eq gen ext for M} of \ref{sec: statement of the main theorem} (with general objects $X_\db$ and arrows that form the data of a generalized extension), or truncated versions of it if the level is less than $k-1$. The references in the text to ``the lowest diagonal", ``the above and right of an object", ``entry $(m,n)$", etc. all refer to this visualization (with the object at entry $(m,n)$ of $(X_\db)$ being $X_{m,n}$). The word ``diagonal"  in this context always refers to the collection of entries $(m,n)$ in the diagram that form a line parallel to the line formed by the $A_r$.

For any generalized extension $(X_\db)$ of any level, for convenience we set $X_{r,r}=0$ for $0\leq r\leq k$. To simplify the writing, as before, if there is no ambiguity we will often simply refer to ``indices in the eligible range" or ``eligible indices"; this simply means that the indices are in the range for which the objects in the equations are available.
\begin{lemma}\label{lem: weight filtration of objects in generalized extensions}
Let $(X_\dbullet)$ be a generalized extension of any level $\ell$.

\noindent (a) For every $m\leq r\leq n$ in the eligible range (depending on the level), we have
\[
Im(X_{m,r}\hookrightarrow X_{m,n}) = W_{p_r}X_{m,n}.
\]
That is, the weight filtration for each object of the diagram is given by the objects directly above it. (Note that in particular, the statement asserts that $W_{p_n}X_{m,n}=X_{m,n}$.)

\noindent (b) The isomorphisms
\begin{equation}\label{eq14}
Gr^WX_{m,n} \cong \bigoplus\limits_{m<r\leq n} A_r
\end{equation}
given by 
\[
Gr^W X_{m,n} \stackrel{\text{(a)}}{\cong} \bigoplus\limits_{m<r\leq n} \frac{X_{m,r}}{X_{m,r-1}} \stackrel{\eqref{eq11}}{\cong} \bigoplus\limits_{m<r\leq n} A_r
\]
for every $(m,n)$ in the eligible range are compatible with the natural injective and surjective maps. That is, we have commutative diagrams
\begin{equation}\label{eq12}
\begin{tikzcd}
Gr^WX_{m,n-1} \ar[equal]{d} 
\arrow[r, hookrightarrow] & Gr^WX_{m,n} \ar[equal]{d}\\
\bigoplus\limits_{m<r\leq n-1} A_r \arrow[r, hookrightarrow] & \bigoplus\limits_{m<r\leq n} A_r
\end{tikzcd}
\hspace{.2in} \text{and} \hspace{.2in} \begin{tikzcd}
Gr^WX_{m,n} \ar[equal]{d} 
\arrow[r, twoheadrightarrow] & Gr^WX_{m+1,n} \ar[equal]{d}\\
\bigoplus\limits_{m<r\leq n} A_r \arrow[r, twoheadrightarrow] & \bigoplus\limits_{m+1<r\leq n} A_r
\end{tikzcd}
\end{equation}
in which the top arrows are $Gr^W$ applied to the structure arrows of $(X_\dbullet)$ and the bottom arrows are the natural embedding and projection maps. The identifications shown as equality in the diagrams are given by the isomorphisms \eqref{eq14}.

\noindent (c) For each $m,n$ in the eligible range, we have an exact sequence
\begin{equation}\label{eq10}
\begin{tikzcd}
   0 \arrow[r] &  A_{m+1}\arrow[r, ] & X_{m,n} \arrow[r, ] &  X_{m+1,n}  \arrow[r] & 0,
\end{tikzcd}
\end{equation}
where the morphism $A_{m+1}\rightarrow X_{m,n}$ is the composition
\[
A_{m+1}=X_{m,m+1}\hookrightarrow X_{m,m+2} \hookrightarrow X_{m,m+3} \hookrightarrow \cdots \hookrightarrow X_{m,n}.
\]

\noindent (d) Let $(f_\dbullet) :(X_\dbullet) \rightarrow (X'_\dbullet)$ be a morphism in $D_\ell(A)$. For each eligible $(m,n)$, the canonical isomorphisms \eqref{eq14} for $(X_\dbullet)$ and $(X'_\dbullet)$ fit into a commutative diagram
\[
\begin{tikzcd}
Gr^WX_{m,n} \ar[equal]{r} \arrow[d, "Gr^Wf_{m,n}"]& \bigoplus\limits_{m<r\leq n} A_r \arrow[d, "(f_{r-1,r})"] \\
Gr^WX'_{m,n} \ar[equal]{r} & \bigoplus\limits_{m<r\leq n} A_r.
\end{tikzcd}
\]
\end{lemma}
\begin{proof}
(a) Fixing $m$, this is seen by induction on $n$ in view of the extension \eqref{eq11}. In the induction step, we first apply the exact functor $W_{p_{n-1}}$ to \eqref{eq11}. Since $W_{p_{n-1}}A_n=0$ we get 
the assertion for $r\leq n-1$. As for when $r=n$, this follows from exactness of $W_{p_n}$ and the fact that $W_{p_n}A_n=A_n$.

\noindent (b) This follows from the construction of the canonical isomorphisms and the commutativity of the diagram of a generalized extension. We leave the details to the reader.

\noindent (c) The exactness of \eqref{eq10} is clear at the first and third object. As for at the middle, tentatively let $B$ be the kernel of $X_{m,n}\twoheadrightarrow X_{m+1,n}$. Applying the exact functor $Gr^W$ to $X_{m,n}\twoheadrightarrow X_{m+1,n}$ it follows from part (b) (see the diagram on the right of equation \eqref{eq12}) that $Gr^WB$ is equal to the subobject $A_{m+1}=Gr^WA_{m+1}$ of $Gr^WX_{m,n}$. It follows that $B=A_{m+1}$.

\noindent (d) Let $m<r\leq n$. By definition of the canonical isomorphism \eqref{eq14} we have a commutative diagram
\[
\begin{tikzcd}[row sep=small]
 & A_r \arrow[dl, twoheadleftarrow] \ar[equal]{dr} &  \\
X_{m,r}  \arrow[rr, twoheadrightarrow] &  & Gr^W_{p_r}X_{m,n} 
\end{tikzcd}
\]
where the horizontal surjective arrow is given by $X_{m,r}\hookrightarrow X_{m,n}$ (mapping $X_{m,r}$ isomorphically to $W_{p_r}X_{m,n}$) and then passing to $Gr^W_{p_r}X_{m,n}$, and the slanted surjective arrow is the composition of the surjective arrows $X_{m',r}\twoheadrightarrow X_{m'+1,r}$ for $m\leq m'<r-1$. The side donated by equality is the identification of \eqref{eq14}. There is an analogous triangle for $(X'_\dbullet)$. The two triangle can be put into a (to be seen to be commutative) diagram:
\begin{equation}\label{eq15}
\begin{tikzcd}[row sep=small]
 & A_r \arrow[dl, twoheadleftarrow] \ar[equal]{dr}  &  \\
X_{m,r}  \arrow[rr, twoheadrightarrow] \arrow[dd, "f_{m,r}"] &  & Gr^W_{p_r}X_{m,n} \arrow[dd, "Gr^W_{p_r}f_{m,n}"] \\
 & A_r \arrow[dl, twoheadleftarrow] \ar[equal]{dr} \arrow[from=uu, crossing over] &  \\
X'_{m,r}  \arrow[rr, twoheadrightarrow] &  & Gr^W_{p_r}X'_{m,n} 
\end{tikzcd}
\end{equation}
where the map $A_r\rightarrow A_r$ on the top is $f_{r-1,r}$. The front and back (triangular) faces are commutative. The rectangular faces on the top left and the bottom are both commutative, the former (resp. latter) by compatibility of morphisms of generalized extensions with the surjective (resp. injective) structure arrows. It follows that the top right face is also commutative.
\end{proof}
Note that in particular, the previous lemma asserts that for every generalized extension $(X_\dbullet)$ of level $k-1$ of $A$ we have 
\[
Gr^WX_{0,n} \cong \bigoplus\limits_{1\leq r\leq n} A_r
\]
for all $1\leq n\leq k$. 

Part (d) of the previous lemma has the following consequence:
\begin{lemma}\label{lem: if identity on A identity everywhere}
Let $(X_\db)$ be a generalized extension of level $\ell$ of $A$. The forgetful map
\[
Aut((X_\db)) \rightarrow Aut(A) = \prod\limits_{1\leq r\leq k}Aut(A_r)
\]
given by $(\sigma_\db)\mapsto (\sigma_{r-1,r})$ is injective. (That is, every automorphism of $(X_\db)$ is uniquely determined by its action on $A$. Of course, the forgetful map above need not be surjective.)
\end{lemma}
\begin{proof}
Suppose $(\sigma_\dbullet)$ is an automorphism of $(X_\db)$ which is identity on $A$. Applying Lemma \ref{lem: weight filtration of objects in generalized extensions}(d) with $(X'_\dbullet)=(X_\dbullet)$ and $(f_\db)=(\sigma_\db)$ we obtain that each $Gr^W\sigma_{m,n}$ and hence $\sigma_{m,n}$ is identity.
\end{proof}

Before we proceed any further, let us introduce a notation: 

\begin{notation}\label{notation: horizontal and vertical extensions}
Given a generalized extension $(X_\dbullet)$ of $A$ of any level, we denote the two extensions \eqref{eq11} and \eqref{eq10} of Definition \ref{defn: generalized extensions} and Lemma \ref{lem: weight filtration of objects in generalized extensions} respectively by $\sX_{m,n}^v$ and $\sX_{m,n}^h$. Thus $\sX_{m,n}^v$ (resp. $\sX_{m,n}^h$) is the extension induced by the vertical arrow into (resp. the horizontal arrow coming out of) $X_{m,n}$.
\end{notation}

The following two lemmas will be useful in constructing morphisms between generalized extensions. The first lemma asserts that given two generalized extensions $(X_\dbullet)$ and $(X'_\dbullet)$ of the same level, every morphism from an object of $(X_\dbullet)$ to an object of $(X'_\dbullet)$ at the same entry extends (or spreads) uniquely to the part of the diagrams to the above and right of that entry.
\begin{lemma} \label{lem: morphisms spread to the top right}
Let $(X_\dbullet)$ and $(X'_\dbullet)$ be generalized extensions of level $\ell$ of $A$. 

\noindent (a) Given any $i,j$ in the eligible range with $j-i>1$ (so that $X_{i,j}$ is below the diagonal of the $A_r$) and a morphism $f: X_{i,j}\rightarrow  X'_{i,j}$, there exists a unique collection of morphisms 
\[
f_{m,n}: X_{m,n}\rightarrow  X'_{m,n} \hspace{.3in} (i\leq m<n\leq j)
\]
such that $f_{i,j}=f$ and the $f_{m,n}$ commute with the morphisms in $(X_\dbullet)$ and $(X'_\dbullet)$. (That is, such that the $f_{m,n}$ give a morphism between the parts of $(X_\dbullet)$ and $(X'_\dbullet)$ between $A_{i+1}$ and $A_j$. See \eqref{eq8}.)

\noindent (b) The extension of $f$ to $(f_{m,n})$ as above behaves well with respect to compositions: if $(X''_\dbullet)$ is also a generalized extension of level $\ell$ of $A$ and $f': X'_{i,j}\rightarrow  X''_{i,j}$ is a morphism, then
\[
(f'\circ f)_{m,n} = f'_{m,n}\circ f_{m,n}
\]
for all eligible $m,n$.
\end{lemma}
\begin{proof}
We may assume that the level is $k-1$ and $(i,j)=(0,k)$, as the seemingly more general statements follow from applying this case to the cropped generalized extensions (cropped in the sense of \eqref{eq8}).

By Lemma \ref{lem: weight filtration of objects in generalized extensions}, for any $n\leq k$ we have $W_{p_n}X_{0,k}\cong X_{0,n}$ and $W_{p_n}X'_{0,k}\cong X'_{0,n}$ (identifications via the injective structure arrows of the generalized extensions). By functoriality of the weight filtration, the morphism $f$ restricts to morphisms $f_{0,n}:X_{0,n}\rightarrow X'_{0,n}$ compatible with each other. So far, we have extended $f$ to the first column of the two generalized extensions.

Assume $f$ has been extended in the desired way to the $f_{m,n}$ for all $m<m'$. For each eligible $n$ we have a commutative diagram
\[
\begin{tikzcd}
  0 \arrow[r] &  A_{m'}\arrow[r, ] \arrow[d, "f_{m'-1,m'}"]& X_{m'-1,n} \arrow[r, ] \arrow[d, "f_{m'-1,n}"] &  X_{m',n}  \arrow[r] & 0\\
  0 \arrow[r] &  A_{m'}\arrow[r, ] & X'_{m'-1,n} \arrow[r, ] &  X'_{m',n}  \arrow[r] & 0
\end{tikzcd}
\]
with rows $\sX_{m'-1,n}^h$ and ${\sX'}_{m'-1,n}^h$. We thus get a morphism $f_{m',n}: X_{m',n} \rightarrow X'_{m',n}$ making this diagram commute.  

We have extended $f$ to the column of entries $(m',n)$ (with $m'$ fixed), and the commutativity with the surjective arrows from the column for $m'-1$ to the one for $m'$ is known by construction. We have a diagram
\[
\begin{tikzcd}[row sep=small, column sep=small]
& X_{m'-1,n-1}  \arrow[rr, twoheadrightarrow] \arrow[dd] & & X_{m',n-1}  \arrow[dd, "f"] \\ 
X_{m'-1,n} \arrow[ur, hookleftarrow] \arrow[rr, crossing over, twoheadrightarrow] \arrow[dd, "f"] & & X_{m',n} \arrow[ur, hookleftarrow] \\
& X'_{m'-1,n-1}  \arrow[rr, twoheadrightarrow] & & X'_{m',n-1} \\
X'_{m'-1,n} \arrow[ur, hookleftarrow] \arrow[rr, twoheadrightarrow ] & & X'_{m',n} \arrow[from=uu, crossing over] \arrow[ur, hookleftarrow]\\
\end{tikzcd}
\]
where the top (resp. bottom) face is a part of $(X_\dbullet)$ (resp. $(X'_\dbullet)$), and the downward maps are the maps induced by $f$. We know the top and bottom faces as well as the left, back and front faces are commutative. In view of the surjectivity of $X_{m'-1,n-1}\twoheadrightarrow X_{m',n-1}$ we get the commutativity of the right face. This completes the proof of the fact that $f$ extends to a morphism of generalized extensions.

As for uniqueness, the map $f_{0,k}$ determines every map $f_{m,n}$ with $0\leq m<n\leq k$ because of the commutativity requirements. Indeed, $f_{0,k}$ determines each $f_{0,n}$ and in turn, $f_{0,n}$ determines each $f_{m,n}$ by the commutativity of the following two diagrams:
\[
\begin{tikzcd}
X_{0,n} \arrow[r, hookrightarrow] \arrow[d, "f_{0,n}"] &  X_{0,k}  \arrow[d, "f_{0,k}"]  \\
X'_{0,n} \arrow[r, hookrightarrow] &  X'_{0,k}  
\end{tikzcd}
\hspace{.3in}
\begin{tikzcd}
X_{0,n} \arrow[r, twoheadrightarrow] \arrow[d, "f_{0,n}"] &  X_{m,n}  \arrow[d, "f_{m,n}"]  \\
X'_{0,n} \arrow[r, twoheadrightarrow] &  X'_{m,n}  
\end{tikzcd}
\]

Part (b) is easily seen from the construction of the $f_{m,n}$ given above.
\end{proof}

The next lemma says that morphisms between the lowest diagonals of two generalized extensions of the same level glue together to give a morphism of generalized extensions if and only if they agree on the diagonal just above the lowest.

\begin{lemma}\label{lem: gluing morphisms between objects on the lowest diagonal}
Suppose $(X_\dbullet)$ and $(X'_\dbullet)$ are generalized extensions of level $\ell$ of $A$. Suppose that for each eligible pair $(i,j)$ with $j-i=\ell+1$ (i.e. on the lowest diagonal) we have a morphism $f_{i,j}: X_{i,j}\rightarrow X'_{i,j}$. For each such $(i,j)$, let $f^v_{i,j-1}$ and $f^h_{i+1,j}$ be the unique morphisms fitting in the commutative diagrams\footnote{The intuition behind the notation is the following: $f^v_{i,j-1}$ (resp. $f^h_{i+1,j}$) is the map between the objects in $(i,j-1)$ (resp. $(i+1,j)$) entries determined by $f_{i,j}$ and the vertical arrows into (resp. horizontal arrows coming out of) the $(i,j)$ entries. That $f^v_{i,j-1}$ and $f^h_{i+1,j}$ exist is by Lemma \ref{lem: weight filtration of objects in generalized extensions} and functoriality of the weight filtration.}
\[
\begin{tikzcd}
X_{i,j-1} \arrow[d, "f^v_{i,j-1}"] \arrow[r, hookrightarrow] & X_{i,j} \arrow[d, "f_{i,j}"]\\
X'_{i,j-1} \arrow[r, hookrightarrow] & X'_{i,j}
\end{tikzcd}
\hspace{.4in}
\begin{tikzcd}
X_{i,j} \arrow[d, "f_{i,j}"] \arrow[r, twoheadrightarrow] & X_{i+1,j} \arrow[d, "f^h_{i+1,j}"]\\
X'_{i,j} \arrow[r, twoheadrightarrow] & X'_{i+1,j}.
\end{tikzcd}
\]
Then there exists a morphism of generalized extensions $(X_\dbullet)\rightarrow (X'_\dbullet)$ extending the given $f_{i,j}$ from the lowest diagonals to the entire diagrams if and only if for every eligible $(m,n)$ with $n-m=\ell$ (i.e. on the diagonal just above the lowest), we have $f^v_{m,n}=f^h_{m,n}$. Moreover, when such an extension exists, it is unique. (In other words, if two morphisms of generalized extensions agree on the lowest diagonal, then the two morphisms are the same.)
\end{lemma} 
\begin{proof}
The uniqueness is immediate from the uniqueness statement in Lemma \ref{lem: morphisms spread to the top right}. That the given condition is necessary follows from the same result and its proof. The new assertion here is that the compatibility condition on the diagonal just above the lowest is sufficient for the morphisms between the lowest diagonals to glue to make a morphism of generalized extensions. 

By Lemma \ref{lem: morphisms spread to the top right} the morphism $f_{0,\ell+1}$ extends to the right and top of entry $(0,\ell+1)$. Assume that the morphisms on the lowest diagonal glue all the way up to and including the entry $(i,j)$ on the lowest diagonal, so that for all eligible pairs $(m,n)$ of indices with $n\leq j$ we already have morphisms $f_{m,n}: X_{m,n}\rightarrow X'_{m,n}$ commuting with the structure injections and surjections. We will argue that if $j\neq k$, we can also glue $f_{i+1,j+1}$ to the current data.

Consider the map $f_{i+1,j}: X_{i+1,j}\rightarrow X'_{i+1,j}$ (which is already available). Then $f_{i+1,j}$ is induced by $f_{i,j}$ and hence is $f^h_{i+1,j}$. By Lemma \ref{lem: morphisms spread to the top right}, $f_{i+1,j+1}$ extends uniquely to a morphism of generalized extensions between the parts of the diagrams above and to the right of entry $(i+1,j+1)$, in particular, inducing the map $f^v_{i+1,j}: X_{i+1,j}\rightarrow X'_{i+1,j}$. By the compatibility condition, $f^h_{i+1,j}=f^v_{i+1,j}$. Applying the uniqueness statement of Lemma \ref{lem: morphisms spread to the top right} at entry $(i+1,j)$ it follows that the map induced by $f_{i+1,j+1}$ at every entry $(m,n)$ with $m\geq i+1$ and $n\leq j$ coincides with the map $f_{m,n}$ already there. Together with the new maps $f_{m,j+1}$ for $m\geq i+1$ induced by $f_{i+1,j+1}$ we have extended the maps between the diagrams one row further down.

There is nothing new to check for the commutativity with the structure injections and surjections: every square that needs to commute is already known to commute.
\end{proof}

\subsection{Equivalence relations on generalized extensions}\label{sec: equiv rels on gen exts}
In this subsection we define the sets $S'_\ell(A)$ and $S_\ell(A)$ and the maps from level $\ell$ to $\ell-1$ in Theorem \ref{thm: main thm 1}. Recall from \S \ref{sec: gen exts, defn} that for each integer $\ell$ with $1\leq \ell\leq k-1$ the collection (also the category) of generalized extensions of level $\ell$ of $A$ is denoted by $D_\ell(A)$. As before, we may sometimes drop the phrase ``of $A$".

There are two natural equivalence relations on each $D_\ell(A)$. The first is simply given by isomorphisms in the category $D_\ell(A)$, and the second is the finer equivalence given by isomorphisms that are identity on $A$.
\begin{notation}
Let $(X_\dbullet)$ and $(X'_\dbullet)$ be generalized extensions of the same level.

\noindent (a) We write $(X_\dbullet) \sim (X'_\dbullet)$ if there exists an isomorphism of generalized extensions $(X_\dbullet) \rightarrow (X'_\dbullet)$. That is, if there an isomorphism $f_{m,n}: X_{m,n}\rightarrow X'_{m,n}$ for each pair $(m,n)$  in the eligible range such that the diagrams \eqref{eq diagrams for morphisms of gen exts} commute for all $(m,n)$. 

\noindent (b) We write $(X_\dbullet) \sim' (X'_\dbullet)$ if there exists an isomorphism $(f_\db): (X_\dbullet) \rightarrow (X'_\dbullet)$ that is identity on $A$, i.e. such that for every $r$ the isomorphism $f_{r-1,r}: A_r\rightarrow A_r$ is the identity map. 

\noindent (c) For each $1\leq \ell\leq k-1$, denote the set of equivalence classes of objects of $D_\ell(A)$ with respect to $\sim$ (resp. $\sim'$) by $S_\ell(A)$ (resp. $S'_\ell(A)$).
\end{notation}

There is a natural surjection $S'_\ell(A)\rightarrow S_\ell(A)$ induced by the identity map on $D_\ell(A)$. Note that by Lemma \ref{lem: if identity on A identity everywhere}, if $(X_\dbullet)$ and $(X'_\dbullet)$ are generalized extensions of level $\ell$ such that $(X_\dbullet) \sim' (X'_\dbullet)$, then there is only one isomorphism $(X_\dbullet) \rightarrow (X'_\dbullet)$ that is identity on $A$.

Recall that for each $2\leq \ell\leq k-1$ we have a truncation functor $\Theta_\ell: D_\ell(A)\rightarrow D_{\ell-1}(A)$, simply erasing the lowest diagonal of a generalized extension. The truncation functors clearly preserve both $\sim$ and $\sim'$, inducing maps $S'_\ell(A)\rightarrow S'_{\ell-1}(A)$ and $S_\ell(A)\rightarrow S_{\ell-1}(A)$ both of which we shall also refer to as truncation maps and (with abuse of notation) denote by $\Theta_\ell$. For each $\ell$, we have a commutative diagram
\begin{equation}\label{eq101}
\begin{tikzcd}
D_{\ell}(A) \arrow[d, twoheadrightarrow] \arrow[r, "\Theta_{\ell}"] & D_{\ell-1}(A)\arrow[d, twoheadrightarrow] \\
S'_\ell(A) \arrow[d, twoheadrightarrow] \arrow[r, "\Theta_{\ell}"] & S'_{\ell-1}(A) \arrow[d, twoheadrightarrow] \\
S_\ell(A) \arrow[r, "\Theta_{\ell}"] & S_{\ell-1}(A)
\end{tikzcd}
\end{equation}
where the vertical arrows are the natural maps: modding out by $\sim'$ first and then further by $\sim$.

The sets $S_\ell(A)$ and the maps $\Theta_\ell: S_\ell(A)\rightarrow S_{\ell-1}(A)$ were characterized differently in the statement of Theorem \ref{thm: main thm 1}. We now describe the action of $Aut(A)$ on $S'_\ell(A)$ and discuss the equivalence of definitions given here and the ones in the statement of Theorem \ref{thm: main thm 1}.

The group $Aut(A)$ acts on $D_\ell(A)$ by twisting the arrows to and from the $A_r$: making the action a left action as usual, $\sigma=(\sigma_r)$ in $Aut(A)=\prod_r Aut(A_r)$ sends the diagram below on the left to the one on the right. Note that the rest of the arrows in the diagram remain unchanged.
\begin{equation}\label{eq17}
\begin{tikzcd}[column sep=small, row sep=small]
A_1 \arrow[d, hookrightarrow, "j_1" ']   & & & &&& \\
X_{0,2} \arrow[r, twoheadrightarrow, "\omega_2"] & A_2 \arrow[d, hookrightarrow, "j_2" '] & & & &&\\
\vdots  & X_{1,3} \arrow[r, twoheadrightarrow, "\omega_3"]  & A_3 \arrow[d, hookrightarrow, "j_3" '] & & &&\\
  &   & \ddots &\ddots &&\\
   && ~ & \hspace{-.2in} X_{k-3,k-1} \arrow[r, twoheadrightarrow, "\omega_{k-1}"]  & A_{k-1} \arrow[d, hookrightarrow, "j_{k-1}" ']&\\
    & &  &\cdots  & \hspace{-.2in} X_{k-2,k} \arrow[r, twoheadrightarrow, "\omega_{k}"] & A_k 
\end{tikzcd}
\hspace{-1.5in}
\begin{tikzcd}[column sep=0.3in, row sep=small]
A_1 \arrow[d, hookrightarrow, "j_1\sigma_1^{-1}" ']   & & & &&& \\
X_{0,2} \arrow[r, twoheadrightarrow, "\sigma_2\omega_2" ] & A_2 \arrow[d, hookrightarrow, "j_2\sigma_2^{-1}" '] & & & &&\\
\vdots  & X_{1,3} \arrow[r, twoheadrightarrow, "\sigma_3\omega_3"]  & A_3 \arrow[d, hookrightarrow, "j_3\sigma_3^{-1}" '] & & &&\\
  &   & \ddots &\ddots &&\\
   && ~ & \hspace{-.2in} X_{k-3,k-1} \arrow[r, twoheadrightarrow, "\sigma_{k-1}\omega_{k-1}"]  & A_{k-1} \arrow[d, hookrightarrow, "j_{k-1}\sigma_{k-1}^{-1}" ']&\\
    & &  & \cdots & \hspace{-.2in} X_{k-2,k} \arrow[r, twoheadrightarrow, "\sigma_{k}\omega_{k}"] & A_k 
\end{tikzcd}
\end{equation}
We use the notation $\sigma\cdot (X_\db)$ for the image of $(X_\db)$ under $\sigma\in Aut(A)$.

If $(f_\dbullet): (X_\dbullet)\rightarrow (X'_\dbullet)$ is an isomorphism that is identity on $A$, then $(f_\dbullet)$ is also such an isomorphism $\sigma\cdot (X_\dbullet)\rightarrow \sigma\cdot (X'_\dbullet)$. Indeed, if either of the two diagrams
\[
\begin{tikzcd}
A_r \ar[equal]{d} \arrow[r, hookrightarrow, "j_r"] & X_{r-1,r+1} \arrow[d, "f_{r-1,r+1}"]\\
A_r \arrow[r, hookrightarrow, "j'_r"] & X'_{r-1,r+1}
\end{tikzcd}
\hspace{.3in}\text{or}
\hspace{.3in}
\begin{tikzcd}
X_{r-2,r} \arrow[d, "f_{r-2,r}"] \arrow[r, twoheadrightarrow, "\omega_r"] & A_r \ar[equal]{d}\\
X'_{r-2,r} \arrow[r, twoheadrightarrow, "\omega'_r"] & A_r
\end{tikzcd}
\]
commutes, then so it does after twisting the horizontal arrows by $\sigma_r$ or $\sigma_r^{-1}$ in $Aut(A_r)$. Thus the action of $Aut(A)$ on $D_\ell(A)$ descends to an action on $S'_\ell(A)$. Moreover, it is clear from the definitions that the top horizontal map in \eqref{eq101} is $Aut(A)$-equivariant; thus so is the truncation map $S'_\ell(A)\rightarrow S'_{\ell-1}(A)$.

The equivalence of the two definitions of $S_\ell(A)$ and $\Theta_\ell:S_\ell(A)\rightarrow S_{\ell-1}(A)$ given in this subsection and the ones given in the statement of Theorem \ref{thm: main thm 1} can now be seen from the following lemma:
\begin{lemma}\label{lem: relation betweem oD and oD'}
Let $1\leq \ell\leq k-1$. The natural surjection
\[
S'_\ell(A) \twoheadrightarrow S_\ell(A)
\]
descending from the identity map on $D_\ell(A)$ descends further to a bijection
\[
S'_\ell(A)/Aut(A) \cong S_\ell(A).
\]
\end{lemma}
\begin{proof}
We will show that two generalized extensions $(X_\db)$ and $(X'_\db)$ are $\sim$-equivalent (i.e. isomorphic) if and only if there exists $\sigma\in Aut(A)$ such that $\sigma\cdot (X_\db)\sim' (X'_\db)$.

Suppose $(f_\db): (X_\db) \rightarrow (X'_\db)$ is an isomorphism. Then $\sigma_r:=f_{r-1,r}$ is an automorphism of $A_r$. For each eligible pair $(m,n)$, let $g_{m,n}=f_{m,n}$ if $n-m>1$. Let $g_{r-1,r}$ be the identity on $A_r$. Set $\sigma=(\sigma_r)\in Aut(A)$. Then $(g_\db): \sigma\cdot (X_\db) \rightarrow (X'_\db)$ is an isomorphism of generalized extensions that is identity on $A$. Thus $\sigma\cdot (X_\db)\sim'  (X'_\db)$.

Conversely, assume $\sigma\cdot (X_\db)\sim'  (X'_\db)$. Let $(g_\db): \sigma\cdot (X_\db) \rightarrow (X'_\db)$ be an isomorphism that is identity on $A$. Then $f_\dbullet: (X_\dbullet) \rightarrow (X'_\dbullet)$ defined by $f_{m,n}=g_{m,n}$ if $n-m>1$ and $f_{r-1,r}=\sigma_r$ is an isomorphism.
\end{proof}

In all that will follow, we adopt the definitions given in this subsection for the sets $S_\ell(A)$ and the maps between them.

\subsection{Generalized extensions in levels 1 and $k-1$}\label{sec: levels 1 and k-1}
We now study the sets $S_\ell(A)$ and $S'_\ell(A)$ when $\ell$ is 1 or $k-1$. The goal is to establish the characterizations given in parts (a) and (d) of Theorem \ref{thm: main thm 1}.

We start by recalling the action of
\[
Aut(A) = \prod\limits_{r}Aut(A_r)
\]
on
\[
\prod\limits_{r} Ext^1(A_{r+1},A_{r})
\]
that appeared in the statement of Theorem \ref{thm: main thm 1}(d). Given $\sigma=(\sigma_r)\in Aut(A)$ and 
\[(\mathcal{E}_r) \in \prod\limits_{r} Ext^1(A_{r+1},A_{r}),\]
the element $\sigma\cdot (\mathcal{E}_r)$ is the element whose $r$-entry is ${(\sigma_{r})}_\ast{(\sigma_{r+1}^{-1})}^\ast \, \sE_r$. Note that after taking a representative for $\sE_r$, i.e. lifting it an extension 
\[
\begin{tikzcd}
0 \arrow[r] & A_r \arrow[r, "j_r"] & E_r \arrow[r, "\omega_{r+1}"] & A_{r+1} \arrow[r] & 0,
\end{tikzcd}
\]
${(\sigma_{r})}_\ast{(\sigma_{r+1}^{-1})}^\ast \, \sE_r$ is the class of the extension obtained by replacing $\omega_{r+1}$ (resp. $j_r$) by $\sigma_{r+1}\omega_{r+1}$ (resp. $j_r \sigma_r^{-1}$).

As we already observed, by definition, the data of a generalized extension of level 1 of $A$ is equivalent to the data of a collection of objects $X_{0,2},X_{1,3},\ldots, X_{k-2,k}$ and short exact sequences
\begin{equation}\label{eq16}
\begin{tikzcd}
0\arrow[r] & A_r \arrow[r] & X_{r-1,r+1} \arrow[r] & A_{r+1} \arrow[r] & 0
\end{tikzcd}\hspace{.5in}( 1\leq r\leq k-1).
\end{equation}
Referring to the notation earlier introduced (see Notation \ref{notation: horizontal and vertical extensions}), the extension above is both $\sX^h_{r-1,r+1}$ and $\sX^v_{r-1,r+1}$. A morphism $(f_\dbullet) :(X_\dbullet)\rightarrow (X'_\dbullet)$ of generalized extensions of level 1 is the data of morphisms $f_{r-1,r+1}: X_{r-1,r+1} \rightarrow  X'_{r-1,r+1}$ and $f_{r-1,r}: A_r\rightarrow A_r$ (for each $r$) making the diagrams
\[
\begin{tikzcd}
0\arrow[r] & A_r \arrow[d, "f_{r-1,r}"] \arrow[r] & X_{r-1,r+1} \arrow[r] \arrow[d, "f_{r-1,r+1}"] & A_{r+1} \arrow[d, "f_{r,r+1}"] \arrow[r] & 0\\0\arrow[r] & A_r \arrow[r] & X_{r-1,r+1} \arrow[r] & A_{r+1} \arrow[r] & 0
\end{tikzcd}
\]
commute. By definition, two generalized extensions $(X_\dbullet)$ and $(X'_\dbullet)$ of level 1 are $\sim'$-equivalent if and only if there are morphisms $f_{r-1,r+1}$ that together with the identity maps on the $A_r$ make the diagrams commute. In other words, $(X_\dbullet) \sim' (X'_\dbullet)$ if and only if for every $r$ the extensions $\sX^h_{r-1,r+1}$ and ${\sX'}^h_{r-1,r+1}$ represent the same element in the corresponding $Ext^1$ group. This is summarized in part (a) below. Part (b) of the statement follows from the fact that the actions of $Aut(A)$ on both $S'(A)$ and $\prod_r Ext^1(A_{r+1},A_r)$ are given by twisting the same arrows in the same way, and the fact (just observed) that $\sim'$ translates to the usual equivalence of 1-extensions.
\begin{lemma}\label{lem: D_2}
(a) Two generalized extensions $(X_\dbullet)$ and $(X'_\dbullet)$ of level 1 are $\sim'$-equivalent if and only if for each $r$ the extensions $\sX^v_{r-1,r+1}$ and ${\sX'}^v_{r-1,r+1}$ (i.e. \eqref{eq16} and its counterpart for $(X'_\dbullet)$) coincide in $Ext^1(A_{r+1},A_r)$. The association $(X_\dbullet)\mapsto (\sX^h_{r-1,r+1})_r$ induces a bijection
\[
S'_1(A) \xrightarrow{ \ \simeq \ } \prod\limits_r Ext^1(A_{r+1},A_r).
\]

\noindent (b) Considering the previously defined actions of $Aut(A)$ on $S'_1(A)$ and $\prod\limits_r Ext^1(A_{r+1},A_r)$, the bijection of part (a) is $Aut(A)$-equivariant and descends to a bijection
\[
S_1(A) \xrightarrow{ \ \simeq \ } \left(\prod\limits_r Ext^1(A_{r+1},A_r)\right)\bigm/ Aut(A).
\]
\end{lemma}
We now turn our attention to the case $\ell=k-1$. Recall from \S \ref{sec: gen exts, defn} that for every pair $(X,\phi)$ of an object $X$ of $\bT$ whose associated graded is isomorphic to $A$ and an isomorphism $\phi:Gr^WX\rightarrow A$, we have an associated generalized extension $ext(X,\phi)$ of level $k-1$. The object at $(m,n)$ entry of $ext(X,\phi)$ is $W_{p_n}X/W_{p_m}X$, with the graded component $W_{p_r}X/W_{p_{r-1}}X$ identified with $A_r$ via $\phi$. The structure morphisms in $ext(X,\phi)$ are the natural injections and projections. The following statement is easily seen from the definitions: 
\begin{lemma}\label{lem: can iso between Gr(X) and A for ext(M,phi)}
For any pair $(X,\phi)$ as above, the canonical isomorphism \eqref{eq14} of Lemma \ref{lem: weight filtration of objects in generalized extensions} for $ext(X,\phi)$ with $(m,n)=(0,k)$ is $\phi$.
\end{lemma}
It is easily seen from the constructions that for every $\sigma\in Aut(A)$,
\begin{equation}\label{eq19}
ext(\sigma\cdot(X,\phi)) = ext(X,\sigma\phi) = \sigma\cdot ext(X,\phi).
\end{equation}
Note that here $\sigma\cdot(X,\phi)$ refers to the action of $Aut(A)$ on the collection of pairs $(X,\phi)$. This action was defined by $\sigma\cdot(X,\phi)=(X,\sigma\phi)$ (see \S \ref{sec: statement of problems for S(A) and S'(A)}).

Recall from \S \ref{sec: statement of problems for S(A) and S'(A)} that two pairs $(X,Gr^WX\xrightarrow{\phi,\simeq}A)$ and $(X',Gr^WX'\xrightarrow{\phi',\simeq}A)$ are said to be equivalent if there exists an isomorphism $f: X\rightarrow X'$ such that $\phi'\circ Gr^Wf=\phi$. Also recall that we denoted the set of equivalence classes of such pairs by $S'(A)$, and that the action of $Aut(A)$ on the collection of pairs $(X,\phi)$ descends to an action on $S'(A)$.
\begin{lemma}\label{lem: D_k}
(a) Two pairs $(X,\phi)$ and $(X',\phi')$ as above are equivalent if and only if the generalized extensions $ext(X,\phi)$ and $ext(X',\phi')$ are $\sim'$-equivalent. 

\noindent (b) Let $(X_\db)$ be a generalized extension of level $k-1$. Let $\phi: Gr^WX_{0,k}\rightarrow A$ be the canonical isomorphism \eqref{eq14} (of Lemma \ref{lem: weight filtration of objects in generalized extensions}) for $(X_\db)$ and $(m,n)=(0,k)$. Then the identity map on $X_{0,k}$ extends to an isomorphism $(X_\db)\rightarrow ext(X_{0,k}, \phi)$ that is identity on $A$.

\noindent (c) The association $(X,\phi)\mapsto ext(X,\phi)$ descends to a bijection
\[
S'(A) \xrightarrow{\ \simeq \ } S'_{k-1}(A).
\]

\noindent (d) Considering the previously defined actions of $Aut(A)$ on $S'(A)$ and $S'_{k-1}(A)$, the bijection of part (c) is $Aut(A)$-equivariant and it descends to a bijection   
\[
S'(A)/Aut(A) \xrightarrow{\ \simeq \ } S'_{k-1}(A)/Aut(A).
\]

\noindent (e) Using the bijections of Lemmas \ref{lem: relation between S'(A) and S(A)} and \ref{lem: relation betweem oD and oD'} to translate the bijection of part (d) to a map 
\[S(A) \xrightarrow{ \ \simeq \ } S_{k-1}(A),\]
this bijection is described as follows: It sends the isomorphism class of $X$ (an object whose associated graded is isomorphic to $A$) to the isomorphism class (i.e. image in $S_{k-1}(A)$) of the generalized extension $ext(X,\phi)$ for any choice of isomorphism $\phi:Gr^WX\rightarrow A$.
\end{lemma}
\begin{proof}
(a) Suppose $(X,\phi)$ and $(X',\phi')$ are equivalent, with $f:X\rightarrow X'$ an isomorphism for which $\phi'\circ Gr^Wf=\phi$. By Lemma \ref{lem: morphisms spread to the top right}, $f$ extends uniquely to an isomorphism $(f_\db): ext(X,\phi)\rightarrow ext(X',\phi')$. In view of Lemmas \ref{lem: can iso between Gr(X) and A for ext(M,phi)} and \ref{lem: weight filtration of objects in generalized extensions}(d) (the latter applied with $(m,n)=(0,k)$), the fact that $\phi'\circ Gr^Wf=\phi$ implies that $(f_\db)$ is identity on $A$.

Conversely, suppose $(f_\dbullet):ext(X,\phi)\rightarrow ext(X',\phi')$ is an isomorphism that is identity on $A$. Then $f_{0,k}: X\rightarrow X'$ is an isomorphism that satisfies $\phi' Gr^Wf_{0,k}=\phi$. This again follows from Lemmas \ref{lem: can iso between Gr(X) and A for ext(M,phi)} and \ref{lem: weight filtration of objects in generalized extensions}(d).

\noindent (b) By Lemma \ref{lem: morphisms spread to the top right}, the identity map on $X_{0,k}$ extends uniquely to an isomorphism $(f_\dbullet): (X_\dbullet) \rightarrow ext(X_{0,k},\phi)$. Now apply Lemma \ref{lem: weight filtration of objects in generalized extensions}(d) with $(m,n)=(0,k)$, $(X'_\dbullet)=ext(X_{0,k},\phi)$ and $(f_\dbullet)$ as in here. The left arrow of the diagram is $Gr^Wf_{0,k}=Gr^WId$ which is just the identity. The top and bottom canonical isomorphisms are both $\phi$. Hence the arrow on the right is the identity map.

\noindent (c) By part (a), $(X,\phi)\mapsto ext(X,\phi)$ descends to an injection $S'(A) \rightarrow S'_{k-1}(A)$, which is also surjective by part (b).

\noindent (d) We have a commutative diagram
\[
\begin{tikzcd}[column sep=large]
\bigm\{(X,\phi)\bigm\} \arrow[r, "ext(-)"] \arrow[d, twoheadrightarrow] & D_{k-1}(A) \arrow[d, twoheadrightarrow]\\
S'(A) \arrow[r, "\text{part (b)}" ,  "\simeq" '] & S'_{k-1}(A),
\end{tikzcd}
\]
where $\{(X,\phi)\}$ means the collection of all pairs $(X,Gr^WX\xrightarrow{\phi,\simeq} A)$. By \eqref{eq19}, the top arrow is $Aut(A)$-equivariant. By definition of the $Aut(A)$-actions on $S'(A)$ and $S'_{k-1}(A)$, so are the two side arrows. It follows that the map of part (b) is $Aut(A)$-equivariant and it descends to a map 
\[
S'(A)/Aut(A) \rightarrow S'_{k-1}(A)/Aut(A),
\]
which is surjective thanks to part (b). It remains to show that it is also injective.

Consider two pairs $(X,\phi)$ and $(X',\phi')$ such that the classes of $ext(X,\phi)$ and $ext(X',\phi')$ in $S'_{k-1}(A)$ are in the same $Aut(A)$-orbit. In view of \eqref{eq19} this means that there exists $\sigma\in Aut(A)$ such that $ext(X,\sigma\phi)$ is $\sim'$-equivalent to $ext(X',\phi')$. That is, there exists an isomorphism $(f_\dbullet)$ of generalized extensions $ext(X,\sigma\phi)\rightarrow  ext(X',\phi')$ such that for each $r$, the morphism $f_{r-1,r}$ is the identity map on $A_r$. 

Consider $f_{0,k}: X\rightarrow X'$. We claim that $\phi'\circ Gr^Wf_{0,k}=\sigma\phi$; this would show that pairs $(X,\sigma\phi)=\sigma\cdot (X,\phi)$ and $(X',\phi')$ are equivalent, so that the classes of $(X,\phi)$ and $(X',\phi')$ in $S'(A)$ are in the same $Aut(A)$-orbit. 

To see the claim, apply Lemma \ref{lem: weight filtration of objects in generalized extensions}(d) to $(f_\dbullet)$ for $(m,n)=(0,k)$. Since the canonical isomorphism \eqref{eq14} for $ext(X,\sigma\phi)$ (resp. $ext(X', \phi')$) is simply $\sigma\phi$ (resp. $\phi'$), we get that the diagram
\[
\begin{tikzcd}
Gr^WX \arrow[d, "Gr^Wf_{0,k}" '] \arrow[r, "\sigma\phi" , "\simeq" ' ] & A \arrow[d, "(f_{r-1,r}) = Id"] \\
Gr^WX' \arrow[r, "\phi'" , "\simeq" ' ] & A
\end{tikzcd}
\]
commutes.

\noindent (e) The given description is clear from the definitions of the other three arrows of the commutative diagram
\[
\begin{tikzcd}
S'(A)/Aut(A) \arrow[d, "\simeq" ' ] \arrow[r, "\text{part (d)}", "\simeq" '] &  S'_{k-1}(A)/Aut(A) \arrow[d, "\simeq" ]\\
S(A) \arrow[r, "\simeq" '] & S_{k-1}(A).
\end{tikzcd}
\]
\end{proof}

The following diagram summarizes our picture so far.
\begin{equation}\label{eq20}
\begin{tikzcd}
\{(X,\phi)\} \arrow[d, twoheadrightarrow] \arrow[r, "ext(-)"] & D_{k-1}(A) \arrow[d, twoheadrightarrow] \arrow[r, "\Theta_{k-1}"] & D_{k-2}(A)\arrow[d, twoheadrightarrow]  \arrow[r, "\Theta_{k-2}"] & ~ \cdots \arrow[r, "\Theta_2"] &  D_1(A)\arrow[d, twoheadrightarrow] & \hspace{-.45in} \cong  \prod\limits_r EXT(A_{r+1},A_r)\arrow[d, twoheadrightarrow] \\
S'(A) \arrow[d, twoheadrightarrow] \arrow[r, "\cong"] & S'_{k-1}(A) \arrow[d, twoheadrightarrow] \arrow[r, "\Theta_{k-1}"] & S'_{k-2}(A) \arrow[d, twoheadrightarrow] \arrow[r, "\Theta_{k-2}"] & ~ \cdots \arrow[r, "\Theta_2"] &  S'_1(A) \arrow[d, twoheadrightarrow] & \hspace{-.45in}  \cong  \prod\limits_r Ext^1(A_{r+1},A_r) \arrow[d, twoheadrightarrow] \\
S(A) \arrow[r, "\cong"] & S_{k-1}(A) \arrow[r, "\Theta_{k-1}"] & S_{k-2}(A) \arrow[r, "\Theta_{k-2}"] & ~ \cdots \arrow[r, "\Theta_2"] &  S_1(A) & \hspace{-.4in}  \cong  \displaystyle{\frac{\prod\limits_r Ext(A_{r+1},A_r)}{Aut(A)}}
\end{tikzcd}
\end{equation}
Here, $\{(X,\phi)\}$ means the collection of all pairs $(X, \phi)$ of consisting of an object $X$ of $\bT$ and an isomorphism $\phi: Gr^WX\rightarrow A$. The map $\{(X,\phi)\}\twoheadrightarrow S'(A)$ sends a pair to its $\sim'$-equivalence class. The map $S'(A)\twoheadrightarrow S(A)$ is induced by $(X,\phi)\mapsto X$. The maps $D_\ell(A)\twoheadrightarrow S'_\ell(A)$ and $S'_\ell(A)\twoheadrightarrow S_\ell(A)$, respectively, are given by modding out by $\sim'$ and (further) by $\sim$. All the maps between the middle and bottom rows can also be thought of as modding out by the action of $Aut(A)$. What remains of Theorem \ref{thm: main thm 1} to be established is the assertions about the structure of the fibers. This will be the subject of the rest of the section.

\subsection{Fibers of truncation maps I: Torsor structures}\label{sec: fibers 1}
Assume $2\leq \ell\leq k-1$. In this subsection we fix a generalized extension $(X_{m,n})_{n-m\leq \ell}$ of level $\ell-1$ and first describe the fiber of the truncation functor $\Theta_\ell: D_\ell(A)\rightarrow D_{\ell-1}(A)$ above it, i.e., the collection of all generalized extensions of level $\ell$ that become $(X_{m,n})_{n-m\leq \ell}$ once their lowest diagonal is erased. We then consider the equivalence relation $\sim'$ on this fiber.

Recall that given objects $X,Y,Z$ of $\bT$, $\sN\in EXT(Z,Y)$ and $\sL\in EXT(Y,X)$, the notation $EXTPAN(\sN,\sL)$ means the collection of all blended extensions of $\sN$ by $\sL$ (with no identification made). Recall also the notations $\sX_{m,n}^v$ and $\sX_{m,n}^h$ for extensions respectively coming from the arrows $X_{m,n-1}\hookrightarrow X_{m,n}$ and $X_{m,n}\twoheadrightarrow X_{m+1,n}$ of a generalized extension $(X_\dbullet)$ (i.e. extensions given by \eqref{eq11} and \eqref{eq10}, respectively, see Notation \ref{notation: horizontal and vertical extensions}).
\begin{construction}\label{cons: fibers to prod of EXTPANs}
There is a natural (to be seen to be bijective) map  
\begin{equation}\label{eq21}
\Theta_\ell^{-1}((X_{m,n})_{n-m\leq \ell}) = \left\{\begin{array}{l}
\text{fiber of} \\ D_\ell(A)\xrightarrow{\Theta_\ell} D_{\ell-1}(A) \\ \text{above $(X_{m,n})_{n-m\leq \ell}$} \end{array}\right\} \xrightarrow{ \ \  \ } \prod\limits_r EXTPAN(\sX^v_{r,r+\ell}, \sX^h_{r-1, r+\ell-1})
\end{equation}
(where the index $r$ on the right runs through the integers $1\leq r\leq k-\ell$) described as follows. Consider an element $(X_{m,n})_{n-m\leq \ell+1}$ of $D_\ell(A)$ in the fiber above $(X_{m,n})_{n-m\leq \ell}$. For each $X_{r-1, r+\ell}$ on its lowest diagonal, the morphisms
\[\begin{tikzcd}
X_{r-1, r+\ell-1} \arrow[d, hookrightarrow] \arrow[r, twoheadrightarrow] & X_{r, r+\ell-1} \arrow[d, hookrightarrow]\\
X_{r-1, r+\ell} \arrow[r, twoheadrightarrow] &  X_{r, r+\ell}
\end{tikzcd}\]
lead to a blended extension
\begin{equation} \label{eq23}
\begin{tikzcd}
   & & 0 \arrow{d} & 0 \arrow{d} &\\
   0 \arrow[r] & A_r \ar[equal]{d} \arrow[r, ] & X_{r-1, r+\ell-1}  \arrow[d, ] \arrow[r] &  X_{r, r+\ell-1} \arrow{d} \arrow[r] & 0 \\
   0 \arrow[r] & A_r \arrow[r] & X_{r-1, r+\ell} \arrow[d] \arrow[r] &  X_{r, r+\ell} \arrow{d}  \arrow[r] & 0, \\
   & & A_{r+\ell} \arrow{d} \ar[equal]{r} & A_{r+\ell} \arrow{d} & \\
   & & 0 & 0 &   
\end{tikzcd}
\end{equation}
in which every map is a composition (uniquely determined by the indices) of the structure arrows. The extensions on the top and right are respectively $\sX^h_{r-1,r+\ell-1}$ and $\sX_{r,r+\ell}^v$. The map \eqref{eq21} sends $(X_{m,n})_{n-m\leq \ell+1}$ to the tuple with this blended extension in its $r$-entry. 
\end{construction}
\begin{lemma}\label{lem: bijectivity of map from fibers at D level and prod of EXTPANs}
The map \eqref{eq21} is bijective.
\end{lemma}
\begin{proof}
We construct the inverse of \eqref{eq21}. Note that for each $r$ the two extensions $\sX^v_{r,r+\ell}$ and $\sX^h_{r-1, r+\ell-1}$ come from the data of the generalized extension $(X_{m,n})_{n-m\leq \ell}$ of level $\ell-1$. For each $r$, consider a blended extension of $\sX^v_{r,r+\ell}$ by $\sX^h_{r-1, r+\ell-1}$. It is given by a diagram of the form \eqref{eq23}, with the top and right extensions being $\sX^h_{r-1, r+\ell-1}$ and $\sX^v_{r,r+\ell}$, respectively. Now enlarge $(X_{m,n})_{n-m\leq \ell}$ to a generalized extension of level $\ell$ by adding to its data, for each $r$, the object $X_{r-1,r+\ell}$ in entry $(r-1,r+\ell)$ and the morphisms $X_{r-1,r+\ell-1}\hookrightarrow X_{r-1,r+\ell}$ and $X_{r-1,r+\ell}\twoheadrightarrow X_{r,r+\ell}$ of the corresponding blended extension. The augmented data $(X_{m,n})_{n-m\leq \ell+1}$ is a generalized extension of level $\ell$. Indeed, the only new squares formed by the structure arrows are the ones in the top rights of our blended extensions. So axiom (i) of the definition of a generalized extension holds. As for axiom (ii) (the exactness of the sequences \eqref{eq11}), the new sequences we must consider are exactly the sequences in the middle columns of our blended extensions.

By sending the tuple of blended extensions we started with to the generalized extension $(X_{m,n})_{n-m\leq \ell+1}$ we obtain a map
\begin{equation}\label{eq22}
\prod\limits_r EXTPAN(\sX^v_{r,r+\ell}, \sX^h_{r-1, r+\ell-1}) \xrightarrow{ \ \ \ } 
\Theta_\ell^{-1}((X_{m,n})_{n-m\leq \ell}).
\end{equation}
The reader easily sees that this maps is the inverse to \eqref{eq21}.
\end{proof}
It is convenient to have a notation for blended extensions of the form \eqref{eq23}:
\begin{notation} For a generalized extension $(X_\dbullet)$ of level $\ell \geq 2$, for any $r$ we denote the blended extension \eqref{eq23} with middle object $X_{r-1, r+\ell}$ by $\sX_{r-1,r+\ell}$ (without a superscript $h$ or $v$).
\end{notation}
Recall from \S \ref{sec: background on blended extensions} that given blended extensions $\sX$ and $\sX'$ of an extension $\sN$ by an extension $\sL$, a morphism (automatically an isomorphism) of blended extensions from $\sX$ to $\sX'$ is a morphism from the middle object $X$ of $\sX$ to the middle object $X'$ of $\sX'$ that induces identity maps on objects in $\sL$ and $\sN$; that is, a morphism $X\rightarrow X'$ that together with the identity maps on the objects in $\sL$ and $\sN$ commute with the arrows in $\sX$ and $\sX'$ (in the obvious sense). The collection of isomorphism classes of blended extensions of $\sN$ by $\sL$ is denoted by $Extpan(\sN,\sL)$.

\begin{lemma}\label{lem: when are two gen ext in a fiber sim' equivalent}
(a) Two elements of the fiber of $\Theta_\ell: D_\ell(A)\rightarrow D_{\ell-1}(A)$ above $(X_\dbullet)$ are $\sim'$-equivalent if and only if their images under the map \eqref{eq21} coincide in 
\[
\prod\limits_r Extpan(\sX^v_{r,r+\ell}, \sX^h_{r-1, r+\ell-1}).
\]
\noindent (b) The map \eqref{eq21} descends to a bijection
\[
\bigm(\Theta_\ell^{-1}((X_\dbullet))\bigm)\bigm/ \sim'  \ \xrightarrow{ \ \simeq \ }  \prod\limits_r Extpan(\sX^v_{r,r+\ell}, \sX^h_{r-1, r+\ell-1}).
\]
In particular,
\[ \bigm(\Theta_\ell^{-1}((X_\dbullet))\bigm)\bigm/ \sim' \]
is either empty or a torsor over
\[
\prod\limits_r Ext^1(A_{r+\ell}, A_r).
\]
Moreover, it is nonempty if and only if for each $r$, the image of the Yoneda product of $\sX^v_{r,r+\ell}$ and $\sX^h_{r-1, r+\ell-1}$ in $Ext^2(A_{r+\ell}, A_r)$ vanishes.
\end{lemma}
\begin{proof}
We first note that part (b) follows immediately from part (a), Lemma \ref{lem: bijectivity of map from fibers at D level and prod of EXTPANs}, and the general theory of blended extensions (see \S \ref{sec: background on blended extensions}, in particular, Lemma \ref{lem: criteria for compatibility of extension pairs}(a)). So we will focus on part (a).

Suppose $(Y_\dbullet)$ and $(Y'_\dbullet)$ are in the fiber above $(X_\dbullet)$, so that dropping the lowest diagonals, $(Y_\dbullet)$ and $(Y'_\dbullet)$ are just $(X_\dbullet)$. The blended extensions \eqref{eq23} for $(Y_\dbullet)$ and $(Y'_\dbullet)$ (respectively, denoted by $\sY _{r-1,r+\ell}$ and $\sY'_{r-1,r+\ell}$) have the same top rows and the same right columns, coming from $(X_\db)$. Since every arrow in $\sY _{r-1,r+\ell}$ (resp. $\sY'_{r-1,r+\ell}$) is the appropriate composition of the structure arrows of $(Y_\dbullet)$ (resp. $(Y'_\dbullet)$), every morphism $(f_\dbullet):  (Y_\dbullet)\rightarrow (Y'_\dbullet)$ of generalized extensions includes the data of a collection of maps between the corresponding objects of $\sY _{r-1,r+\ell}$ and $\sY'_{r-1,r+\ell}$ that commute with the arrows in the two blended extensions.

Let $(f_\dbullet):  (Y_\dbullet)\rightarrow (Y'_\dbullet)$ be a morphism that is identity on $A$. Both $(Y_\dbullet)$ and $(Y'_\dbullet)$ truncate to $(X_\dbullet)$, so that $(f_\dbullet)$ must be identity on $(X_\dbullet)$ (see Lemma \ref{lem: if identity on A identity everywhere}). Combining with the earlier comments, it follows that for each $r$, the isomorphism $f_{r-1,r+\ell}$ from $Y_{r-1,r+\ell}$ to $Y'_{r-1,r+\ell}$ gives an isomorphism of blended extensions from $\sY _{r-1,r+\ell}$ to $\sY'_{r-1,r+\ell}$.

Conversely, suppose that for each $r$, the classes of blended extensions $\sY _{r-1,r+\ell}$ and $\sY'_{r-1,r+\ell}$ coincide in $Extpan(\sX^v_{r,r+\ell}, \sX^h_{r-1, r+\ell-1})$. Let $f_{r-1,r+\ell}$ be the morphism $Y_{r-1,r+\ell}\rightarrow Y'_{r-1,r+\ell}$ that gives an isomorphism of blended extensions $\sY _{r-1,r+\ell}\rightarrow \sY'_{r-1,r+\ell}$. Then $f_{r-1,r+\ell}$ induces identity on $X_{r-1,r+\ell-1}$ and $X_{r,r+\ell}$. By Lemma \ref{lem: gluing morphisms between objects on the lowest diagonal} the collection of morphisms $f_{r-1,r+\ell}$ glues together to give a morphism $(f_\dbullet): (Y_\dbullet)\rightarrow (Y'_\dbullet)$. This morphism is identity on the diagonal just above the lowest, and hence is identity on all of $(X_\dbullet)$. 
\end{proof}

\subsection{Fibers of truncation maps II: Canonicity of the torsor structures}\label{sec: fibers of truncations II}
We continue to assume $2\leq \ell\leq k-1$. Let $(X_\db)$ be a generalized extension of level $\ell-1$. Denote the class of $(X_\db)$ in $S'_{\ell-1}(A)$ by $[(X_\db)]_{\sim'}$. In the previous subsection we studied $\Theta_\ell^{-1}((X_\db))/\sim'$ and saw that it has a torsor structure. To finish the proof of Theorem \ref{thm: main thm 1}(b,c) we need to show that 
\[
\Theta_\ell^{-1}([(X_\db)]_{\sim'}) \cong \Theta_\ell^{-1}((X_\db))/\sim', 
\]
and that moreover the torsor structure obtained this way on the fiber of $S'_\ell(A)\rightarrow S'_{\ell-1}(A)$ above $\epsilon:= [(X_\db)]_{\sim'}$ is independent of the choice of representative $(X_\db)$ for $\epsilon$. By the end of this subsection these will be established and the proof of Theorem \ref{thm: main thm 1}(b,c) will be completed. 

The following simple definition is convenient.
\begin{defn}\label{def: transports}
Let $(X_\dbullet)$ be a generalized extension of any level. Let $(i,j)$ be an eligible pair, and $f:X_{i,j}\rightarrow X'$ an isomorphism.

\noindent (a) The transport of $(X_\dbullet)$ along $f$, denoted by $tr((X_\db), f)$, is the generalized extension obtained from $(X_\dbullet)$ by replacing the object $X_{i,j}$ at entry $(i,j)$ by $X'$ via $f$. That is, by replacing $X_{i,j}$ by $X'$, and the arrows to (resp. from) $X_{i,j}$ by their composition with $f$ (resp. $f^{-1}$).

\noindent (b) The collection of morphisms $(f_\db)$ given by $f_{m,n}=Id_{X_{m,n}}$ for every eligible pair $(m,n)\neq (i,j)$ and $f_{i,j}=f$ is called the isomorphism given by the transport datum $f$. 
\end{defn}

One easily sees that the transport defined above is indeed a generalized extension\footnote{We note that in the generality of Definition \ref{def: transports}, the transport may not be a generalized extension {\it of $A$}. However, in all applications of this construction in the paper, whenever $i=j-1$ (so that $X_{i,j}=A_j$) we will also take $X'$ to be $A_j$, so that the transport will always be indeed a generalized extension of $A$ as well.}, and that the collection of isomorphisms $(f_\db)$ of (b) commutes with the structure arrows of $(X_\db)$ and its transport $tr((X_\db), f)$. One also easily sees that the transport construction behaves well with respect to compositions: given $(X_\db)$ and isomorphisms $X_{i,j}\xrightarrow{f}Y \xrightarrow{g} Z$, we have
\begin{equation}\label{eq109}
tr(tr((X_\db),f),g) = tr((X_\db),gf)
\end{equation}
and the composition of the isomorphisms given by the transport datum $f$ first and then $g$
\[
(X_\db) \rightarrow tr((X_\db),f) \rightarrow tr(tr((X_\db),f),g)
\]
is just the isomorphism given by the transport datum $gf$.

More generally, given a generalized extension $(X_\db)$ of any level, a set $I$ of eligible pairs of indices and for each $(m,n)\in I$ an isomorphism $f_{m,n}:X_{m,n}\rightarrow X'_{m,n}$, we may talk about the transport $tr((X_\db),(f_{m,n})_{(m,n)\in I})$ of $(X_\dbullet)$ along $(f_{m,n})_{(m,n)\in I}$. Making the transport all at once is the same as making it step by step for one $(m,n)\in I$ at a time (note that the transport operations along morphisms from $X_{m,n}$ for different $(m,n)$ commute with one another). The collection of morphisms $(g_\db)$ where $g_{m,n}=Id_{X_{m,n}}$ for every eligible $(m,n)\notin I$ and $g_{m,n}=f_{m,n}$ if $(m,n)\in I$ commutes with the structure arrows of $(X_\db)$ and $tr((X_\db),(f_{m,n})_{(m,n)\in I})$; in line with Definition \ref{def: transports}(b), we call $(g_\db)$ the isomorphism given by the transport data $(f_{m,n})_{(m,n)\in I}$.

Now let $(X_\db)$ and $(X'_\db)$ be generalized extensions of level $\ell-1$ of $A$, and $(f_\db):(X_\db)\rightarrow (X'_\db)$ an isomorphism of generalized extensions. Given $(Y_\db)\in \Theta_\ell^{-1}((X_\db))$, it follows from the commutativity of the diagrams of \eqref{eq diagrams for morphisms of gen exts} that the transport of $(Y_\db)$ along $(f_\db)$ is in the fiber of $\Theta_\ell$ above $(X'_\db)$. The map
\begin{equation}\label{eq24}
\Theta_\ell^{-1}((X_\db)) \rightarrow \Theta_\ell^{-1}((X'_\db)) \hspace{.3in}(Y_\db)\mapsto tr((Y_\db), (f_\db))
\end{equation}
is a bijection, with its inverse given by transport along $(f^{-1}_\db)$. For every $(Y_\db)$, the isomorphism given by the transport data $(f_\db)$ is an isomorphism of generalized extensions $(Y_\db)\rightarrow tr((Y_\db), (f_\db))$. Thus \eqref{eq24} descends to a bijection between the $\sim$-equivalence classes. 

If $(Y_\db^{(1)})$ and $(Y_\db^{(2)})$ in $\Theta_\ell^{-1}((X_\db))$ are $\sim'$-equivalent with $(g_\db): (Y_\db^{(1)})\rightarrow (Y_\db^{(2)})$ an isomorphism that is identity on $A$, then the composition 
\[\begin{tikzcd}[column sep = large]
tr((Y_\db^{(1)}), (f_\db)) \ \ \arrow[r, "\text{iso. given by}", "\text{tr. data $(f_\db^{-1})$}" '] & \ \ (Y_\db^{(1)}) \arrow[r, "(g_\db)"] & (Y_\db^{(2)}) \ \  \arrow[r, "\text{iso. given by}" , "\text{tr. data $(f_\db)$}" '] & \ \ tr((Y_\db^{(2)}), (f_\db))
\end{tikzcd}
\] 
is identity on $A$ as well (even if $(f_\db)$ is not identity on $A$). Hence \eqref{eq24} also descends to a bijection between the $\sim'$-equivalence classes, with its inverse induced by transport along $(f^{-1}_\db)$. 

If the morphism $(f_\db):(X_\db)\rightarrow (X'_\db)$ is identity on $A$, then every $(Y_\db)$ above $(X_\db)$ is $\sim'$-equivalent to its transport along $(f_\db)$, with the isomorphism given by the transport data giving the $\sim'$-equivalence. 

In particular, we note from the above that if $(X_\db)$ and $(X'_\db)$ in $D_{\ell-1}(A)$ are $\sim$-equivalent (resp. $\sim'$-equivalent), then for every $(Y_\db)\in \Theta^{-1}((X_\db))$ there exists $(Y'_\db)\in \Theta^{-1}((X'_\db))$ that is $\sim$-equivalent (resp. $\sim'$-equivalent) to $(Y_\db)$.

We obtain the following lemma regarding the fibers of truncation maps $S'_\ell(A)\rightarrow S'_{\ell-1}(A)$ and $S_\ell(A)\rightarrow S_{\ell-1}(A)$ (recall that we refer to both of these also by $\Theta_\ell$).
\begin{lemma}\label{lem: fibers of truncations between oDs are projections of fibers of truncations between Ds}
Let $(X_\db)\in D_{\ell-1}(A)$. Denote the classes of $(X_\db)$ in $S_\ell(A)$ and $S'_\ell(A)$ respectively by $[(X_\db)]_\sim$ and $[(X_\db)]_{\sim'}$.

\noindent (a) The natural injection 
\[
\bigm(\Theta_\ell^{-1}((X_\dbullet))\bigm)\bigm/ \sim \ \rightarrow \ \Theta_\ell^{-1}([(X_\db)]_\sim)
\]
is bijective. If $(X'_\db)\in D_{\ell-1}(A)$ and $(f_\db): (X_\db)\rightarrow (X'_\db)$ is an isomorphism, then we have a commutative diagram
\[
\begin{tikzcd}[column sep = 0in]
\bigm(\Theta_\ell^{-1}((X_\dbullet))\bigm)\bigm/ \sim \arrow[rr, "\text{tr. along $(f_\db)$}" , "\simeq" ' ] \arrow[dr, "\simeq"] & & \bigm(\Theta_\ell^{-1}((X'_\dbullet))\bigm)\bigm/ \sim \arrow[dl, "\simeq"] \\
& \Theta_\ell^{-1}([(X_\db)]_\sim) &
\end{tikzcd}
\]
where the horizontal arrow is given by transport along $(f_\db)$ (descending from \eqref{eq24}) and the other two arrows are the natural maps.

\noindent (b) The natural injection 
\begin{equation}\label{eq25}
\bigm(\Theta_\ell^{-1}((X_\dbullet))\bigm)\bigm/ \sim' \ \rightarrow \ \Theta_\ell^{-1}([(X_\db)]_{\sim'})
\end{equation}
is bijective. If $(X'_\db)\in D_{\ell-1}(A)$ and $(f_\db): (X_\db)\rightarrow (X'_\db)$ is an isomorphism that is identity on $A$, then we have a commutative diagram
\begin{equation}\label{eq106}
\begin{tikzcd}[column sep = 0in]
\bigm(\Theta_\ell^{-1}((X_\dbullet))\bigm)\bigm/ \sim' \arrow[rr, "\text{tr. along $(f_\db)$}" , "\simeq" ' ] \arrow[dr, "\simeq" '] & & \bigm(\Theta_\ell^{-1}((X'_\dbullet))\bigm)\bigm/ \sim' \arrow[dl, "\simeq"] \\
& \Theta_\ell^{-1}([(X_\db)]_{\sim'}) &
\end{tikzcd}
\end{equation}
where the horizontal arrow is given by transport along $(f_\db)$ (descending from \eqref{eq24}) and the other two arrows are the natural maps. If $(f_\db): (X_\db)\rightarrow (X'_\db)$ is an isomorphism that is not necessarily identity on $A$, we still have a bijection
\begin{equation}\label{eq26}
\bigm(\Theta_\ell^{-1}((X_\dbullet))\bigm) \bigm/ \sim'  \ \xrightarrow{\text{tr. along $(f_\db)$}} \ \bigm(\Theta_\ell^{-1}((X'_\dbullet))\bigm)\bigm/ \sim'
\end{equation}
which forms a commutative diagram if we pass along on both sides to $\Theta_\ell^{-1}([(X_\db)]_\sim)$.
\end{lemma}
Combining part (b) with Lemma \ref{lem: when are two gen ext in a fiber sim' equivalent}(b) we obtain torsor structures on the nonempty fibers of $S'_\ell(A)\rightarrow S'_{\ell-1}(A)$. At the moment however, the torsor structure on the fiber of $S'_\ell(A)\rightarrow S'_{\ell-1}(A)$ above the class of $(X_\db)$ appears to depend on the choice of representative $(X_\db)$. Our next task is to rule out this dependence.

The group $Aut(A)$ acts on
\begin{equation}\label{eq27}
\prod\limits_r Ext^1(A_{r+\ell}, A_r)
\end{equation}
similarly to the action we already considered when $\ell=1$, i.e. by pushforwards and pullbacks. An element $\sigma=(\sigma_r)$ (with $\sigma_r\in Aut(A_r)$) sends a tuple of extension classes $\mathcal{E}=(\mathcal{E}_r)$ to the tuple that has $(\sigma_r)_\ast (\sigma_{r+\ell}^{-1})^\ast \mathcal{E}_r$ in its $r$-entry. Denoting the image of $\mathcal{E}$ under the action by $\sigma$ by $\sigma\cdot \mathcal{E}$, thus the $r$-entry of
$\sigma\cdot \mathcal{E}$ is obtained, after taking a representative for $\mathcal{E}_r$ in $EXT(A_{r+\ell}, A_r)$, by composing the arrow coming out of $A_r$ by $\sigma_r^{-1}$ and the arrow going to $A_{r+\ell}$ by $\sigma_{r+\ell}$.

\begin{lemma}\label{lem: compatibility of torsor structures}
Let $(X_\db)$ and $(X'_\db)$ be in $D_{\ell-1}(A)$ and $(f_\db):(X_\db)\rightarrow (X'_\db)$ an isomorphism. Suppose that the fiber of $\Theta_\ell: D_\ell(A)\rightarrow D_{\ell-1}(A)$ above $(X_\db)$ (and hence $(X'_\db)$) is nonempty. Consider 
\[
\bigm(\Theta_\ell^{-1}((X_\dbullet))\bigm) \bigm/ \sim' \hspace{.2in}\text{and} \hspace{.2in} \bigm(\Theta_\ell^{-1}((X'_\dbullet))\bigm)\bigm/ \sim' 
\]
as torsors for \eqref{eq27} via the canonical bijection of Lemma \ref{lem: when are two gen ext in a fiber sim' equivalent}(b) (for $(X_\db)$ and $(X'_\db)$, respectively). Then the bijection \eqref{eq26} satisfies the following identity: denoting the action of the group \eqref{eq27} on the torsors above by $\ast$, the map \eqref{eq26} (as well as \eqref{eq24}) by $tr(-,(f_\db))$, and the restriction of $(f_\db)$ to $A$ by $f_A$, then for every $(Y_\db) \in \Theta_\ell^{-1}((X_\dbullet))$ and for every tuple of extension classes $\mathcal{E}=(\mathcal{E}_r)$ in \eqref{eq27} we have
\begin{equation}\label{eq29}
tr(\mathcal{E}\ast [(Y_\db)]_{\sim'}, \, (f_\db))  \ = \  (f_A\cdot \mathcal{E}) \ast \, tr([(Y_\db)]_{\sim'}, (f_\db)),
\end{equation}
where (with abuse of notation) $[(Y_\db)]_{\sim'}$ here means the image of $(Y_\db)$ in $(\Theta_\ell^{-1}((X_\db)))/\sim'$. 

In particular, if the restriction of $(f_\db)$ to $A$ is a scalar multiple of the identity map, then the bijection \eqref{eq26} is an isomorphism of torsors.
\end{lemma}
\begin{proof}
Consider the commutative diagram
\begin{equation}\label{eq107}
\begin{tikzcd}[column sep = huge]
\Theta_\ell^{-1}((X_\dbullet)) \arrow[r, "\eqref{eq21}", "\simeq" '] \arrow[d, twoheadrightarrow] & \prod\limits_r EXTPAN(\sX^v_{r,r+\ell}, \sX^h_{r-1, r+\ell-1}) \arrow[d, twoheadrightarrow] \\
\bigm(\Theta_\ell^{-1}((X_\dbullet))\bigm) \bigm/ \sim' \arrow[r, "\text{Lem. \ref{lem: when are two gen ext in a fiber sim' equivalent}(b)}", "\simeq" '] & \prod\limits_r Extpan(\sX^v_{r,r+\ell}, \sX^h_{r-1, r+\ell-1}).
\end{tikzcd}
\end{equation}
The torsor structures over the group \eqref{eq27} on the lower level descend from a map
\[
\prod\limits_r EXT(A_{r+\ell},A_r)\times \prod\limits_r EXTPAN(\sX^v_{r,r+\ell}, \sX^h_{r-1, r+\ell-1})\rightarrow \prod\limits_r EXTPAN(\sX^v_{r,r+\ell}, \sX^h_{r-1, r+\ell-1}).
\]
We use the symbol $\ast$ for the operation given by the latter map on the top right object of \eqref{eq107}, as well as the operation on $\Theta_\ell^{-1}((X_\dbullet))$ induced by it, and the operations on the lower level of the diagram descended from it. We also use the same notation for the analogous operations for $(X'_\db)$.

Let $(Y_\db)$ be in $\Theta_\ell^{-1}((X_\dbullet))$ and $\sE=(\sE_r)$ a tuple in \eqref{eq27}. Lift each $\sE_r$ to an element of $EXT(A_{r+\ell},A_r)$, which with abuse of notation we also denote by $\sE_r$. 

Set
\[(Z_\db) := \sE\ast (Y_\db) \in \Theta_\ell^{-1}((X_\dbullet))\]
and $(Z'_\db)=tr((Z_\db), (f_\db))$. Suppose that the blended extension of $\sX_{r,r+\ell}^v$ by $\sX_{r-1,r+\ell-1}^h$ associated to $(Z_\db)$ by Construction \ref{cons: fibers to prod of EXTPANs} is given by the diagram on the left below. Then the blended extension of ${\sX'}_{r,r+\ell}^v$ by ${\sX'}_{r-1,r+\ell-1}^h$ associated to $(Z'_\db)$  is given by the diagram on the right below. We have dropped the indices from the $f_{m,n}$ to save space (they are determined by the indices of the objects). 
\begin{equation}\label{eq35}
\begin{tikzcd}[row sep=small, column sep=small]
   & & 0 \arrow{d} & 0 \arrow{d} &\\
   0 \arrow[r] & A_r \ar[equal]{d} \arrow[r, "j"] & X_{r-1, r+\ell-1}  \arrow[d, "\overline{\iota}" ] \arrow[r, "\pi"] &  X_{r, r+\ell-1} \arrow[d, "\iota"] \arrow[r] & 0 \\
   0 \arrow[r] & A_r \arrow[r, "\overline j"] & Z_{r-1, r+\ell} \arrow[r, "\overline \pi"] \arrow[d, "\overline \omega"] &  X_{r, r+\ell} \arrow[d, "\omega"]  \arrow[r] & 0 \\
   & & A_{r+\ell} \arrow{d} \ar[equal]{r} & A_{r+\ell} \arrow{d} & \\
   & & 0 & 0 &   
\end{tikzcd} 
\hspace{.1in}
\begin{tikzcd}[row sep=small, column sep=small]
    & & 0 \arrow{d} & 0 \arrow{d} &\\
   0 \arrow[r] & A_r \ar[equal]{d} \arrow[r, "fjf^{-1}"] & X'_{r-1, r+\ell-1}  \arrow[d, "\overline{\iota} f^{-1}" ] \arrow[r, "f\pi f^{-1}"] &  X'_{r, r+\ell-1} \arrow[d, "f \iota f^{-1}"] \arrow[r] & 0 \\
   0 \arrow[r] & A_r \arrow[r, "\overline j f^{-1}"] & Z_{r-1, r+\ell} \arrow[r, "f\overline \pi"] \arrow[d, "f\overline \omega"] &  X'_{r, r+\ell} \arrow[d, "f\omega f^{-1}"]  \arrow[r] & 0 \\
   & & A_{r+\ell} \arrow{d} \ar[equal]{r} & A_{r+\ell} \arrow{d} & \\
   & & 0 & 0 &  
\end{tikzcd} 
\end{equation}
Bringing the indices back to avoid confusion, comparing the second rows of the two diagrams we have
\begin{equation}\label{eq28}
(f_{r-1,r})_\ast  \sZ_{r-1,r+\ell}^h = f_{r,r+\ell}^\ast {\sZ'}_{r-1,r+\ell}^h,
\end{equation}
where $\sZ_{r-1,r+\ell}^h$ and ${\sZ'}_{r-1,r+\ell}^h$ refer to the second horizontal extensions in the two diagrams respectively (see Notation \ref{notation: horizontal and vertical extensions}). The equality here as well as all the other equalities in the rest of this argument take place in the corresponding $Ext^1$ groups.

By definition of the torsor structure on the set of isomorphism classes of blended extensions (see \S \ref{sec: background on blended extensions}), we have 
\[
\sZ_{r-1,r+\ell}^h = \sY_{r-1,r+\ell}^h+\omega^\ast\sE_r
\]
in $Ext^1(X_{r,r+\ell}, A_r)$. Combining the last two equations we obtain
\[
f_{r,r+\ell}^\ast {\sZ'}_{r-1,r+\ell}^h = (f_{r-1,r})_\ast\sY_{r-1,r+\ell}^h + (f_{r-1,r})_\ast \, \omega^\ast\sE_r = f_{r,r+\ell}^\ast {\sY'}_{r-1,r+\ell}^h + (f_{r-1,r})_\ast \, \omega^\ast\sE_r,
\]
where $(Y'_\db)=tr((Y_\db),(f_\db))$. (We have used the analogue of \eqref{eq28} for $(Y_\db)$.) Thus
\begin{equation}\label{eq108}
{\sZ'}_{r-1,r+\ell}^h = {\sY'}_{r-1,r+\ell}^h  + (\omega f^{-1}_{r,r+\ell})^\ast (f_{r-1,r})_\ast \sE_r
\end{equation}
as pushforward commutes with pullback. 

Let 
\[
(Z''_\db) :=  (f_A\cdot \mathcal{E}) \ast \, tr((Y_\db), (f_\db)) = (f_A\cdot \mathcal{E}) \ast (Y'_\db),
\]
so that the right hand side of \eqref{eq29} is the class of $(Z''_\db)$ in 
\[
\bigm(\Theta_\ell^{-1}((X'_\dbullet))\bigm)\bigm/ \sim'. 
\]
On recalling the definition of $f_A\cdot \sE$ and the fact that the map $X'_{r,r+\ell}\twoheadrightarrow A_{r+\ell}$ is $f_{r+\ell-1, r+\ell}\omega f_{r,r+\ell}^{-1}$, we have
\[
{\sZ''}_{r-1,r+\ell}^h = {\sY'}_{r-1,r+\ell}^h + (f_{r+\ell-1, r+\ell} \, \omega f_{r,r+\ell}^{-1})^\ast (f_{r-1,r})_\ast (f^{-1}_{r+\ell-1,r+\ell})^\ast \, \sE_r \stackrel{\eqref{eq108}}{=} {\sZ'}_{r-1,r+\ell}^h. 
\] 
Thus by Lemma \ref{lem: X mapsto X^h is injective when Hom(A_2,A_1)=0}, $(Z'_\db)$ and $(Z''_\db)$ coincide in  
\[
\prod\limits_r Extpan({\sX'}^v_{r,r+\ell}, {\sX'}^h_{r-1, r+\ell-1})
\]
and hence are in the same $\sim'$-equivalence class. This completes the proof of the identity. 

The statement about the special case when $(f_\db)$ is a scalar automorphism of $A$ is immediate from the identity.
\end{proof}
Putting the results of this subsection and the previous one together, we obtain parts (b) and (c) of Theorem \ref{thm: main thm 1}:

\begin{proof}[Proof of Theorem \ref{thm: main thm 1}(b,c)]
Consider an element $\epsilon$ of $S'_{\ell-1}(A)$, i.e. a $\sim'$-equivalence class of generalized extensions of level $\ell-1$. Let $(X_\db)$ be a representative of the class. Use the bijection \eqref{eq25} and the bijection of Lemma \ref{lem: when are two gen ext in a fiber sim' equivalent}(b) to give the fiber of $S'_\ell(A)\rightarrow S'_{\ell-1}(A)$ above $\epsilon$ the structure of a torsor for 
\[
\prod\limits_r Ext^1(A_{r+\ell}, A_r)
\]
when this fiber is nonempty. This torsor structure is independent of the choice of $(X_\db)$. Indeed, suppose one chooses another representative $(X'_\db)$ for $\epsilon$. Let $(f_\db):(X_\db)\rightarrow (X'_\db)$ be an isomorphism that is identity on $A$. By Lemma \ref{lem: compatibility of torsor structures}, transport along $(f_\db)$ gives an isomorphism of torsors \eqref{eq26}. In view of the commutative diagram \eqref{eq106} of Lemma \ref{lem: fibers of truncations between oDs are projections of fibers of truncations between Ds}(b), the induced torsor structures on $\Theta_\ell^{-1}(\epsilon)$ thus coincide.

The assertion in Theorem \ref{thm: main thm 1}(c) giving a sufficient condition for surjectivity of the truncation map $S'_\ell(A)\rightarrow S'_{\ell-1}(A)$ is immediate from the constructions and Lemma \ref{lem: criteria for compatibility of extension pairs}(a). (In fact, a more precise statement about the image of $S'_\ell(A)\rightarrow S'_{\ell-1}(A)$ can be obtained from the last sentence of Lemma \ref{lem: when are two gen ext in a fiber sim' equivalent}(b).)
\end{proof}

\subsection{Fibers of truncation maps III: $\Gamma$-actions}\label{sec: fibers 3}
So far, we have studied the fibers of the truncation maps $S'_\ell\rightarrow S'_{\ell-1}$. The proof of the final remaining part of Theorem \ref{thm: main thm 1} involves an additional ingredient, namely the group actions that allow us to pass on from the fibers on the second row of equation \eqref{eq20} to the fibers on its third row. This is the subject of this subsection. In particular, we will deduce part (e) of Theorem \ref{thm: main thm 1}.

Fix $2\leq \ell\leq k-1$ and $(X_\db)\in D_{\ell-1}(A)$. The group $Aut((X_\db))$ of the automorphisms of $(X_\db)$ (as a generalized extension) acts on the fiber of $\Theta_\ell: D_\ell(A)\rightarrow  D_{\ell-1}(A)$ above $(X_\db)$ by transport. That is, an element $\sigma=(\sigma_\db) \in Aut((X_\db))$ acts by sending $(Y_\db)\in \Theta_\ell^{-1}((X_\db))$ to the transport of $(Y_\db)$ along $\sigma$. The fact that this transport is also in $\Theta_\ell^{-1}((X_\db))$ is because $\sigma$ is an automorphism of generalized extensions. Denote the image of $(Y_\db)$ under this action by $\sigma\cdot (Y_\db)$. Then $\sigma\cdot (Y_\db) = tr((Y_\db), \sigma)$ is obtained from $(Y_\db)$ by twisting {\it only} the arrows between the two lowest diagonals of $(Y_\db)$. More explicitly, each horizontal arrow $Y_{r-1, r+\ell}\twoheadrightarrow X_{r, r+\ell}$ (resp. vertical arrow $X_{r-1, r+\ell-1}\hookrightarrow Y_{r-1, r+\ell}$) between the lowest two diagonals gets composed with $\sigma_{r, r+\ell}$ (resp. $\sigma_{r-1, r+\ell-1}^{-1}$). The rest of the diagram remains unchanged. Note that since the action of $Aut((X_\db))$ on $\Theta_\ell^{-1}((X_\db))$ is given by transports, we have $(Y_\db)\sim \sigma\cdot (Y_\db)$ for all $(Y_\db)$ and $\sigma$.

In view of Lemma \ref{lem: fibers of truncations between oDs are projections of fibers of truncations between Ds}(b) (applied with $(X'_\db)=(X_\db)$ and $(f_\db)=\sigma$), the action of $Aut((X_\db))$ on the fiber $\Theta_\ell^{-1}((X_\db))$ descends to an action on $\bigm( \Theta_\ell^{-1}((X_\db)) \bigm) \bigm/ \sim'$. This action captures the passing from $\sim'$ to $\sim$ on $\Theta_\ell^{-1}((X_\db))$:
\begin{lemma}\label{lem: orbits of Gamma action}
As above, let $2\leq \ell\leq k-1$ and $(X_\db)\in D_{\ell-1}(A)$. Let $(Y_\db)$ and $(Y'_\db)$ be in $\Theta_\ell^{-1}((X_\db))$. Then $(Y_\db)$ and $(Y'_\db)$ are $\sim$-equivalent if and only if the classes of $(Y_\db)$ and $(Y'_\db)$ mod $\sim'$ are in the same orbit of the action of $Aut((X_\db))$ on  
\begin{equation}\label{eq30}
\bigm( \Theta_\ell^{-1}((X_\db)) \bigm) \bigm/ \sim'.
\end{equation}
\end{lemma}
\begin{proof} 
Suppose that $(Y_\db)$ and $(Y'_\db)$ are isomorphic. Let $(f_\db): (Y_\db)\rightarrow (Y'_\db)$ be an isomorphism. Then $(f_\db)$ restricts to an automorphism $\sigma$ of $(X_\db)$. Consider the composition of isomorphisms
\[
\begin{tikzcd}[row sep=large]
\sigma\cdot (Y_\db) \arrow[r] & (Y_\db) \arrow[r, "(f_\db)"] &  (Y'_\db),
\end{tikzcd}
\]
where the first arrow is the isomorphism given by the transport data $\sigma^{-1}$ (so is identity on the objects of the lowest diagonal and $\sigma^{-1}$ on $(X_\db)$, see \S \ref{sec: fibers of truncations II}). This composition isomorphism is identity on $A$, so that $\sigma\cdot (Y_\db) \sim' (Y'_\db)$.
 
Conversely, suppose that the classes of $(Y_\db)$ and $(Y'_\db)$ in \eqref{eq30} are in the same orbit of the action of $Aut((X_\db))$. Thus there is $\sigma\in Aut((X_\db))$ such that $\sigma\cdot (Y_\db) \sim' (Y'_\db)$. Then $(Y_\db)\sim \sigma\cdot (Y_\db)\sim (Y'_\db)$.
\end{proof}

The group actions above on the fibers of $D_\ell(A)\rightarrow D_{\ell-1}(A)$ descend canonically to group actions on the fibers of $S'_\ell(A)\rightarrow S'_{\ell-1}(A)$. Indeed, let $\epsilon'\in S'_{\ell-1}(A)$. For every $(X_\db)$ and $(X'_\db)$ in $D_{\ell-1}(A)$ that belong to $\epsilon'$, by Lemma \ref{lem: if identity on A identity everywhere} there exists a unique isomorphism $(f_\db):(X_\db)\rightarrow (X'_\db)$ that is identity on $A$. This isomorphism induces an isomorphism $Aut((X_\db))\rightarrow Aut((X'_\db))$ by conjugation. These distinguished isomorphisms $Aut((X_\db))\rightarrow Aut((X'_\db))$ for various pairs $((X_\db),(X'_\db))$ of representatives of $\epsilon'$ are compatible with respect to compositions with one another, so that they form a projective system of groups. Set
\[
\Gamma(\epsilon') := \lim\limits_{(X_\db)\in\epsilon'} Aut((X_\db)).
\]
An element of $\Gamma(\epsilon')$ is the data of an automorphism of $(X_\db)$ for each representative $(X_\db)\in\epsilon'$, such that the automorphisms for different representatives corresponds to one another under the distinguished isomorphisms between the automorphism groups. For each $(X_\db)$, the projection map $\Gamma(\epsilon')\rightarrow Aut((X_\db))$ is an isomorphism.

Given any $(X_\db)$ and $(X'_\db)$ representing $\epsilon'$, we have a commutative diagram
\[
\begin{tikzcd}[column sep = 0in]
\bigm(\Theta_\ell^{-1}((X_\dbullet))\bigm)\bigm/ \sim' \arrow[rr, "\simeq"  ] \arrow[dr, "\simeq" '] & & \bigm(\Theta_\ell^{-1}((X'_\dbullet))\bigm)\bigm/ \sim' \arrow[dl, "\simeq"] \\
& \Theta_\ell^{-1}(\epsilon'), &
\end{tikzcd}
\]
where the slanted arrows are the natural maps and the top arrow is the bijection given by transport along the distinguished isomorphism $(X_\db)\rightarrow (X'_\db)$ (see Lemma \ref{lem: fibers of truncations between oDs are projections of fibers of truncations between Ds}(b)). Identifying $Aut((X_\db))$ and $Aut((X'_\db))$ via the distinguished isomorphism between them, the top map of this diagram is compatible with their actions (this is easily verified by \eqref{eq109}). We thus obtain a well-defined action of $\Gamma(\epsilon')$ on $\Theta_\ell^{-1}(\epsilon')$, which can be computed by taking any representative $(X_\db)$ of $\epsilon'$: via the canonical isomorphism $\Gamma(\epsilon')\rightarrow \Gamma((X_\db))$ and the left bijection of the diagram above, the action of $\Gamma(\epsilon')$ on $\Theta_\ell^{-1}(\epsilon')$ is the action of $Aut((X_\db))$ on $(\Theta_\ell^{-1}((X_\dbullet)))/ \sim'$.

We have all the necessary components to establish part (e) of Theorem \ref{thm: main thm 1}. Let us restate the result.
\begin{prop}\label{prop: main thm 1, part e}(Theorem \ref{thm: main thm 1}, part (e))
For every $\epsilon\in S_{\ell-1}(A)$ and $\epsilon'\in S'_{\ell-1}(A)$ above $\epsilon$, the natural map 
\[
\Theta_\ell^{-1}(\epsilon') \rightarrow \Theta_\ell^{-1}(\epsilon)
\]
(restricted from $S'_\ell(A)\rightarrow S_\ell(A)$) is surjective, and it descends to a bijection
\[
\Theta_\ell^{-1}(\epsilon')/\Gamma(\epsilon') \ \cong \ \Theta_\ell^{-1}(\epsilon).
\]
\end{prop}
\begin{proof}
Choose a generalized extension $(X_\db)$ in $\epsilon'$ (and hence $\epsilon$). We have a commutative diagram
\[
\begin{tikzcd}
\bigm(\Theta_\ell^{-1}((X_\dbullet))\bigm)\bigm/ \sim' \arrow[r, "\simeq"] \arrow[d, twoheadrightarrow] & \Theta_\ell^{-1}(\epsilon') \arrow[d]\\
\bigm(\Theta_\ell^{-1}((X_\dbullet))\bigm)\bigm/ \sim \arrow[r, "\simeq"] & \Theta_\ell^{-1}(\epsilon)
\end{tikzcd}
\]
where all maps are the natural ones and the horizontal maps are bijective by Lemma \ref{lem: fibers of truncations between oDs are projections of fibers of truncations between Ds}. Thus the right vertical arrow is also surjective. 

By Lemma \ref{lem: orbits of Gamma action}, two elements of $(\Theta_\ell^{-1}((X_\dbullet)))/ \sim'$ coincide in $(\Theta_\ell^{-1}((X_\dbullet)))/ \sim$ if and only if they are in the same orbit of the action of $Aut((X_\db))$ on $(\Theta_\ell^{-1}((X_\dbullet)))/ \sim'$. In view of the commutative diagram above and the definition of the action of $\Gamma(\epsilon')$ on $\Theta_\ell^{-1}(\epsilon')$, it follows that two elements of $\Theta_\ell^{-1}(\epsilon')$ coincide in $\Theta_\ell^{-1}(\epsilon)$ if and only if they are in the same $\Gamma(\epsilon')$-orbit.
\end{proof}

We end this discussion with a result about the stabilizers of the $\Gamma$-actions, which will be used in the totally nonsplit case studied in \S \ref{sec: tot nonsplit case} below.
\begin{lemma}\label{lem: stablizer of the action by automorphisms}
Let $(X_\db)\in D_{\ell-1}(A)$ and $(Y_\db)\in \Theta^{-1}((X_\db))$. Then the stabilizer of the $\sim'$-equivalence class of $(Y_\db)$ for the action of $Aut((X_\db))$ on \eqref{eq30} is the image of the restriction map 
\[Aut((Y_\db)) \hookrightarrow Aut((X_\db)).\]
\end{lemma}
\begin{proof}
Let $\sigma\in Aut((X_\db))$. Then $\sigma$ fixes the $\sim'$-equivalence class of $(Y_\db)$ if and only if there exists an isomorphism $f: \sigma\cdot(Y_\db)\rightarrow (Y_\db)$ that is identity on $A$, or equivalently, on $(X_\db)$ (by Lemma \ref{lem: if identity on A identity everywhere}). If $f$ is such an isomorphism, then let $\tilde{\sigma}$ be the composition $(Y_\db) \rightarrow \sigma\cdot(Y_\db)\xrightarrow{f} (Y_\db)$, where the first arrow is the isomorphism given by the transport data $\sigma$ (thus identity on the lowest diagonal and $\sigma$ on $(X_\db)$, see \S \ref{sec: fibers of truncations II}). Then $\tilde{\sigma}$ is an automorphism of $(Y_\db)$ that restricts to $\sigma$ on $(X_\db)$.

On the other hand, if $\sigma$ extends to an automorphism $\tilde\sigma$ of $(Y_\db)$, let $f$ be the composition $\sigma\cdot(Y_\db) \rightarrow (Y_\db) \xrightarrow{  \tilde\sigma  } (Y_\db)$, where the first arrow is the isomorphism given by the transport data $\sigma^{-1}$. Then $f$ is identity on $(X_\db)$.
\end{proof}

\begin{rem}
Let $(X_\db)\in D_{\ell-1}(A)$ and $(Y_\db)\in \Theta^{-1}((X_\db))$. By the previous lemma, the class of $(Y_\db)$ in $\Theta^{-1}((X_\db))/\sim'$ is a fixed point of the action of $Aut((X_\db))$ if and only if every automorphism of $(X_\db)$ extends to an automorphism of $(Y_\db)$. Using the formula of Lemma \ref{lem: compatibility of torsor structures} one can see that the former statement is also equivalent to linearity of the action of $Aut((X_\db))$ on $\prod_r Ext^1(A_{r+\ell}, A_r)$ obtained as follows: Use the class of $(Y_\db)$ as the base point to get an isomorphism
\[
\Theta^{-1}((X_\db))/\sim' \ \ \simeq \ \prod\limits_r Ext^1(A_{r+\ell}, A_r)
\]
(recall that the left hand side is a torsor over the right hand side). Now transport the action of $Aut((X_\db))$ from the left side to the right.

Already when $k=3$ it is easy to see that the map $Aut((Y_\db)) \hookrightarrow Aut((X_\db))$ does not always have to be surjective. Given $(X_\db)$ with nonempty $\Theta^{-1}((X_\db))$, it would be interesting to see if there always exists a $(Y_\db)\in \Theta^{-1}((X_\db))$ such that $Aut((Y_\db))\cong Aut((X_\db))$. Note that in the basic case of $k=3$ (i.e. the case of blended extensions), given $(X_\db)$ and arbitrary $(Y_\db)\in \Theta^{-1}((X_\db))$, the question of whether a given element of $Aut((X_\db))$ extends to an automorphism of $(Y_\db)$ has appeared previously in the work \cite{BVK16} of Barbieri-Viale and Kahn.\footnote{I thank Daniel Bertrand for bringing this to my attention.} Proposition D.1.5 of Appendix D therein gives a necessary and sufficient condition for this. (In the language of Barbieri-Viale and Kahn, the question is about whether a {\it partial gluing} extends to a {\it gluing}.)
\end{rem}

\section{The totally nonsplit case}\label{sec: tot nonsplit case}
In this section we consider an important special case of Theorem \ref{thm: main thm 1}(e) (or rather, its more explicit version stated as Proposition \ref{prop: main thm 1, part e}). We will be interested in the fibers of the truncation map $S_\ell(A)\rightarrow S_{\ell-1}(A)$ above {\it totally nonsplit} classes  of generalized extensions. This case of the result, which is recorded below as Proposition \ref{prop: weakly totally nonsplit case of Thm 1(e)}, will be used in the application to motives in \S \ref{sec: mixed motives with maximal unipotent radicals}.

From here forward, we shall assume that $\bT$ is a filtered {\it tannakian} category over a field of characteristic 0. We start by recalling some background on totally nonsplit extensions in \S \ref{sec: review of tot nonsplit exts} below. Then in \S \ref{sec: thm B, tot nonsplit case} we will define the analogous notion for generalized extensions and establish the desired result.

\subsection{Review of totally nonsplit extensions}\label{sec: review of tot nonsplit exts}
The following notion, which is due to Bertrand \cite{Ber01} as far as the author knows, will play a crucial role in the remainder of the paper.

\begin{defn}\label{def: tot nonsplit exts}
Let $X$ and $Y$ be objects of $\bT$.

\noindent (a) An extension or an extension class $\mathcal{E}$ of $\mathbbm{1}$ by $X$ is called totally nonsplit if for every proper subobject $X'$ of $X$ the pushforward of $\mathcal{E}$ along the quotient map $X\rightarrow X/X'$ is nonsplit.

\noindent (b) An extension or an extension class $\mathcal{E}$ of $Y$ by $X$ is called totally nonsplit if the extension class of $\mathbbm{1}$ by $\inHom(Y,X)$ corresponding to $\mathcal{E}$ under the canonical isomorphism
\begin{equation}\label{eq48}
Ext^1(Y,X) \cong Ext^1(\mathbbm{1}, \inHom(Y,X))
\end{equation}
is totally nonsplit.
\end{defn}
We first make a remark about the special case of the definition about extensions of $\mathbbm{1}$ by $X$. In view of the long exact sequence obtained by applying the functor $Hom(\mathbbm{1},-)$ to the sequence
\[
\begin{tikzcd}
   0 \arrow[r] & X' \arrow[r, ] & X \arrow[r, ] &  X/X'  \arrow[r] & 0,
\end{tikzcd}
\]
an extension (or extension class) $\mathcal{E}$ of $\mathbbm{1}$ by $X$ is totally nonsplit if and only if for every proper subobject $X'$ of $X$, the extension class of $\mathcal{E}$ is not in the image of the pushforward map
\[
Ext^1(\mathbbm{1}, X') \rightarrow Ext^1(\mathbbm{1}, X).
\]
We also make a cautionary remark about the general case of the definition, when $\mathcal{E}$ is an extension of $Y$ by $X$. The reader should keep in mind that in this definition, for $\mathcal{E}$ to be totally nonsplit we need to first consider $\mathcal{E}$ as an extension of $\mathbbm{1}$ (note that here we make use of the tensor structure and rigidity). The notion would remain the same if we considered the extension as an extension by $\mathbbm{1}$, with the statements being dualized. 

The following property of totally nonsplit extensions, which follows immediately from \cite[Theorem 1.2]{Es23}, will be used in the next subsection. Note that the proof of said result in {\it loc. cit.} makes use of the assumption that $\bT$ is a tannakian category over a field of characteristic zero.
\begin{lemma}\label{lem: automorphisms of totally nonsplit extensions}
Suppose that $\bT$ is a filtered tannakian category over a field of characteristic zero. Let $X$ and $Y$ be two nonzero objects of $\bT$. Let $\mathcal{E}$ be a totally nonsplit extension of $Y$ by $X$, with its middle object denoted by $E$. Suppose that every weight of $X$ is less than every weight of $Y$. Then the only endomorphisms of $E$ are the scalar maps.
\end{lemma}

\subsection{Totally nonsplit generalized extensions}\label{sec: thm B, tot nonsplit case}
In what follows, our notation (in particular, the graded object $A$) and terminology are as in \S \ref{sec: objects in filtered tan cats with given gr}, with the added assumption that $\bT$ is a filtered tannakian category over a field $K$ of characteristic 0. The following definition generalizes the notion of total nonsplitting to generalized extensions.

\begin{defn}\label{def: tot nonsplit gen exts}
We say a generalized extension $(X_\db)$ of any level is totally nonsplit if the extensions $\sX^v_{m,n}$ and $\sX^h_{m,n}$ (see Notation \ref{notation: horizontal and vertical extensions}) are totally nonsplit for every pair $(m,n)$ in the eligible range.
\end{defn}

In level 1, where a generalized extension is the data of an extension $\sE_r$ of $A_{r+1}$ by $A_r$ for each $r$, the generalized extension $(\sE_r)$ is totally nonsplit if and only if all the $\sE_r$ are totally nonsplit. The notion of totally nonsplit generalized extensions clearly descends to isomorphism classes of generalized extensions. 

Proposition \ref{prop: main thm 1, part e} gives a description of the fibers of the truncation map $S_\ell(A)\rightarrow S_{\ell-1}(A)$. We will show in this section that this description greatly simplifies for the fibers above totally nonsplit elements of $S_{\ell-1}(A)$. The following consequence of Lemma \ref{lem: automorphisms of totally nonsplit extensions} will be key:

\begin{lemma}\label{lem: aut group of tot nonsplit gen exts}
Let $(X_\db)$ be a totally nonsplit generalized extension of positive level. Then every automorphism of $(X_\db)$ is a scalar map, i.e. $Aut((X_\db))\cong K^\times$.
\end{lemma}
\begin{proof}
Let $\ell-1$ be the level. By restricting to the action on the objects on the lowest diagonal we have an injection
\begin{equation}\label{eq38}
Aut((X_\db)) \hookrightarrow \prod\limits_r Aut (X_{r,r+\ell}).
\end{equation}
By Lemma \ref{lem: gluing morphisms between objects on the lowest diagonal}, the image of this map is the set of those elements $(\sigma_{r,r+\ell})$ (with $\sigma_{r,r+\ell}$ an automorphism of $X_{r,r+\ell}$) the entries of which are compatible on the diagonal just above the lowest, i.e. the set of elements $(\sigma_{r,r+\ell})$ such that for each $r$, the two automorphisms of $X_{r,r+\ell-1}$ induced by $\sigma_{r-1,r+\ell-1}$ and $\sigma_{r,r+\ell}$ respectively via the arrows $X_{r-1,r+\ell-1}\twoheadrightarrow X_{r,r+\ell-1}$ and $X_{r,r+\ell-1}\hookrightarrow X_{r,r+\ell}$ coincide (when they are both available). Thanks to Lemma \ref{lem: automorphisms of totally nonsplit extensions}, for every $r$ the total nonsplitting of the extension $\sX^v_{r,r+\ell}$ or $\sX^h_{r,r+\ell}$ implies that every automorphism of $X_{r,r+\ell}$ is a scalar map. Thus each factor $Aut (X_{r,r+\ell})$ of the codomain of \eqref{eq38} is $K^\times$. Compatibility on the diagonal just above the lowest implies that the image of \eqref{eq38} is the diagonal copy of $K^\times$.
\end{proof}

Combining Lemma \ref{lem: aut group of tot nonsplit gen exts} with the description of the stabilizers of the $\Gamma$-actions in Lemma \ref{lem: stablizer of the action by automorphisms}, we obtain the following special case of Proposition \ref{prop: main thm 1, part e}:

\begin{prop}\label{prop: weakly totally nonsplit case of Thm 1(e)}
Let $2\leq \ell\leq k-1$ and $\epsilon \in S_{\ell-1}(A)$. Let $\epsilon'\in S'_{\ell-1}(A)$ be a lift of $\epsilon$ to $S'_{\ell-1}(A)$. Suppose $\epsilon$ is totally nonsplit. Then the action of $\Gamma(\epsilon')$ on $\Theta_\ell^{-1}(\epsilon')$ is trivial. The map $S'_\ell(A)\rightarrow S_\ell(A)$ restricts to a bijection
\[
\Theta_\ell^{-1}(\epsilon') \cong \Theta_\ell^{-1}(\epsilon).
\]
In particular, the choice of $\epsilon'$ above $\epsilon$ makes $\Theta_\ell^{-1}(\epsilon)$ a torsor for $\prod\limits_{r}Ext^1(A_{r+\ell},A_r)$.
\end{prop}

\begin{proof}
We need to argue that the action of $\Gamma(\epsilon')$ on $\Theta_\ell^{-1}(\epsilon')$ is trivial. The rest of the statement then follows from Proposition \ref{prop: main thm 1, part e} and Theorem \ref{thm: main thm 1}(b).

Let $(X_\db)\in D_{\ell-1}(A)$ be a representative of $\epsilon'$. By Lemma \ref{lem: aut group of tot nonsplit gen exts}, the automorphism group of $(X_\db)$ consists only of the nonzero scalar maps. Thus for every $(Y_\db)$ in $\Theta_\ell^{-1}((X_\db))$, the restriction map $Aut((Y_\db)) \rightarrow Aut((X_\db))$ is surjective, so that by Lemma \ref{lem: stablizer of the action by automorphisms} the action of $Aut((X_\db))$ on $\Theta^{-1}((X_\db))\bigm/ \sim'$ fixes the class of $(Y_\db)$. That is, the action of $Aut((X_\db))$ on $\Theta_\ell^{-1}((X_\db))/\sim'$ and hence the action of $\Gamma(\epsilon')$ on $\Theta_\ell^{-1}(\epsilon')$ is trivial.
\end{proof}

\begin{rem}\label{rem: comparing torsor structures on tot nonsplit fibers of level of S_l->S_(l-1) for different choices of lifts}
Let $\epsilon\in S_{\ell-1}(A)$ be a totally nonsplit element. Then the torsor structures on the fiber $\Theta_\ell^{-1}(\epsilon)$ corresponding to different choices of lifts of $\epsilon$ to $S'_{\ell-1}(A)$ are related to each other in a canonical way, as follows. Let $\epsilon'_1$ and $\epsilon'_2$ be two elements of $S_{\ell-1}'(A)$ above $\epsilon$. Use the notation $\ast_1$ (resp. $\ast_2$) for the translation operation for the torsor structure on $\Theta_\ell^{-1}(\epsilon)$ corresponding to the choice of $\epsilon'_1$ (resp. $\epsilon'_2$) as the lift of $\epsilon$. Then there exists an automorphism $\phi$ of $\prod_r Ext^1(A_{r+\ell},A_r)$ such that for every $\tilde{\epsilon}\in \Theta_\ell^{-1}(\epsilon)$ and $\sE\in \prod_r Ext^1(A_{r+\ell},A_r)$, we have
\begin{equation}\label{eq102}
\sE \ast_1 \tilde{\epsilon} = \phi (\sE) \, \ast_2 \, \tilde{\epsilon}.
\end{equation}
Indeed, choosing representatives $(X^{(1)}_\db)$ and $(X^{(2)}_\db)$ in $D_{\ell-1}$ for $\epsilon'_1$ and $\epsilon'_2$, and an isomorphism $(f_\db): (X^{(1)}_\db)\rightarrow (X^{(2)}_\db)$, denote the restriction of $(f_\db)$ to $A$ by $f_A$. The desired map $\phi$ is the automorphism of $\prod_r Ext^1(A_{r+\ell},A_r)$ induced by $f_A$ by pullbacks and pushforwards (i.e. given by $\phi(\sE)=f_A\cdot \sE$ in our earlier notation). Note that since $\epsilon$ is totally nonsplit, by Lemma \ref{lem: aut group of tot nonsplit gen exts} the isomorphism $(f_\db)$ is unique up to scaling, so that $\phi$ will not depend on the choice of $(f_\db)$. One can use Proposition \ref{prop: weakly totally nonsplit case of Thm 1(e)} and the formula of Lemma \ref{lem: compatibility of torsor structures} to obtain \eqref{eq102}.
\end{rem}

\begin{rem}
Let us say a generalized extension $(X_\db)$ of $A$ of level $\ell-1$ is {\it almost} totally nonsplit if for every $0\leq r\leq k-\ell$, at least one of the extensions $\sX^v_{r,r+\ell}$ or $\sX^h_{r,r+\ell}$ is totally nonsplit (note that $X_{r,r+\ell}$ is on the lowest diagonal of $(X_\db)$). Everything said in \S \ref{sec: thm B, tot nonsplit case} (but not in \S \ref{sec: classification of mixed motives with maximal unipotent radicals in graded independent case}) about totally nonsplit generalized extensions remains valid for almost totally nonsplit generalized extensions, with the exact same proofs. 
\end{rem}

\section{Motives with maximal unipotent radicals}\label{sec: mixed motives with maximal unipotent radicals}
\subsection{Setting and background}\label{sec: setting for section on application to motives}
In this final section we assume that $\bT$ is a filtered tannakian category over a field $K$ of characteristic zero such that the pure objects of $\bT$ are semisimple. The prototype examples of this are the category of graded-polarizable mixed Hodge structures over $\QQ$, and any reasonable tannakian category of mixed motives over a subfield of $\CC$ (e.g. those of Ayoub \cite{Ay14} and Nori \cite{HM17}, those defined via systems of realizations by Deligne \cite{De89} and Jannsen \cite{Ja90}, and Levine's category of mixed Tate motives over a number field \cite{Le93}). Inspired by the latter set of examples, we will often refer to an object of $\bT$ as a motive.

Let $X$ be an object of $\bT$. Let $\ffu(X)$ be the object of $\bT$ associated with the Lie algebra of the unipotent radical of the tannakian group of $X$. The object $\ffu(X)$ has been studied in various contexts by many, including Deligne (\cite{De89} and \cite[Appendix]{Jo14}), Andr\'{e} \cite{An92}, Bertrand \cite{Ber01}, Bertolin (\cite{Be02} and \cite{Be03}), Hardouin (\cite{Har06} and \cite{Har11}), Jossen \cite{Jo14}, and the author and Murty (\cite{EM1} and \cite{EM2}). We will take the background material on the definition of $\ffu(X)$ for granted, referring the reader to \S 4.2 and \S 2 of \cite{EM2} for a detailed review of this background. We just recall here that $\ffu(X)$ is the canonical subobject of $W_{-1}\inEnd(X)$ (where $\inEnd(X)$ is the internal Hom $\inHom(X,X)$) with the following property: for every fiber functor $\omega$ from $\bT$ to the category of vector space over $K$, if we consider the tannakian group of $X$ with respect to $\omega$ as a subgroup of $GL(\omega X)$, then
\[
\omega \ffu(X) \subset \omega W_{-1}\inEnd(X) = W_{-1}End(\omega X)
\]
is the Lie algebra of the kernel of the natural surjection from the tannakian group of $X$ with respect to $\omega$ to the tannakian group of $Gr^WX$ with respect to $\omega$. The fact that this kernel is the unipotent radical of the tannakian group of $X$ is because $Gr^WX$ is semisimple.

\begin{defn}\label{def: maximal unipotent radical}
We say $\ffu(X)$ is maximal\footnote{In \cite{EM2} instead of ``maximal" we used the word ``large" for this.} or $X$ has a maximal unipotent radical if
\[
\ffu(X) = W_{-1}\inEnd(X). 
\]
\end{defn}

The aim of this section is to use the earlier results of the paper to study isomorphism classes of motives with maximal unipotent radicals and associated graded isomorphic to a given semisimple object $A$, when $A$ is {\it graded-independent}. We start by reviewing some generalities in \S \ref{sec: maximality and total nonsplitting}. There is little new information in this discussion, and its inclusion is for referencing purposes and to put the later results in a better context. In \S \ref{sec: maximality criterion} we define the notion of graded-independence (Definition \ref{defn: graded independence}), and give a simple criterion for maximality of unipotent radicals in the case of graded-independent motives (Theorem \ref{thm: maximality criteria}). We then combine this in \S \ref{sec: classification of mixed motives with maximal unipotent radicals in graded independent case} with the earlier results of the paper to give a classification result for motives with maximal unipotent radicals and a given graded-independent associated graded (Theorem \ref{thm: classification of motives with give gr in graded independent case}). Finally, in \S \ref{sec: examples} we consider the example of 4-dimensional mixed Tate motives over $\QQ$ to illustrate the results.

\subsection{Maximality and total nonsplitting}\label{sec: maximality and total nonsplitting}
The following lemma summarizes some aspects of the relationship between maximality of unipotent radicals (Definition \ref{def: maximal unipotent radical}) and total nonsplitting of extensions (Definition \ref{def: tot nonsplit exts}) coming from the weight filtration. Recall that throughout, $\bT$ is a filtered tannakian category over a field of characteristic zero in which pure objects are semisimple. The objects of $\bT$ are referred to as motives.

\begin{lemma}\label{lem: maximality and total nonsplitting}
Let $X$ be a motive with more than one weight.

\noindent (a) If $\ffu(X)$ is maximal, then for every integers $\ell<m<n$ the extension
\begin{equation}\label{eq50}
\begin{tikzcd}
0 \arrow[r] & W_mX/W_\ell X \arrow[r] & W_nX/W_\ell X \arrow[r] & W_nX/W_m X \arrow[r] & 0
\end{tikzcd}
\end{equation}
is totally nonsplit.

\noindent (b) If $Gr^WX=Gr^W_mX \oplus Gr^W_nX$ with $m<n$, then $\ffu(X)$ is maximal if and only if the extension
\begin{equation}\label{eq49}
\begin{tikzcd}
0 \arrow[r] & Gr^W_mX \arrow[r] & X \arrow[r] & Gr^W_nX \arrow[r] & 0
\end{tikzcd}
\end{equation}
is totally nonsplit. (Thus when $X$ has only two weights, the necessary condition of part (a) for maximality of $\ffu(X)$ is also sufficient.)

\noindent (c) Taking $\bT$ to be the category of graded-polarizable rational mixed Hodge structures, there exists an object $X$ of $\bT$ with 3 weights for which all of the extensions of part (a) are totally nonsplit, but $\ffu(X)$ is not maximal. (Thus in general, the condition of part (a) for maximality is not sufficient.)
\end{lemma}
 
We first discuss parts (b) and (c). The former is by \cite[Corollary 3.4.1]{EM1}, which is a slightly more general version of a result of Hardouin \cite[Theorem 2]{Har11} (see also Theorem 2.1 of the unpublished article \cite{Har06}). In the setting of part (b) with two weights $m<n$ (and semisimple $Gr^WX$), \cite[Corollary 3.4.1]{EM1} asserts that $\ffu(X)$ is the smallest subobject of 
\[W_{-1}\inEnd(X)=\inHom(Gr^W_nX, Gr^W_mX)\] 
with the following property: the extension \eqref{eq49}, considered as an element of 
\[
Ext^1(\mathbbm{1}, \inHom(Gr^W_nX, Gr^W_mX)) 
\]
via the isomorphism \eqref{eq48} of Definition \ref{def: tot nonsplit exts}, splits after pushing forward along the quotient map 
\[
\inHom(Gr^W_nX, Gr^W_mX) \rightarrow \inHom(Gr^W_nX, Gr^W_mX)/\ffu(X).
\]
On recalling the definition of total nonsplitting (see Definition \ref{def: tot nonsplit exts}), the statement in (b) is now immediate.

As for part (c), we refer the reader to \cite[\S 5.2]{Es23} for an easy example with $X$ being the Hodge realization of a (Deligne \cite{De74}) 1-motive. One can further give an example where $X$ is the Hodge realization of a deficient 1-motive (deficient in the sense of \cite{Be02}, i.e. $W_{-2}\ \ffu(X)=0$). See \cite[\S 6.10]{EM2}.
\medskip\par 
Turning our attention to part (a), let us first recall a result from \cite{EM2}.
\begin{thm}[Theorem 4.9.1 of \cite{EM2}]\label{thm: characterization of u_p from the ANT paper}
Let $X$ be an object of $\bT$. Let $\sE_m(X)$ be the class of the extension
\begin{equation}\label{eq67}
\begin{tikzcd}
0 \arrow[r] & W_mX \arrow[r] & X \arrow[r] & X/W_mX \arrow[r] & 0
\end{tikzcd}
\end{equation}
considered as an element of 
\[
Ext^1(\mathbbm{1}, \inHom(X/W_mX, W_mX))
\]
via the canonical isomorphism \eqref{eq48}. Let 
\[\ffu_m(X) : = \ffu(X)\cap \inHom(X/W_mX, W_mX),\]
where $\inHom(X/W_mX, W_mX)$ is considered as a subobject of $W_{-1}\inEnd(X)$ in the natural way. Then for every subobject $L$ of $\inHom(X/W_mX, W_mX)$ one has 
\[
\ffu_m(X) \subset L
\] 
if and only if the pushforward $\sE_m(X)/L$ of $\sE_m(X)$ along the quotient map
\[
\inHom(X/W_mX, W_mX) \rightarrow \inHom(X/W_mX, W_mX)/L
\] 
is in the image of the obvious map
\[
Ext^1_{\langle (X/W_mX)\oplus W_mX \rangle}\bigm(\mathbbm{1}, \inHom(X/W_mX, W_mX)/L\bigm) \ \hookrightarrow \  Ext^1\bigm(\mathbbm{1}, \inHom(X/W_mX, W_mX)/L\bigm).
\]
\end{thm}
Here, the notation $\langle Y\rangle$ means the tannakian subcategory generated by an object $Y$ of $\bT$ (i.e., the smallest full tannakian subcategory of $\bT$ containing $Y$ and closed under subquotients). The notation $Ext^1_{\langle Y\rangle}$ means the $Ext^1$ group for the category $\langle Y\rangle$.

Referring to the setting of the theorem, note in particular that the result implies that if a subobject $L$ of $\inHom(X/W_mX, W_mX)$ has the property that the pushforward $\sE_m(X)/L$ splits, then $L$ contains $\ffu_m(X)$. Thus one obtains the following corollary:
\begin{cor}\label{cor: if u_p is maximal then E_p is totally nonsplit}
With the notation and setting as in Theorem \ref{thm: characterization of u_p from the ANT paper}, if 
\[\ffu_m(X)=\inHom(X/W_mX, W_mX),\] 
then $\sE_m(X)$ (or equivalently, the extension \eqref{eq67}) is totally nonsplit.
\end{cor}
We are ready to deduce part (a) of Lemma \ref{lem: maximality and total nonsplitting}.
\begin{proof}[Proof of Lemma \ref{lem: maximality and total nonsplitting}(a)]
Let $\ffu(X)=W_{-1}\inEnd(X)$. The inclusion $\langle W_nX/W_\ell X\rangle \subset \langle X\rangle$ induces a surjection from the tannakian group of $X$ to the tannakian group of $W_nX/W_\ell X$, which in turn induces a surjection $\ffu(X)\rightarrow \ffu(W_nX/W_\ell X)$. This map fits in a commutative diagram
\[
\begin{tikzcd}
\ffu(X) \ar[equal]{d} \arrow[r, twoheadrightarrow] & \ffu(W_nX/W_\ell X)  \arrow[d, hookrightarrow ]\\
W_{-1}\inEnd(X) \arrow[r, twoheadrightarrow ] & W_{-1}\inEnd(W_nX/W_\ell X),
\end{tikzcd}
\]
where the vertical arrows are the canonical inclusions and the bottom horizontal arrow is the map that after applying a fiber functor $\omega$, it sends an element of $W_{-1}End(\omega X)$ to its induced endomorphism of $W_n \omega X/W_\ell \omega X$. It follows that $\ffu(W_nX/W_\ell X)$ is also maximal. The assertion now follows from the previous corollary (applied to $W_nX/W_\ell X$ instead of $X$).
\end{proof}
\begin{rem}
The hypothesis of semisimplicity of pure objects of $\bT$ is not needed for Theorem \ref{thm: characterization of u_p from the ANT paper} (and is not assumed in Theorem 4.9.1 of \cite{EM2}). Subsequently, Corollary \ref{cor: if u_p is maximal then E_p is totally nonsplit} and Lemma \ref{lem: maximality and total nonsplitting}(a) are true in arbitrary filtered tannakian categories over fields of characteristic zero. However, Lemma \ref{lem: maximality and total nonsplitting}(b) needs the hypothesis of semisimplicity of pure objects.
\end{rem}

\subsection{Maximality of unipotent radicals for graded-independent motives}\label{sec: maximality criterion}
We now define the notion of graded-independence.

\begin{defn}\label{defn: graded independence}
Let $X$ be a nonzero motive. Denote the weights of $X$ by $p_1<\cdots<p_k$. Consider the $k$ objects
\begin{equation}\label{eq110}
C_r := \inHom(Gr^W_{p_{r+1}}X, Gr^W_{p_{r}}X) \hspace{.3in} (1\leq r\leq k-1)
\end{equation}
and
\begin{equation}\label{eq111}
C_0 := \bigoplus_{j-i>1} \inHom(Gr^W_{p_j}X, Gr^W_{p_{i}}X).
\end{equation}
We say that the motive $X$ is graded-independent if for every distinct $r,r'$ ($0\leq r,r'\leq k-1$), the two objects $C_r$ and $C_{r'}$ do not have any nonzero isomorphic subquotients (or equivalently, subobjects since $Gr^WX$ is semisimple).
\end{defn}

The property is an ``independence axiom"  in the spirit of such axioms in \cite{EM2}. In general, the reason such properties are of interest is twofold. Firstly, they force $Gr^W\ffu(X)$ to decompose in a way that makes it easier to study. Secondly, they are not far too restrictive, in the sense that they are satisfied in some very interesting situations. For instance, the property above (as well as all the independence axioms in \cite{EM2}) is satisfied as long as the weights of $X$ are sufficiently spread out so that the numbers $p_i-p_j$ are all distinct as the integers $i,j$ vary in $1\leq i<j\leq k$.
\medskip\par 
Recall that for an arbitrary motive $X$, maximality of $\ffu(X)$ is a stronger condition than the total nonsplitting of the canonical extensions coming from the weight filtration on $X$ (see Lemma \ref{lem: maximality and total nonsplitting}(c)). But this is not the case for graded-independent motives:

\begin{thm}\label{thm: maximality criteria}
Let $X$ be a graded-independent motive. Let $p_1<\cdots<p_k$ be the weights of $X$. Set $p_0:=p_1-1$ (so that $W_{p_0}(X)=0$). Then the following two statements are equivalent:
\begin{itemize}
\item[(i)] $\underline{\fu}(X)$ is maximal.
\item[(ii)] For every integer $1\leq r \leq k-1$, the extension
\begin{equation}\label{eq7}
\begin{tikzcd}
   0 \arrow[r] & Gr^W_{p_{r}}X \arrow[r, ] & \displaystyle{\frac{W_{p_{r+1}}X}{W_{p_{r-1}}X}} \arrow[r, ] &  Gr^W_{p_{r+1}}X \arrow[r] & 0
\end{tikzcd}
\end{equation}
is totally nonsplit.
\end{itemize}
\end{thm}
\begin{proof}
The fact that (i) implies (ii) is clear and in fact, does not require the graded-independence condition; see Lemma \ref{lem: maximality and total nonsplitting}(a). We will prove that (ii) implies (i). 

It is enough to show that 
\[
Gr^W\ffu(X) = Gr^WW_{-1}\inEnd(X).
\]
The graded-independence condition (on recalling that $\langle Gr^WX\rangle $ is semisimple) guarantees that the subobject $Gr^W\ffu(X)$ of
\[
Gr^WW_{-1}\inEnd(X) = \bigoplus\limits_{i<j}\inHom(Gr^W_jX,Gr^W_iX)
\]
decomposes according to the decomposition of $Gr^WW_{-1}\inEnd(X)$ as the direct sum of the $k$ objects \eqref{eq110} and \eqref{eq111} of Definition \ref{defn: graded independence}. For each $1\leq r \leq k-1$, consider the composition 
\[
\ffu(X) \twoheadrightarrow \ffu(W_{p_{r+1}}X/W_{p_{r-1}}X) \hookrightarrow \inHom(Gr^W_{p_{r+1}}X,Gr^W_{p_{r}}X),
\]
where the first map is induced by the inclusion of $\langle W_{p_{r+1}}X/W_{p_{r-1}}X\rangle$ in $\langle X\rangle$  and the second map is the natural inclusion. Since \eqref{eq7} is totally nonsplit, by Lemma \ref{lem: maximality and total nonsplitting}(b) this second arrow is an equality, so that the composition is surjective. Applying $Gr^W$ we get a surjection of $Gr^W\ffu(X)$ onto $\inHom(Gr^W_{p_{r+1}}X,Gr^W_{p_{r}}X)$. In view of the decomposition of $Gr^W\ffu(X)$ as the direct sum of its intersections with the $k$ objects of Definition \ref{defn: graded independence} and the fact that these $k$ objects do not have any nonzero isomorphic subquotients, it follows that $Gr^W\ffu(X)$ contains
\[
\inHom(Gr^W_{p_{r+1}}X,Gr^W_{p_{r}}X)
\]
for all $1\leq r\leq k-1$. Since $Gr^WW_{-1}\inEnd(X)$ is generated as a Lie algebra by these $k-1$ objects, the result follows.
\end{proof}

\subsection{Motives with maximal unipotent radicals and a prescribed graded-independent associated graded}\label{sec: classification of mixed motives with maximal unipotent radicals in graded independent case}
Fix $k\geq 2$ and nonzero pure motives $A_1,\ldots, A_k$, with $A_r$ of weight $p_r$ and $p_1<\cdots<p_k$. Set
\[
A : =\bigoplus\limits_{1\leq r\leq k} A_r.
\]
We will use the notation and language of \S \ref{sec: objects in filtered tan cats with given gr}. Thus $S(A)$ denotes the set of isomorphism classes of motives whose associated graded is isomorphic to $A$, and for each $1\leq \ell\leq k-1$ by $S_\ell(A)$ we denote the set of $\sim$-equivalence classes (i.e. isomorphism classes) of generalized extensions of level $\ell$ of $A$. The truncation map from level $\ell$ to level $\ell-1$ is denoted by $\Theta_\ell$. 

\begin{notation} Let $S^\ast(A)$ be the subset of $S(A)$ consisting of the isomorphism classes of motives $X$ with maximal $\ffu(X)$ and $Gr^W X$ isomorphic to $A$.
\end{notation}

In this subsection we will combine Theorem \ref{thm: maximality criteria} with our work in \S \ref{sec: objects in filtered tan cats with given gr} and \S \ref{sec: thm B, tot nonsplit case} to give a characterization of the set $S^\ast(A)$ when $A$ is graded-independent.
\medskip\par 
Recall from Definition \ref{def: tot nonsplit gen exts} that we say a generalized extension $(X_\db)$ of $A$ of any level is totally nonsplit if the extensions $\sX^v_{m,n}$ and $\sX^h_{m,n}$ (see Notation \ref{notation: horizontal and vertical extensions}) are totally nonsplit for every pair $(m,n)$ in the eligible range. This notion descends to isomorphism classes of generalized extensions. It is clear from the definition that if a generalized extension is totally nonsplit, then so are its truncations. In the graded-independent case, thanks to Theorem \ref{thm: maximality criteria} the converse is also true:

\begin{lemma}\label{lem: when A is graded-independent, a gen ext is totally nonsplit if and only if its first level is}
Let $A$ be graded-independent and $(X_\db)$ a generalized extension of $A$ of level $\ell\geq 2$. If the truncation of $(X_\db)$ to level one (i.e. $(X_{i,j})_{j-i\leq 2}$) is totally nonsplit, then so is $(X_\db)$.
\end{lemma}

\begin{proof}
Suppose $(X_{i,j})_{j-i\leq 2}$ is totally nonsplit. By definition, this means every extension 
\begin{equation}\label{eq41}
\begin{tikzcd}
0 \arrow[r] & A_i \arrow[r] & X_{i-1, i+1} \arrow[r] & A_{i+1}\arrow[r] & 0
\end{tikzcd}  
\end{equation}
is totally nonsplit. Consider an object $X_{r-1,r+\ell}$ on the lowest diagonal of $(X_\db)$. By Lemma \ref{lem: weight filtration of objects in generalized extensions} we have a canonical isomorphism from $Gr^W X_{r-1,r+\ell}$ to the direct sum of $A_r$, $A_{r+1}$, $\ldots$, $A_{r+\ell}$. Denoting this canonical isomorphism by $\phi$, from Lemma \ref{lem: D_k}(b) we know that the two generalized extensions (i) the part of $(X_\db)$ to the above and right of $X_{r-1,r+\ell}$ (i.e. consisting of the $X_{i,j}$ with $i\geq r-1$ and $j\leq r+\ell$) and (ii) $ext(X_{r-1,r+\ell}, \phi)$ (i.e. the generalized extension associated with $(X_{r-1,r+\ell}, \phi)$, see \S \ref{sec: gen exts, defn}) are $\sim'$-equivalent. It follows that after identifying
\[Gr^WX_{r-1,r+\ell} \cong \bigoplus\limits_{r\leq i\leq r+\ell} A_i\]
via $\phi$, for each $i$ with $r\leq i<r+\ell$ the class of the extension
\[
\begin{tikzcd}
0 \arrow[r] & Gr^W_{p_i}X_{r-1,r+\ell} \arrow[r] & W_{p_{i+1}}X_{r-1,r+\ell}/W_{p_{i-1}}X_{r-1,r+\ell} \arrow[r] & Gr^W_{p_{i+1}}X_{r-1,r+\ell} \arrow[r] & 0
\end{tikzcd}  
\]
coming from the weight filtration on $X_{r-1,r+\ell}$ is equal to the class of the extension \eqref{eq41}. In particular, the former extension is totally nonsplit for all $r\leq i<r+\ell$. Since $A$ is graded-independent, so is $X_{r-1,r+\ell}$, so that by Theorem \ref{thm: maximality criteria}, $X_{r-1,r+\ell}$ has a maximal unipotent radical. It now follows from Lemma \ref{lem: maximality and total nonsplitting}(a) that all of the extensions coming from the weight filtration on $X_{r-1,r+\ell}$ are totally nonsplit, so that $ext(X_{r-1,r+\ell}, \phi)$ is totally nonsplit. Being isomorphic to $ext(X_{r-1,r+\ell}, \phi)$, the part of $(X_\db)$ to the above and right of $X_{r-1,r+\ell}$ is thus also totally nonsplit. This is true for all $r$, hence $(X_\db)$ is totally nonsplit.
\end{proof}

We are ready to give the main result of this section.

\begin{thm}\label{thm: classification of motives with give gr in graded independent case}
For each $1\leq \ell\leq k-1$, denote the set of all totally nonsplit elements of $S_\ell(A)$ by $S^\ast_\ell(A)$. Recall that $S^\ast(A)$ is the subset of $S(A)$ consisting of the isomorphism classes of objects with maximal $\ffu$.

\noindent (a) We have a succession of maps
\[
S^\ast_{k-1}(A) \xrightarrow{\Theta_{k-1}}S^\ast_{k-2}(A) \xrightarrow{\Theta_{k-2}}S^\ast_{k-3}(A)\xrightarrow{\Theta_{k-3}} \cdots \xrightarrow{ \ \Theta_{3} \ } S^\ast_2(A) \xrightarrow{ \ \Theta_{2} \ }S^\ast_{1}(A)
\]
given by truncation. There is a natural bijection
\begin{equation}\label{eq68}
S^\ast_1(A) \cong \, \biggm\{(\sE_r)\in \prod\limits_{r}Ext^1(A_{r+1},A_r): \text{each $\sE_r$ is totally nonsplit}\biggm\} \biggm/ \hspace{-.05in} Aut(A),
\end{equation}
with the action of $Aut(A)$ given by pushforwards and pullbacks.

\noindent (b) Suppose that $A$ is graded-independent. Then the following statements are true:
\begin{itemize}
\item[(i)] There is a canonical bijection 
\[
S^\ast(A) \cong S^\ast_{k-1}(A).
\]
\item[(ii)] Let $\ell\geq 2$. The fiber of $\Theta_\ell: S^\ast_{\ell}(A)\rightarrow S^\ast_{\ell-1}(A)$ above any $\epsilon\in S^\ast_{\ell-1}(A)$ agrees with the fiber of $\Theta_\ell: S_{\ell}(A)\rightarrow S_{\ell-1}(A)$ above $\epsilon$. In particular, this fiber is a torsor over
\[
\prod\limits_{r} Ext^1(A_{r+\ell}, A_r)
\]
if it is nonempty.
\item[(iii)] For every $\ell\geq 2$, if the $Ext^2$ groups 
\[
Ext^2(A_{r+\ell}, A_r)
\]
vanish for all $1\leq r\leq k-\ell$, then $\Theta_\ell: S_\ell^\ast(A)\rightarrow S_{\ell-1}^\ast(A)$ is surjective.
\end{itemize}
\end{thm}

\begin{proof}
(a) The truncation of a totally nonsplit generalized extension is totally nonsplit, so the truncation maps between the $S_\ell(A)$ restrict to maps between the $S_\ell^\ast(A)$. The bijection \eqref{eq68} is the restriction of the bijection of Lemma \ref{lem: D_2}(b).
\medskip\par 
\noindent (b) The canonical bijection between $S^\ast(A)$ and $S^\ast_{k-1}(A)$ is the restriction of the bijection 
\begin{equation}\label{eq42}
S(A)\rightarrow S_{k-1}(A)
\end{equation}
of Lemma \ref{lem: D_k}(e), which is given by sending the isomorphism class of $X$ to the isomorphism class of the generalized extension $ext(X,\phi)$ associated with a pair $(X,\phi)$, where $\phi$ is any isomorphism $Gr^WX\rightarrow A$. The fact that \eqref{eq42} maps $S^\ast(A)$ into $S^\ast_{k-1}(A)$ is by Lemma \ref{lem: maximality and total nonsplitting}(a). The fact that the restricted map $S^\ast(A)\rightarrow S^\ast_{k-1}(A)$ is surjective follows from Theorem \ref{thm: maximality criteria} and the surjectivity of \eqref{eq42}: Given $\epsilon\in S^\ast_{k-1}(A)$, let $X$ be a motive whose isomorphism class gets mapped to $\epsilon$ by \eqref{eq42}. Then $X$ is graded-independent (because so is $A$). Since $\epsilon$ is totally nonsplit, so is $ext(X,\phi)$ for any isomorphism $\phi: Gr^WX\rightarrow A$. Thus in particular, all the extensions in Theorem \ref{thm: maximality criteria}(ii) are totally nonsplit. Theorem \ref{thm: maximality criteria} implies that $\ffu(X)$ is maximal. This finishes the proof of statement (i).

The first assertion in (b)(ii) is by Lemma \ref{lem: when A is graded-independent, a gen ext is totally nonsplit if and only if its first level is}. Indeed, given $\epsilon\in S^\ast_{\ell-1}(A)$ and $\tilde{\epsilon}\in S_{\ell}(A)$ above $\epsilon$, since $\epsilon$ is totally nonsplit, the truncation of $\tilde{\epsilon}$ to level 1 ( = the truncation of $\epsilon$) is totally nonsplit. Thus by Lemma \ref{lem: when A is graded-independent, a gen ext is totally nonsplit if and only if its first level is}, $\tilde{\epsilon}$ is totally nonsplit.

The second assertion in (ii) and the assertion in (iii) now follow from Proposition \ref{prop: weakly totally nonsplit case of Thm 1(e)} and Theorem \ref{thm: main thm 1}(c).
\end{proof}
\begin{rem} The special case of Theorem \ref{thm: classification of motives with give gr in graded independent case} when $k=3$, $A_3=\mathbbm{1}$, and $Ext^1(\mathbbm{1}, A_1)=0$ was proved in \cite[\S 6]{EM2} (see in particular, \S 6.7 therein).
\end{rem}

\subsection{Example: Mixed Tate motives with four weights and maximal unipotent radicals}\label{sec: examples}
In this final subsection, as an example that illustrates Theorem \ref{thm: classification of motives with give gr in graded independent case} in a case with $k>3$, we take $\bT$ to be the category $\mathbf{MT}(\QQ)$ of mixed Tate motives over $\QQ$ defined by Levine \cite{Le93}, and consider the problem of classifying isomorphism classes of 4-dimensional graded-independent objects of $\mathbf{MT}(\QQ)$ with four weights and maximal unipotent radicals.

Recall the description of Ext groups in $\mathbf{MT}(\QQ)$ (see \cite{DG05}, for instance): The $Ext^2$ groups all vanish, and
\begin{align}
\dim_\QQ Ext^1(\mathbbm{1}, \QQ(n)) =& \ \begin{cases} 1 \hspace{.3in} \text{if $n$ is odd and $> 1$}\\
0 \hspace{.3in} \text{if $n$ is even or $\leq 0$} \end{cases} \notag\\
Ext^1(\mathbbm{1}, \QQ(1)) \cong& \ \ \QQ^\times\otimes \QQ. \label{eq44}
\end{align}
For any odd integer $n>1$, the middle objects of all nonsplit extensions of $\mathbbm{1}$ by $\QQ(n)$ are isomorphic, as all nonzero elements of $Ext^1(\mathbbm{1}, \QQ(n))$ are in the same $\QQ^\times$-orbit. Denote this middle object, which is unique up to isomorphism, as $Z_n$ (one may think of this as ``the motive of $\zeta(n)$", because by Beilinson \cite{Bei84} it has a period matrix with entries $1$, $\zeta(n)/(2\pi i)^n$, $0$ and $1/(2\pi i)^n$). The extensions of $\mathbbm{1}$ by $\QQ(1)$ are given by Kummer motives. Under the isomorphism \eqref{eq44}, for any rational $r>1$ the extension class corresponding to $r\otimes 1$ arises from the weight filtration on the relative homology $L_r:=H_1(\mathbb{G}_m, \{1,r\})$ (this is ``the motive of $\log(r)$", with a period matrix with entries $1$, $\log(r)/(2\pi i)$, $0$ and $1/(2\pi i)$).

Set
\begin{equation}\label{eq45}
A := \QQ(a+b+c) \oplus \QQ(a+b) \oplus \QQ(a) \oplus \mathbbm{1},
\end{equation}
where $a,b$ and $c$ are positive integers. A motive with associated graded isomorphic to $A$ is graded-independent if and only if $a,b,c$ are distinct, $a+b\neq c$ and $b+c\neq a$. Assume said conditions hold. By Theorem \ref{thm: classification of motives with give gr in graded independent case}(b) we have truncation maps
\begin{equation}\label{eq46}
S^\ast(A) \cong S^\ast_3(A) \twoheadrightarrow S^\ast_2(A) \twoheadrightarrow S^\ast_1(A),
\end{equation}
which are surjective because the $Ext^2$ groups vanish. In view of Lemma \ref{lem: D_2} and the fact that the automorphism group of each $\QQ(n)$ is canonically isomorphic to $\QQ^\times$ we easily see
\begin{align}
S_1(A) & \cong Ext^1(\QQ(a+b), \QQ(a+b+c))/\QQ^\times \, \times \, Ext^1(\QQ(a), \QQ(a+b))/\QQ^\times \, \times \, Ext^1(\mathbbm{1}, \QQ(a))/\QQ^\times \notag\\
& \cong Ext^1(\mathbbm{1}, \QQ(c))/\QQ^\times \, \times \, Ext^1(\mathbbm{1}, \QQ(b))/\QQ^\times \, \times \, Ext^1(\mathbbm{1}, \QQ(a))/\QQ^\times.\label{eq70}
\end{align}
Taking the orbits of nonsplit ( = totally nonsplit) extensions in each factor above we get the subset $S_1^\ast(A)$ of $S_1(A)$ consisting of totally nonsplit elements. By the description of $Ext^1$ groups recalled earlier, the set $S_1^\ast(A)$ and hence $S^\ast(A)$ is nonempty if and only if $a$, $b$, $c$ are all odd (in which case the conditions $a+b\neq c$ and $b+c\neq a$ are automatic).

Assume $a,b,c$ are distinct odd positive integers. With a view towards implications for periods (see Remark \ref{rem: periods} below), the case when $1\in\{a,b,c\}$ is more interesting. Assume $b=1$ (the discussion for when $a$ or $c$ is 1 is similar). We start forming our generalized extensions from the smallest level, working backwards through the maps of \eqref{eq46}. From \eqref{eq70} and the description of the $Ext^1$ groups in $\mathbf{MT}(\QQ)$ we have
\[
S_1^\ast(A) \cong \{Z_c(a+1)\}\times  \{L_r(a) \}_r \times  \{Z_a\},
\]
where $r$ runs through a set of representatives for the nonzero orbits of the $\QQ^\times$-action on the vector space $\QQ^\times\otimes \QQ$. Once we fix $r$, the lifting to $S_2^\ast(A)$ is unique, as each fiber of $S^\ast_2(A) \rightarrow S^\ast_1(A)$ is a torsor over
\[
Ext^1(\QQ(a),\QQ(a+1+c)) \times Ext^1(\mathbbm{1},\QQ(a+1)) \ = \ 0.
\]
This unique lifting is the class of the generalized extensions of the form
\begin{equation}\label{eq47}
\begin{array}{cccc}
\QQ(a+1+c)&&&\\
Z_c(a+1) &\QQ(a+1)&&\\
M'_{c,r}(a) &L_r(a)&\QQ(a)&\\
&M_{a,r}&Z_a&\mathbbm{1}.
\end{array}
\end{equation}
Here, $M_{a,r}$ (resp. $M'_{c,r}$) is the object (unique up to isomorphism) that fits as the middle object of a blended extension of $Z_a$ by $L_r(a)$ (resp. $L_r$ by $Z_c(1)$). We refer the reader to \cite[\S 1.3 \& \S 6.8]{EM2} for a more detailed discussion of these motives. (The notation here for $M_{a,r}$ and $M'_{c,r}$ is consistent with {\it loc. cit.}, except for a slight difference in the indices.) Finally, the fiber of $S^\ast(A)\rightarrow S^\ast_2(A)$ above the class of \eqref{eq47} is a torsor over 
\[
Ext^1(\mathbbm{1}, \QQ(a+1+c))\simeq \QQ \, .
\]

\begin{rem}\label{rem: periods}
We end the paper with two remarks about periods. First, if $1\notin\{a,b,c\}$, then since the relevant Ext groups for $\mathbf{MT}(\QQ)$ coincide with the corresponding groups for the category $\mathbf{MT}(\ZZ)$ of mixed Tate motives over $\ZZ$ (see for instance, \cite{DG05}), all the motives in $S(A)$ will actually be in $\mathbf{MT}(\ZZ)$ and hence by Brown's work \cite{Br12} their periods will be generated over $\QQ$ by multiple zeta values and $2\pi i$. Second, back to the case $b=1$, every motive $X$ whose isomorphism class lives in the fiber of $S^\ast(A)$ above the class of \eqref{eq47} has a period matrix of the form
\[
\begin{pmatrix}
(2\pi i)^{-a-1-c} & (2\pi i)^{-a-1-c} \zeta(c) & (2\pi i)^{-a-1-c}p'_{c,r} & p_{a,r,c}(X)  \\ 
& (2\pi i)^{-a-1} & (2\pi i)^{-a-1} \log(r) & (2\pi i)^{-a-1}p_{a,r}\\
&& (2\pi i)^{-a} & (2\pi i)^{-a} \zeta(a)\\
&&&1
\end{pmatrix}.
\]
Here, the bottom right (resp. top left) $3\times 3$ matrix is a period matrix of $M_{a,r}$ (resp. $M'_{c,r}(a)$) and is independent of $X$. Since $\ffu(X)$ is maximal, the motivic Galois group of $X$ has dimension 7 ( = $1+\dim W_{-1}\inEnd(X)$). Assuming Grothendieck's period conjecture for $X$, the transcendence degree of the field generated by the periods of $X$ is 7, i.e. the numbers $p_{a,r,c}(X)$, $p_{a,r}$, $p'_{c,r}$, $\zeta(c)$, $\zeta(a)$, $\log(r)$, and $\pi$ are algebraically independent. Since $X$ is ramified at $r$ (because of $L_r$), as far as the author knows, aside from the cases of $r=2,3,4,6,8$ (where Deligne's \cite{De10} applies), the nature of the unspecified periods of $X$ is not known. See the Introduction of \cite{EM2} and the references therein for more details. 
\end{rem}

\end{document}